\documentclass[letterpaper,10pt]{article}

\usepackage{amsfonts,amssymb,eucal,amsmath,amsthm}
\usepackage{graphicx}
\usepackage[T1]{fontenc}
\usepackage{latexsym,url}
\usepackage{lscape}

\newtheorem{thm}{Theorem}[section] 
\newtheorem{lemma}[thm]{Lemma}     
\newtheorem*{numless_thm}{Theorem}

\newtheorem{cor}[thm]{Corollary}

\newtheorem{prop}[thm]{Proposition}

\newtheorem*{numless_lemma}{Lemma}

\theoremstyle{remark}
\newtheorem{rmk}[thm]{Remark}

\theoremstyle{definition}
\newtheorem{defn}[thm]{Definition}

\newcommand{\til}[1]{\widetilde{#1}}

\newcommand{\bdry}{\partial}

\title{Detection of incompressible surfaces in hyperbolic punctured torus bundles}
\author{Henry Segerman}

\begin{document}
\maketitle

\begin{abstract}
Culler and Shalen, and later Yoshida, give ways to construct incompressible surfaces in 3-manifolds from ideal points of the character and deformation varieties, respectively. We work in the case of hyperbolic punctured torus bundles, for which the incompressible surfaces were classified by Floyd and Hatcher. We convert non fiber incompressible surfaces from their form to the form output by Yoshida's construction, and run his construction backwards to give (for non semi-fibers, which we identify) the data needed to construct ideal points of the deformation variety corresponding to those surfaces via Yoshida's construction. We use a result of Tillmann to show that the same incompressible surfaces can be obtained from an ideal point of the character variety via the Culler-Shalen construction. In particular this shows that all boundary slopes of non fiber and non semi-fiber incompressible surfaces in hyperbolic punctured torus bundles are strongly detected.\\
\end{abstract}

\section{Introduction}
\label{intro}
It is well known that the shape of an ideal hyperbolic tetrahedron can be specified by a single complex variable: by conjugating we can move the tetrahedron around in the upper half space model of $\mathbb{H}^3$ so that three of the vertices are at 0, 1 and $\infty$. The fourth vertex will then be at some point $z$ on the complex plane boundary of the upper half space model. Looking down on the complex plane, we can also see $z$ as the complex dihedral angle between the faces either side of the $0$ - $\infty$ edge. If we rearrange the tetrahedron so that different vertices are moved to 0, 1 and $\infty$, we also get complex dihedral angles of $\frac{z-1}{z}$ and $\frac{1}{1-z}$.\\

When we glue ideal tetrahedra together to form a manifold $M$ (open at the ideal vertices), we require that the product of the complex dihedral angles around each edge is 1 (with the sum of the arguments of the angles being $2 \pi$). This gives us a \textbf{gluing equation} for each edge, in terms of various $z$, $\frac{z-1}{z}$ and $\frac{1}{1-z}$ factors, with different $z$ variables for each tetrahedron. Multiplying up by denominators will give us polynomial equations, and so we have an affine variety $\mathfrak{T}(M)$ in $(\mathbb{C}\setminus\{0,1\})^N$, where $N$ is the number of tetrahedra in the tetrahedralisation of $M$.\\

$\mathfrak{T}(M)$ depends on the tetrahedralisation of the manifold of course, but also on the choice of one of the three dihedral angles in each tetrahedron.

\begin{defn}
Let $\mathcal{T}$ be a tetrahedralisation of the 3-manifold $M$, with $N$ tetrahedra. Consider the same tetrahedralisation, labelled with a choice of one of the three dihedral angles in each tetrahedron. By an abuse of notation we will also refer to this as $\mathcal{T}$. We define the \textbf{tetrahedron variety} of $M$ with respect to the labelled tetrahedralisation $\mathcal{T}$, $\mathfrak{T}(M) = \mathfrak{T}(M;\mathcal{T})$ to be the affine variety in $(\mathbb{C}\setminus\{0,1\})^N$ defined as the solutions of the gluing equations, where each dimension of the ambient space corresponds to a tetrahedron.
\end{defn}
This is also known as Thurston's \textbf{parameter space}. A closely related variety retains the symmetry of the tetrahedralisation at the cost of using three times as many variables, by not making a choice of dihedral angle in each tetrahedron:
\begin{defn}
Let $\mathcal{T}$ be a tetrahedralisation of the 3-manifold $M$, with $N$ tetrahedra. 
We assign 3 complex variables to each tetrahedron, $z_1$, $z_2$ and $z_3$, one for each dihedral angle, related by the identities (See \cite{thurston}): 
$$z_1 z_2 z_3 = -1$$ 
$$1 - z_1 + z_1 z_2 = 0$$
Each gluing equation now states that a product of positive powers of $z_1$, $z_2$ and $z_3$ must be equal to 1, and so we define the \textbf{deformation variety} of $M$ with respect to the tetrahedralisation $\mathcal{T}$, $\mathfrak{D}(M) = \mathfrak{D}(M;\mathcal{T})$ to be the affine variety in $(\mathbb{C}\setminus\{0,1\})^{3N}$ defined as the solutions of the gluing equations together with the identities between complex dihedral angles within each tetrahedron, where each of the three complex dihedral angles in each tetrahedron corresponds to a dimension of the ambient space.
\end{defn}
We can view $\mathfrak{T}(M;\mathcal{T})$ as a projection of $\mathfrak{D}(M;\mathcal{T})$ onto a certain subset of the variables given by the labelling on the tetrahedralisation. We can of course recover $\mathfrak{D}(M)$ from a given $\mathfrak{T}(M)$ since $z_2 = \frac{z_1-1}{z_1}$ and $z_3 = \frac{1}{1-z_1}$. Thus we will mostly work with a particular $\mathfrak{T}(M)$, but consider $\mathfrak{D}(M)$ as the concrete underlying object.

\begin{defn}
\label{tetra_ideal_point}
An \textbf{ideal point of the tetrahedron variety} is a limit point of $\mathfrak{T}(M)$  at which one or more of the tetrahedra angles converges to 0, 1 or $\infty$. We say that such a tetrahedron \textbf{degenerates}. We retain information about the relative rates at which different tetrahedra degenerate.  \\

More formally: \\

An \textbf{ideal point of the deformation variety} is a limit point $\bar{p}$ on $S^{3N-1}$ of the set of $3N$-tuples in the interior of $B^{3N}$: 
$$\left\{ \left. \frac{\left(\log|z^{(1)}_1|, \log|z^{(1)}_2|, \ldots, \log|z^{(N)}_3|\right)}{\sqrt{1 +\sum\left(\log|z^{(j)}_i|\right)^2}} \right| \left(z^{(1)}_1, z^{(1)}_2, \ldots, z^{(N)}_3 \right) \in \mathfrak{D}(M)  \right\}$$  
\end{defn}

This is Bergman's logarithmic limit set, defined in \cite{bergman71}. We could also define a similar logarithmic limit set for some $\mathfrak{T}(M;\mathcal{T})$, although we would then lose the relative rate information for any degenerations which send a particular labelled complex dihedral angle to 1. \\

\begin{rmk}
Yoshida~\cite{yoshida91} (in his Definition 3.1) has a slightly different definition of an ideal point (of the tetrahedron variety), requiring that a certain slope on the boundary torus of the manifold is non trivial (and thus that the holonomy of any other boundary curve will blow up as the ideal point is approached). Such an ideal point of the deformation variety will necessarily give an incompressible surface via Yoshida's construction after compressions. Yoshida's construction still goes through without this extra condition and we obtain a surface, but it may compress away to the empty surface, or a boundary parallel surface. In all cases in this paper however, ideal points of the deformation variety do satisfy Yoshida's extra condition. We prefer this definition because we feel it is a more natural condition on the deformation variety itself. For example, it may be that some 3-manifolds contain closed incompressible surfaces which can be constructed from an ideal point of the deformation variety (in our sense). Such surfaces have no boundary slopes and the corresponding limit of the deformation variety would not be an ideal point in the Yoshida version of the definition.
\end{rmk}

\begin{defn}
A surface $S$ in a 3-manifold $M$ with $\partial{S} \subset \partial{M}$ is said to be \textbf{incompressible} if $S$ has no sphere or boundary parallel components and if every loop in $S$ that bounds a disk in $M \setminus S$ also bounds a disk in $S$. A surface with boundary is said to be \textbf{$\partial$-incompressible} if every arc $\alpha$ in $S$ (with $\partial(\alpha) \subset \partial S$) which is homotopic to $\partial M$ is homotopic in $S$ to $\partial S$. We will in general refer to surfaces that are incompressible, $\partial$-incompressible and not boundary parallel as incompressible surfaces, and deal in this paper only with oriented surfaces.
\end{defn}

Culler and Shalen~\cite{cullershalen83} 
give a method of constructing an incompressible surface in a 3-manifold from an ideal point of the character variety of that manifold.  The question naturally arises of the degree to which a reverse construction might be possible. That is, given an incompressible surface, does it come from an ideal point?

\begin{defn}
If $M$ is a 3-manifold with a single torus boundary, orientable, irreducible and compact, we say that a boundary slope of an incompressible surface that can be produced by the Culler-Shalen construction from an ideal point of the character variety is \textbf{detected}. If there is no closed surface that comes from the same ideal point then the slope is \textbf{strongly detected}, otherwise it is \textbf{weakly detected}. We also refer to an incompressible surface as being \textbf{detected} (by the character variety) if it can be produced by the Culler-Shalen construction from an ideal point of the character variety.
\end{defn}

Some previous results about detection of surfaces: Ohtsuki~\cite{ohtsuki94} gives a classification of surfaces in 2-bridge knot complements that are detected, and in particular shows that not every incompressible surface can be obtained from the construction. Schanuel and Zhang~\cite{schanuelzhang01} gave the first examples of non-fiber (and non semi-fiber) boundary slopes that are not strongly detected, although they are weakly detected. Chesebro and Tillmann~\cite{chesebrotillmann05} give an infinite family of hyperbolic knots, each of which has at least one boundary slope of an incompressible surface (non-fiber and non semi-fiber) that is not strongly detected.

\begin{thm}
\label{main_result}
All non fiber and non semi-fiber incompressible surfaces in hyperbolic punctured torus bundles over the circle are detected by the character variety.
\end{thm}
\begin{cor}
All boundary slopes of non fiber and non semi-fiber incompressible surfaces in punctured torus bundles are strongly detected.
\end{cor}
\begin{proof}
Due to the classification of the incompressible surfaces of (hyperbolic)  punctured torus bundles by Floyd and Hatcher~\cite{floydhatcher82} (and independently by Culler, Jaco and Rubinstein~\cite{cullerjacorubinstein}), we know that there are no closed incompressible surfaces, so all detected slopes are in fact strongly detected.
\end{proof}
\begin{rmk}
Punctured torus bundles with elliptic or parabolic monodromy should presumably be not too hard to analyse, but restrict ourselves to hyperbolic monodromies here.
\end{rmk} 

In 1982, Floyd and Hatcher~\cite{floydhatcher82}, and Culler, Jaco and Rubinstein~\cite{cullerjacorubinstein} classified the orientable incompressible surfaces in punctured torus bundles. We will work from the Floyd-Hatcher version. Yoshida~\cite{yoshida91} constructs an incompressible surface given an ideal point of the tetrahedra variety of non-zero (i.e. non-trivial) slope. \\

Given an incompressible surface in a punctured torus bundle, obtained from Yoshida's construction, it is not immediately obvious which Floyd-Hatcher surface it is isotopic to. In fact, as constructed, the Yoshida surface may need a number of ambient 2-surgeries and deletions of sphere components before it is incompressible (and boundary incompressible). However, it must of course be reducible to one of the Floyd-Hatcher surfaces.\\

The plan of attack is to reverse this process: to start with a Floyd-Hatcher surface, isotope and add sphere components until it is in the form of a Yoshida surface, and then use the data we obtain from the position of the Yoshida surface to construct the ideal point of the deformation variety, which we can then convert to an ideal point of the character variety using a result of Tillmann~\cite{tillmann_degenerations}.\\

In fact it is generally easier to follow the isotopies from the more convoluted Yoshida form to the simpler Floyd-Hatcher form, so our argument for that part of the proof proceeds in that direction: For a surface given to us in Floyd-Hatcher form, we give the corresponding surface in Yoshida form (in section \ref{translation}) and check that the Yoshida form surface simplifies to the given Floyd-Hatcher form (in section \ref{surfaces_isotope}). \\

Yoshida's construction relates the rates at which various tetrahedra degenerate to the position of the surface with respect to the tetrahedralisation (the number of twisted square pieces of the surface in each tetrahedron give the relative rates of degeneration). We use this relation in reverse: we need to show that the degeneration rates implied by the Yoshida form of our surface correspond to an ideal point $\bar{p}$ of the deformation variety. For surfaces that are not semi-fibers this is possible, and we identify which these are in section \ref{semi-fibers}.\\

In section \ref{tilde} we obtain a solution to some form of the gluing equations that correspond to the (hopefully) ideal point (there may be non-degenerate tetrahedra in the tetrahedralisation, as well as the ones that are degenerate), and show that we can approach this solution as the limit of finite (no angle is 0, $\infty$ or 1) points of the deformation variety.\\

The solution we will find for the ideal point will be explicit, in the sense of giving actual complex values for the angles in the non-degenerate tetrahedra, and the "directions" of degeneration for the degenerate tetrahedra. \\

The equations we find a solution to are not the gluing equations themselves, since a number of variables are supposedly converging to 0 or $\infty$ as we approach the ideal point. Instead we make a number of changes of variable. First we change which dihedral angle within each degenerating tetrahedron is labelled, to standardise so that the labelled angle is the one converging to 0. Secondly we perform a kind of weighted "blow up" (of the algebraic geometry sort), replacing each variable which is now converging to 0 with a "direction" variable, multiplied by some power of a global "convergence variable", which we call $\zeta$. So if $Z$ is the complex angle in some tetrahedron, which is supposedly converging to 0 at rate k, we replace $Z$ in our equations with $\zeta^k y$ (here $y$ is a direction variable). Each gluing equation describes the complex angles around an edge of the tetrahedralisation, and for the gluing equations to be satisfied as our tetrahedra degenerate, we must have the sum of the rates at which angles around this edge converge to 0 to equal the sum of the rates at which angles converge to $\infty$. This is the case for the rates we obtain from our Yoshida form of the surface. The effect this has on our equations is that a power of $\zeta$ factors out from each equation corresponding to an edge at which some dihedral angles are converging to 0 or $\infty$. Deleting this factor and rearranging the equations to form polynomial equations we reach a form of the equations we call \textbf{tilde equations}. \\

It is these equations we will find a solution for, when $\zeta = 0$, corresponding to being at the ideal point. We bring in another idea here, concerning the angle variables (that is, the ones that do not degenerate as we approach the ideal point):

\begin{numless_lemma}
$a_k = \frac{1 - \cos{k \beta}}{1 - \cos{\beta}}$  is a solution of  $a_{k-1}a_{k+1} - (1 - a_k)^2 = 0$
\end{numless_lemma}
This recursion relation comes from the gluing equations at the 4-valent edges throughout one of the two types of fan of a torus bundle, and holds whether or not we are at a point at which some tetrahedra are degenerating. A very similar result, for a closely related recursion relation, holds for the other type of fan. It turns out that when tetrahedra are degenerating, the ends of a fan degenerate in such a way that we can add an extra "fake" variable to each end of the fan, set their values to be 1, and then the recursion relations also hold at the ends of the fans.\\

Using these results we can assign values to variables corresponding to all of the non-degenerate tetrahedra. The values assigned will be real numbers, but not 0 or 1 (or $\infty$), so they correspond to flat but non-degenerate tetrahedra. The remaining direction variables depend on each other in understandable ways: the idea of the method is to set the limiting holonomy of the meridian to be $1$ (actually something closely related to this), and use this normalisation to determine all direction variables, with a choice of sign in some cases. None of the direction variables are assigned a value of 0.\\

Having obtained a solution $\tilde{p}$ of the tilde equations, we then want to show that we also have solutions nearby that correspond to finite points of the tetrahedron variety. There are two parts to this:
\begin{itemize}
\item (i) We show that there are points of the variety defined by the tilde equations arbitrarily near to $\tilde{p}$.
\item (ii) We show that points $q$ close enough to $\tilde{p}$ must have the variable $\zeta \neq 0$.
\end{itemize}
Then by continuity, at $q$ all angle variables are still away from 0 or 1, and all direction variables are away from 0, and so multiplying them by the appropriate power of $\zeta$ (recall that the original complex angle, $Z = \zeta^k y$) we obtain a complex angle that is near, but not equal to 0. Hence all complex angles of tetrahedra for this point are non-degenerate and we have a finite point of the tetrahedron variety.\\

We can see (i) as a consequence of the fact that we have one more variable than we have equations, and provide an algebraic geometry proof of this. (ii) comes from the fact that $\tilde{p}$ is isolated among solutions with $\zeta = 0$. We show this by considering the steps we took to find the solution $\tilde{p}$, and show that we only ever had finitely many choices at each step (there are a finite number of choices for each fan of angle variables and some finite number of choices of sign later on, so sufficiently close to $\tilde{p}$ we would have to choose the exact same parameters as for $\tilde{p}$ itself).\\

We obtain an ideal point of the deformation variety, and then use a result in \cite{segerman_twsq_spnn}, which uses Tillmann's result in \cite{tillmann_degenerations} to show that a corresponding ideal point of the character variety produces the same incompressible surface. \\

I would like to thank my thesis advisor, Steve Kerckhoff, whose guidance has been invaluable. I have also benefited from conversations about this with Daryl Cooper, David Futer, Fran\c{c}ois Gu\'eritaud, Stephan Tillmann, Ravi Vakil and Kirsten Wickelgren.

\section{The Canonical Tetrahedralisation of a Torus Bundle}
\sectionmark{The Canonical Tetrahedralisation}
\label{tetrahedralisation}

We use the tetrahedralisation $\mathcal{H}$, sometimes called the Floyd-Hatcher or monodromy tetrahedralisation. It first appears in Floyd-Hatcher~\cite{floydhatcher82}, based on an idea from \cite{thurston}. It is also known as J\o rgensen's ideal triangulation, and was proved to be the canonical triangulation of punctured torus bundles by Lackenby~\cite{lackenby03}. Gu\'{e}ritaud~\cite{gueritaud} gives an excellent exposition.\\

\begin{figure}[hp]
\begin{center}
\includegraphics[width=0.6\textwidth]{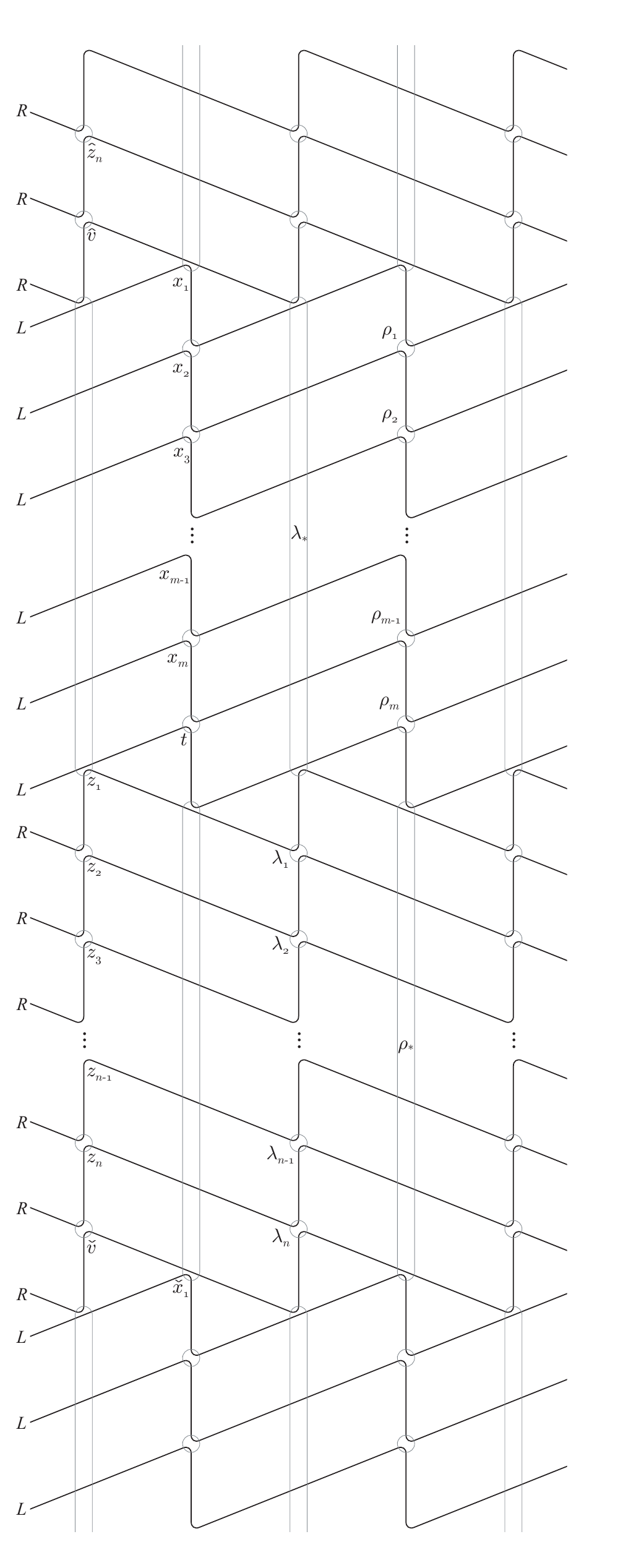}
\caption{Canonical tetrahedralisation of a punctured torus bundle. See section \ref{tilde} for information about the labelling.}
\label{tetrahedralisation_diag}
\end{center}
\end{figure}

Figure \ref{tetrahedralisation_diag} shows a picture of the tetrahedralisation as seen from the torus boundary. $\mathcal{H}$ consists of a stack of tetrahedra, one on top of the next. Each (ideal) tetrahedron has four vertices at infinity, and we have truncated each tetrahedron at each of its four vertices to produce four triangles on the boundary torus. We can see the four triangles in the layers labelled $t$ or $v$. The vertices of the resulting triangulation of the boundary torus are shown on the diagram with circles around them, labelled $\lambda_k$ or $\rho_k$. There are also special vertices, labelled $\lambda_*$ and $\rho_*$ which have been stretched out on this diagram for clarity. We are to imagine collapsing these "long" vertices down to points. Doing this will also change all of the apparently 4 sided polygons in the diagram into triangles, as expected for the truncated ends of a tetrahedron. The edges of those triangles do not quite meet at the vertices in order to highlight which tetrahedron a boundary triangle comes from. A layer of triangles which is "connected" through the vertices are all the truncated boundary of the same tetrahedron.\\

As discussed in section \ref{intro}, the shape of an ideal tetrahedron in $\mathbb{H}^3$ is specified by one of its dihedral angles, together with a scaling factor across that angle. This information is encoded as a single complex number ("complex angle") assigned to one of the dihedral angles in the tetrahedron (See \cite{thurston}). This shows up on the torus boundary as a complex angle at one of the three corners of each triangle. If we label one angle $z$, then moving clockwise around the triangle, the other two angles are $\frac{z - 1}{z}$ and $\frac{1}{1-z}$. We choose the uppermost dihedral angle for the labelling of each tetrahedron in our tetrahedralisation $\mathcal{T}$, where by "uppermost" we mean in relation to the (unlabelled) tetrahedralisation $\mathcal{H}$ of the torus bundle. It turns out that opposite edges of an ideal tetrahedron have the same complex angle, and so the value of the bottom-most dihedral angle is the same as that of the uppermost.\\

On the left side of the diagram we see the boundaries between tetrahedron "layers", labelled with either $L$ or $R$. These are the $L$ and $R$ from the decomposition of the monodromy $\phi$ into the generators. A tetrahedron that lies between an $L$ and an $R$ is called a \textbf{hinge} tetrahedron, and we use the variables $t$ and $v$ to describe the uppermost angle of those tetrahedra. All other tetrahedra are part of \textbf{fans} of tetrahedra, separated from neighbouring fans by the hinge tetrahedra (when the "long" vertices are collapsed, the torus boundary picture of such a sequence of tetrahedra looks like a fan). We use the variables $x_i$ and $z_j$ to refer to the uppermost angles of those tetrahedra. It is possible for there to be no tetrahedra within a fan when the hinge tetrahedra are next to each other. It has been observed that fans of tetrahedra seem to act very much as a unit, and one of the themes of this paper is to make explicit some aspects of this notion.  \\

\section{Various Forms of Surfaces}
\label{forms_of_surface}

As mentioned in the introduction we will be deforming surfaces that begin in Yoshida form into Floyd-Hatcher form. We will divide the torus bundle into sections, roughly at the boundaries between fans (this will be made precise later), and perform the necessary surface isotopies locally (i.e within each section individually). In the following two subsections we will describe the two forms of surfaces, and convert them to our own format, within which we will perform the isotopies.\\

\begin{figure}[hp]
\begin{center}
\includegraphics[width=0.5\textwidth]{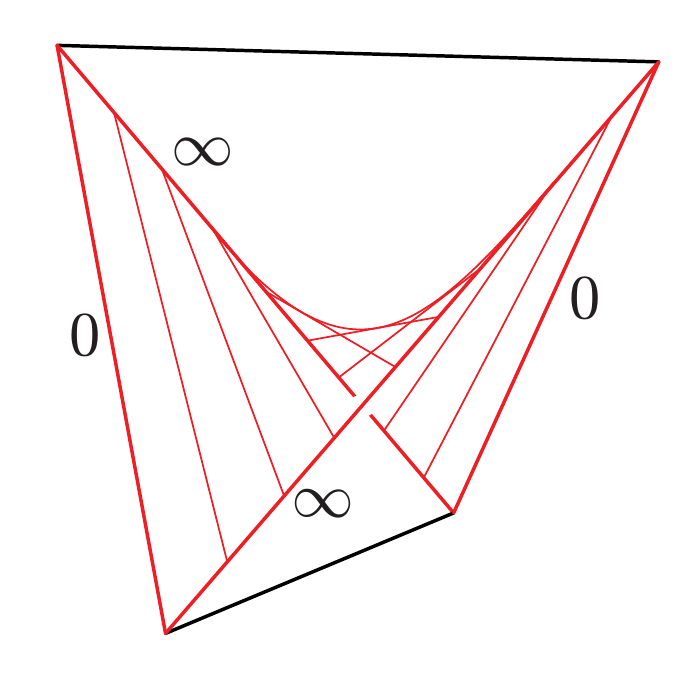}
\caption{A twisted square within a tetrahedron, with 0 and $\infty$ edges labelled (see section \ref{yoshida_form}). $\mathcal{N}_e$ neighbourhoods are not shown.}
\label{fig_twisted_square.pdf}
\end{center}
\end{figure}

Before describing this format we alter our tetrahedralisation of $M$ slightly to give a different subdivision of $M$ into pieces (but continue to refer to it as $\mathcal{H}$): First we truncate all the tetrahedra so that we get an induced triangulation on $\partial M$. We add a layer $\partial M \times I$ to the subdivision, so that $\partial M \times \{1\}$ touches the truncated ends of the tetrahedra and $\partial M \times \{0\}$ is $\partial M$.  We take small disjoint cylindrical neighbourhoods $\mathcal{N}_e$ of the edges $e$ of the tetrahedralisation.
For the tetrahedron with complex angle label $z$, let $\mathcal{T}_z$ be the truncated tetrahedron minus the union of the $\mathcal{N}_e$ for all $e$ that the tetrahedron is incident to. $M$ is then the disjoint union of $\partial M \times I$ with the $\mathcal{N}_e$ and $\mathcal{T}_z$. \\
\begin{defn}
Our surfaces will be made from five kinds of pieces, all of which intersect each other only at their boundaries:
\begin{itemize}
\item A twisted square that sits inside a $\mathcal{T}_z$, with its four edges on four of the six "edge neighbourhood" boundaries. See figure \ref{fig_twisted_square.pdf}.
\item A triangle that sits inside a $\mathcal{T}_z$, parallel and close to one of the 4 faces of the original tetrahedron, with its three edges on the edge neighbourhood boundaries that bound the face of the original tetrahedron.
\item A long thin strip that sits inside an $\mathcal{N}_e$ and respects its product structure.
\item An annulus that sits inside the $\partial M \times I$ and respects its product structure, with boundary curves on each boundary of $\partial M \times I$.
\item A disk that sits within the $\partial M \times I$ and whose boundary curve is on $\partial M \times \{1\}$, the boundary of $\partial M \times I$ that faces the interior of the manifold.
\end{itemize}
\label{piece_types}
\end{defn}

The thin strips serve to glue together the twisted squares and triangles near to edges of the tetrahedralisation. The first three types of surface piece have boundary on $\partial M \times \{1\}$ and so strictly speaking the twisted square is an octagon (it has an edge across $\partial M \times \{1\}$ at each "corner" of the twisted square) and the triangle is a hexagon (also has an edge at each "corner"). The strip has 4 edges: two long edges parallel to the $e$ of the $\mathcal{N}_e$ in which the strip lies, and two short edges on $\partial M \times \{1\}$. See Figures \ref{tetra_with_surfaces1} and \ref{tetra_with_surfaces2} for some pictures of these pieces of surface in a tetrahedron.

\subsection{Incompressible Surfaces in Floyd-Hatcher Form}
\label{fhform}
\begin{figure}[ht]
\begin{center}
\includegraphics[width=0.7\textwidth]{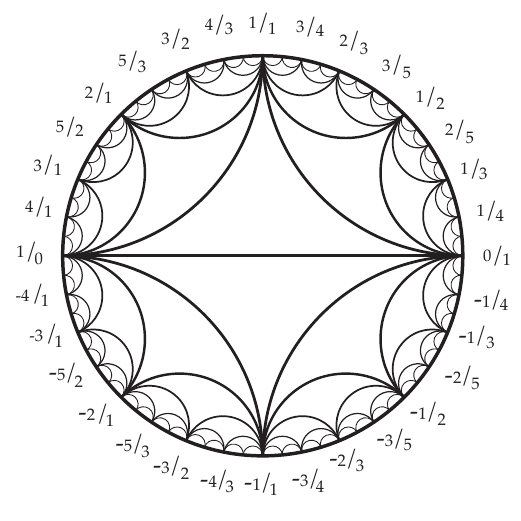}
\caption{The Farey Graph (diagram by Allen Hatcher).}
\label{farey_graph}
\end{center}
\end{figure}

Floyd-Hatcher~\cite{floydhatcher82} classify the connected, orientable, incompressible, $\partial$-incompressible surfaces in a torus bundle (excluding the boundary torus itself and the fiber) by edge paths $\gamma$ in the Farey graph diagram of $\text{PSL}_2(\mathbb{Z})$ (see Figure \ref{farey_graph}) which are invariant by the monodromy $\phi$ and minimal, in the sense that no two successive edges of $\gamma$ lie in the same triangle. See Floyd-Hatcher~\cite{floydhatcher82}, Theorem 1.1. The minimality condition implies that $\gamma$ is in fact constrained to lie on a "\textbf{Farey strip}", that is a subset of the Farey Graph consisting of a connected chain of triangles. See Figure \ref{section_types} for some examples of parts of Farey strips. The minimality condition also implies that $\gamma$ cannot divide a fan of the torus bundle in two. I.e. $\gamma$ can travel along either side of the strip, or cross from one side to the other at a border between fans.\\

The vertices of the Farey graph can be viewed as the rational numbers $\frac{a}{b}$, together with $\frac{1}{0}$. Two vertices $\frac{a}{b}$ and $\frac{c}{d}$ are joined by an edge if $ad-bc = \pm 1$. Putting aside the incompressible surface for a moment, we can see how to read off the tetrahedralisation of the punctured torus bundle from the monodromy $\phi$ using the Farey graph. \\

We begin with the punctured torus bundle seen as a cube $[0,1] \times [0,1] \times [0,1]$, minus its vertical edges and with some identifications: We fix a reference basis for the torus taken from the cube edges. We identify the front face with the back, and the left face with the right by translation to obtain $(T^2 \setminus \{0\}) \times [0,1]$. Then identify the bottom with the top, after applying $\phi$ (seen as a linear transformation preserving $\mathbb{Z} \times \mathbb{Z} \subset \mathbb{R} \times \mathbb{R}$) to the bottom face before gluing.\\

As we build $\phi$ from $L$s and $R$s, we can build the punctured torus bundle as a stack of $(T^2 \setminus \{0\}) \times [0,1]$s, one for each $L$ or $R$. Let $\phi_k$ be the $k$th generator (either $L$ or $R$) in the decomposition of $\phi$, where we count from the bottom of the stack, and $\phi = \phi_1 \phi_2 \ldots \phi_N$ (acting on vectors to its right, as usual). Then the basis vectors for the punctured torus, \scriptsize $\left(\begin{array}{c} 1 \\ 0 \end{array}\right)$\normalsize and \scriptsize $\left(\begin{array}{c} 0 \\ 1 \end{array}\right)$\normalsize at the bottom of the stack map up to the $k$th level (boundary between blocks) by $\phi_1 \phi_2 \ldots \phi_k$. Thus we obtain a list of vectors (and so slopes: $\frac{p}{q}$ corresponds to \scriptsize $\left(\begin{array}{c} q \\ p \end{array}\right)$\normalsize). At each level of the stack we have two vectors, the image of the basis vectors at the bottom of the stack under $\phi_1 \phi_2 \ldots \phi_k$. The corresponding slopes are the points on the boundary of our Farey strip. \\

We can also see the tetrahedralisation $\mathcal{H}$ of the torus bundle as a stack of tetrahedra. The six edges of each tetrahedron in the stack have one of four slopes: the bottom and top edges of the tetrahedron have their own slopes, for the other four "middle" edges, opposite edges will have the same slope. Each pair of neighbouring triangles of the Farey strip corresponds to a tetrahedron. The pair of triangles have 4 vertices, corresponding to these 4 slopes of the tetrahedron, the vertices touching both triangles correspond to the slopes of the two pairs of opposite edges of the tetrahedron (the "middle" slopes) and the vertices not touching both triangles correspond to the top and bottom edges of the tetrahedron.\\

One can retrieve the triangulation of the boundary torus induced by the tetrahedralisation of the punctured torus bundle from the Farey strip. Simply take the Farey strip and reflect it across one of its two sides. We now have a "double thickness" strip. Reflect this across one of its sides (equivalent to taking a translate) to obtain a "quadruple thickness" strip. The strip we obtain is combinatorially identical to the triangulation of the boundary torus. One can see this by considering the relationship between the tetrahedra that make up the tetrahedralisation, and the way they must connect to each other based on the slopes that their edges have.  \\

\begin{figure}[ht]
\begin{center}
\includegraphics[width=0.7\textwidth]{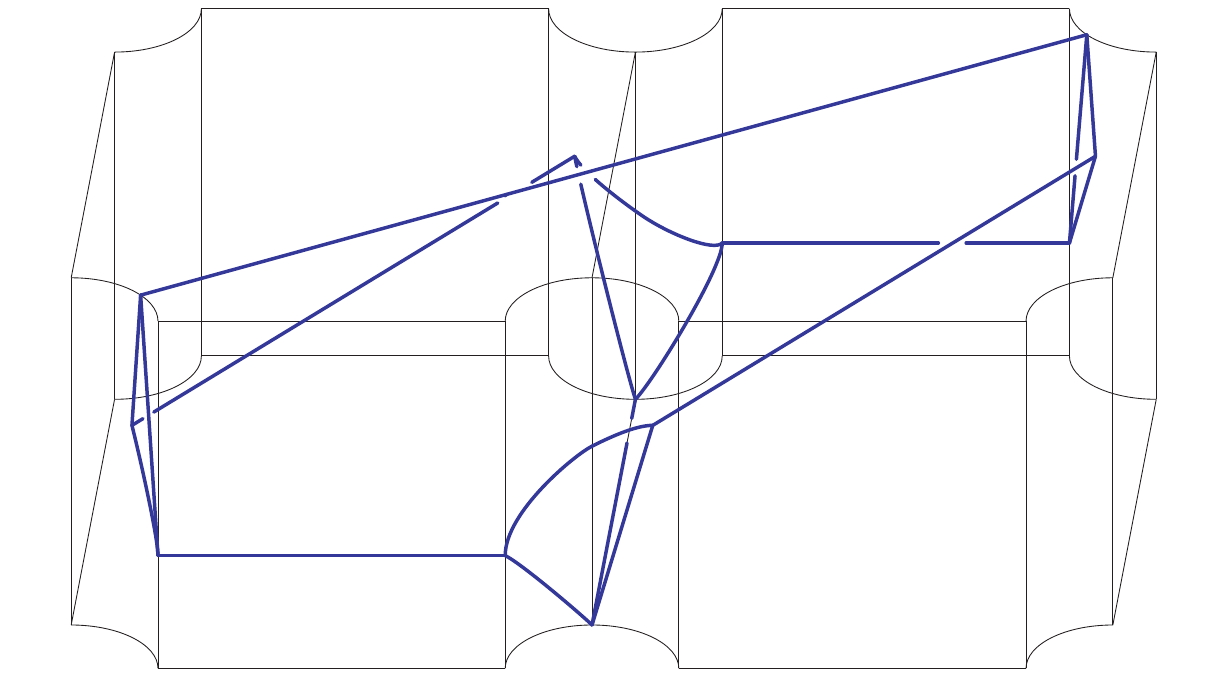}
\caption{Tetrahedron with slopes $\frac{1}{0}$, $\frac{0}{1}$, $\frac{1}{1}$ and $\frac{1}{2}$}
\label{embed_tetra}
\end{center}
\end{figure}

We can map the picture of the torus bundle as a stack of tetrahedra onto the picture of the torus bundle as a stack of $(T^2 \setminus \{0\}) \times [0,1]$ blocks. We require that each edge of a tetrahedron is contained within the first $T^2 \setminus \{0\}$ at the boundary between $(T^2 \setminus \{0\}) \times [0,1]$ blocks that has the correct slope, so we have to allow "stretching out" of the vertices vertically.   Even better, if we cut out a cylindrical neighbourhood around the puncture $\times [0,1]$ and truncate the tetrahedra appropriately we can see them as in Figure \ref{embed_tetra}. We show two copies of a $(T^2 \setminus \{0\}) \times [0,1]$ block for clarity. This diagram actually contains three $(T^2 \setminus \{0\}) \times [0,1]$ layers across which we make our $L$ or $R$ moves, resulting in the four different slopes.\\

Returning to the incompressible surfaces: Floyd-Hatcher index the non-fiber incompressible surfaces by edge paths $\gamma$ in the Farey strip which are invariant by $\phi$ and minimal, in the sense that no two successive edges of $\gamma$ lie in the same triangle.

\begin{figure}[ht]
\begin{center}
\includegraphics[width=0.5\textwidth]{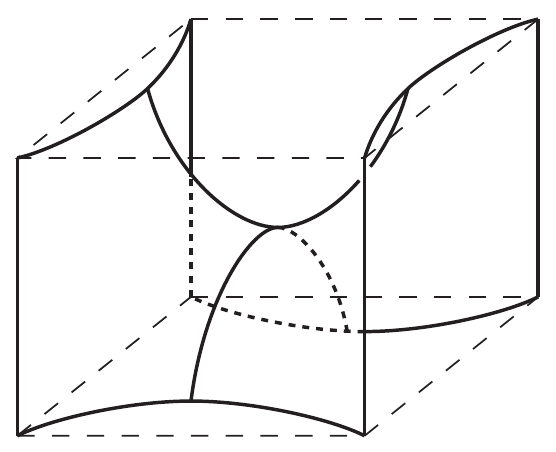}
\caption{Saddle embedded in a cube}
\label{saddle}
\end{center}
\end{figure}

To construct a surface from such a path, we glue together a number of saddles vertically through the stack. See Figure \ref{saddle}. The saddle is embedded in a cube which, after removing the vertical edges and gluing front with back and left with right, we view as $(T^2 \setminus \{0\}) \times [0,1]$. We require one such saddle for each edge of the path $\gamma$. Such an edge joins two vertices of the Farey graph, say $\frac{a_i}{b_i}$ with $\frac{a_{i+1}}{b_{i+1}}$. We transform this saddle by applying in each level the linear transformation:
$$\left(\begin{array}{cc} b_i & b_{i+1} \\ a_i & a_{i+1} \end{array} \right)$$
This has the effect of sending the bottom edges to slope $\frac{a_i}{b_i}$ and the top edges to $\frac{a_{i+1}}{b_{i+1}}$. We insert this block into our stack of $(T^2 \setminus \{0\}) \times [0,1]$ blocks by putting the bottom of the saddle block at the level at which the slope $\frac{a_i}{b_i}$ first appears, and the top of the saddle block at the level at which the slope $\frac{a_{i+1}}{b_{i+1}}$ first appears.\\

We can now see where these surfaces lie with respect to the tetrahedralisation $\mathcal{H}$. There are two cases, depending on if the edge $e$ of $\gamma$ crosses the strip or not:

\subsubsection{$e$ crosses the strip}\label{cross_strip} 
If the edge does cross the strip then there is a pair of neighbouring triangles of the Farey strip, with $e$ as the shared edge. The pair of neighbouring triangles corresponds to a tetrahedron, and the surface has boundary on the four middle edges of the tetrahedron consisting of two pairs with the same slope within each pair. Thus this saddle section of the surface lies as a twisted square in that tetrahedron, separating the top edge from the bottom. In Figure \ref{embed_tetra}, the surface would have boundary on the $\frac{0}{1}$ and $\frac{1}{1}$ edges. The saddle connects (using long thin strips) through to other parts of the surface heading downwards through (the lower down) $\mathcal{N}_\frac{0}{1}$ and upwards through (the higher up) $\mathcal{N}_\frac{1}{1}$. \\

\subsubsection{$e$ does not cross the strip}\label{no_cross_strip}
 The case of an edge $e$ that does not cross the strip is a little more complex. As in the previous case we look for tetrahedra which have edges with slopes which are the boundary of the saddle surface for $e$. There is only one triangle of the Farey strip with $e$ as a boundary, and so two pairs of neighbouring triangles which touch $e$. If we look at the upper of the two pairs of neighbouring triangles, then $e$ joins the bottom slope of the tetrahedron with one of the two middle slopes. We can now see where this saddle is in Figure \ref{embed_tetra} if it joins the bottom slope, $\frac{1}{0}$ to $\frac{1}{1}$: the surface consists of the two lower faces of the tetrahedron, which we push inside the tetrahedron slightly. These two triangles connect to each other through the remaining middle slope edge ($\mathcal{N}_\frac{0}{1}$ in Figure \ref{embed_tetra}) to form the saddle, and connect downwards through the bottom edge ($\mathcal{N}_\frac{1}{0}$) and upwards through the first middle slope edge ($\mathcal{N}_\frac{1}{1}$). We could of course have looked at the position of this saddle on the lower of the two tetrahedra, in which case the saddle would be formed from the upper two faces.\\

The curve $\gamma$ cannot "split apart" a fan of tetrahedra due to its minimality requirement. Thus whenever we have an edge that does not cross the strip, we will have to continue along the side of the entire fan before having the choice to cross the strip instead. The surface section we obtain from going along the entire side of a fan consists of the boundary triangles between each pair of neighbouring tetrahedra in the fan, as well as the boundary triangles between the tetrahedra at the ends of the fan and the hinge tetrahedra next to them.

\begin{defn}
A surface in {\bf Floyd-Hatcher form} is a surface made out of twisted square, triangle and long thin strip pieces as described in sections \ref{cross_strip} and \ref{no_cross_strip}, together with annulus parts in $\partial M \times I$ which extend curves on $\partial M \times \{1\}$ to $\partial M \times \{0\} = \partial M$.
\end{defn}
The pieces used are types of piece allowed by Definition \ref{piece_types}.\\

The last step in the construction of a Floyd-Hatcher incompressible surface is to check if the surface constructed so far is orientable or not. If it is not orientable then it is replaced with the boundary of a small tubular neighbourhood of the original. This has the effect of doubling the number of parallel surfaces in each tetrahedron. In all that follows, the fact that we may actually be manipulating two parallel copies of each piece of surface will not change any of the arguments, and so we will rarely refer to this issue.

\subsection{Surfaces in Yoshida Form}\label{yoshida_form}

In \cite{yoshida91} Yoshida constructs a surface made out of twisted squares and long thin strips from information about the rates and ways in which tetrahedra in the tetrahedralisation are degenerating as we approach an ideal point. The algorithm to generate the surface however only requires the following data:
\begin{itemize}
\item For each $\mathcal{T}_z$ a non-negative integer $k_z$ (which will be the number of parallel twisted squares to put into this tetrahedron).
\item For each $\mathcal{T}_z$ with $k_z \neq 0$, a choice of one of the three pairs of opposite edges of the tetrahedron (this choice is the pair of edges at which the twisted square is not incident).
\end{itemize}
This data is subject to conditions, one condition for each edge of the tetrahedralisation: we insert $k_z$ parallel twisted squares into $\mathcal{T}_z$ oriented according to the choice above, and consider the edges of twisted squares incident at each edge $e$ of the tetrahedralisation. We label the edges of the twisted squares as in figure \ref{fig_twisted_square.pdf}. We require that the number of 0 edges and $\infty$ edges of twisted squares incident at $e$ are equal.\\

To continue the construction, we connect the twisted squares to each other through the $\mathcal{N}_e$ with long thin strips in such a way that we connect 0 edges of twisted squares to $\infty$ edges. Note that there are potentially choices to be made in how the edges of the twisted squares are connected to each other within a $\mathcal{N}_e$. We will choose particular choices here for convenience with the conversion to Floyd-Hatcher. The boundary curves of the twisted squares and long thin strips on $\partial M \times \{1\}$ are either null-homotopic in the torus, or not. We cap off null-homotopic curves with disjoint disk pieces in $\partial M \times I$. We extend any non null-homotopic curves across $\partial M \times I$ with annuli. 

\begin{defn}
A surface in {\bf Yoshida form} is a surface made out of twisted squares, long thin strips, annuli and disks as described above.
\end{defn}

We now describe how Yoshida obtains the data above from an ideal point of the deformation variety: Yoshida restricts to 3-manifolds with one torus boundary component, which implies that the deformation variety is 1 (complex) dimensional, which implies that any ideal point has rational ratios of degeneration of the variables. We multiply up by the least common multiple of the denominators to get integer values, which are the values we choose for the $k_z$. For the degenerate tetrahedra (when $k_z \neq 0$) the choice of pair of edges are those whose complex dihedral angle converges to 1 as we approach the ideal point. Note that the other two pairs converge to 0 and $\infty$, which match up with the labels on the twisted squares as in figure \ref{fig_twisted_square.pdf}. \\

\subsection{Boundary Diagrams of the Surfaces}
\begin{landscape}
\begin{figure}[hp]
\begin{center}
\includegraphics[width=7.0in]{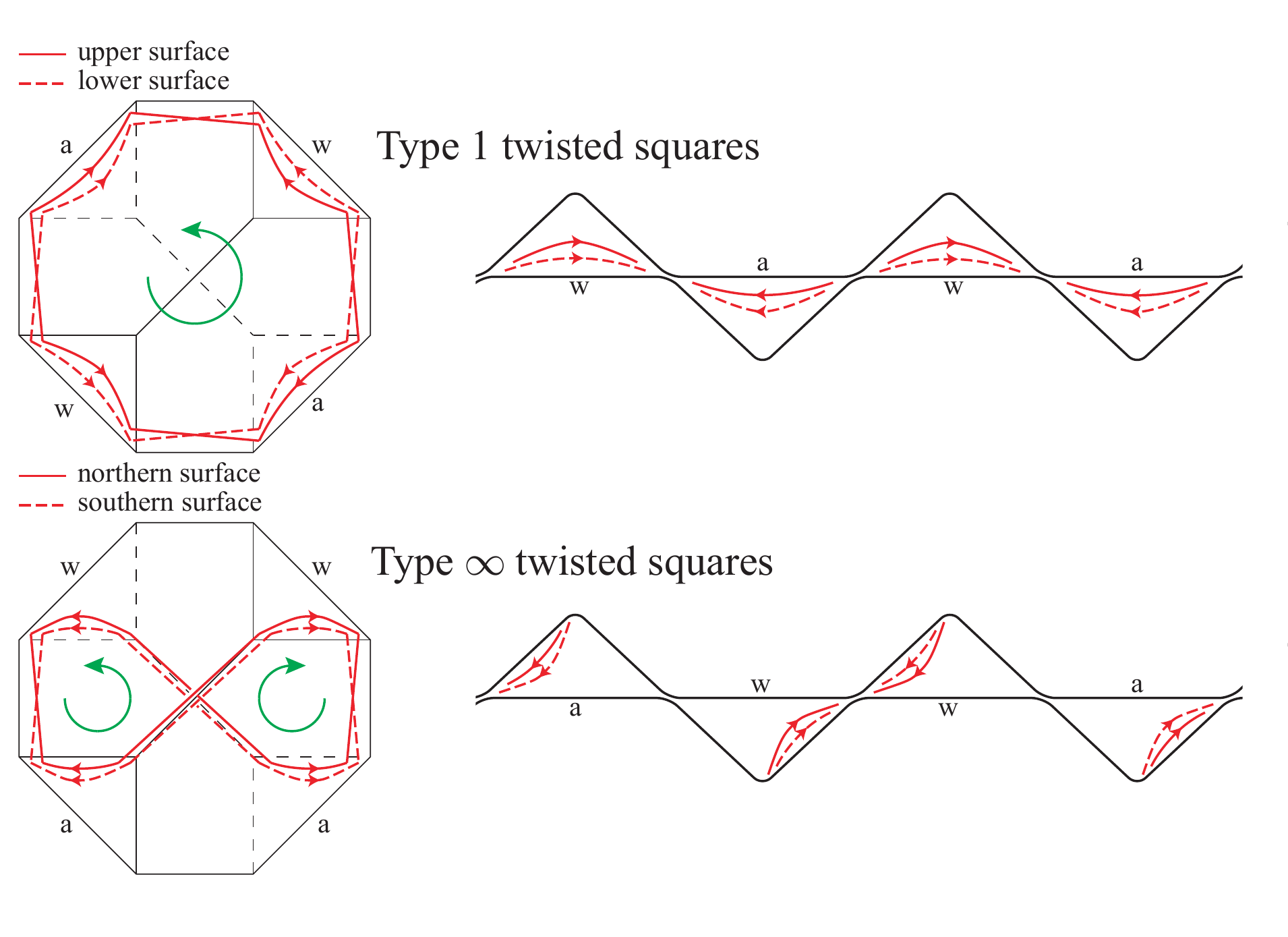}
\caption{Twisted squares of type 1 and $\infty$ in a tetrahedron and the corresponding picture on $\partial M \times \{1\}$.}
\label{tetra_with_surfaces1}
\end{center}
\end{figure}
\end{landscape}

\begin{landscape}
\begin{figure}[hp]
\begin{center}
\includegraphics[width=7.0in]{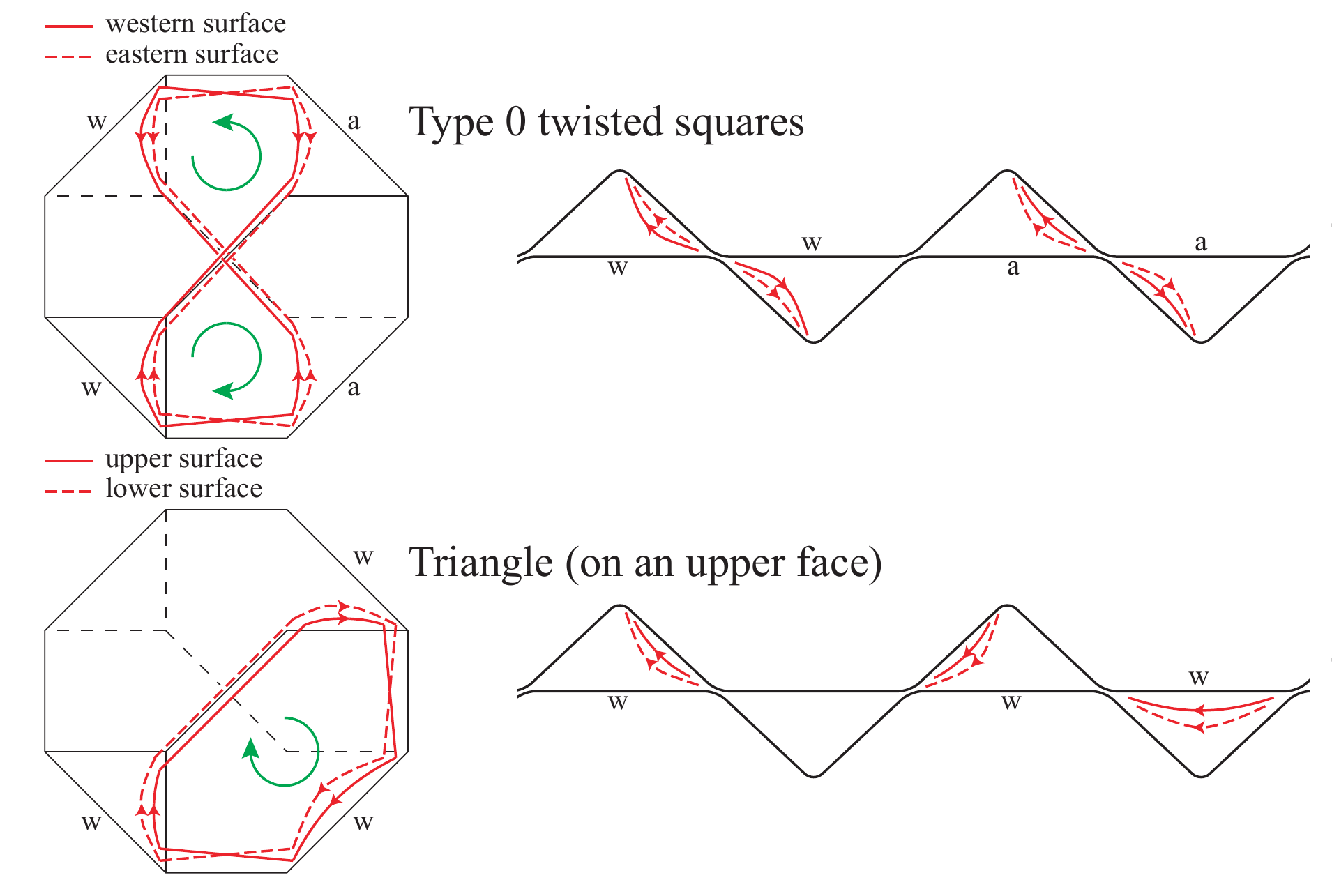}
\caption{Twisted square of type 0 and a triangle in a tetrahedron and the corresponding picture on $\partial M \times \{1\}$.}
\label{tetra_with_surfaces2}
\end{center}
\end{figure}
\end{landscape}

Working with surfaces inside of tetrahedra is difficult and time consuming. Fortunately however, all of the information encoded by a surface in the form given by Definition \ref{piece_types} can be read off from the pattern of the boundary of the pieces of surface on $\partial M \times \{1\}$. We analyse this in Figures \ref{tetra_with_surfaces1} and \ref{tetra_with_surfaces2}.\\

In the figures, we see tetrahedra viewed from above, looking down on the torus bundle. There are three ways to put a twisted square in a tetrahedron, named \textbf{types} 1, $\infty$ and 0, with reference to the complex angle at the top and bottom edges (with respect to the torus bundle) of the tetrahedron. We also show a triangle piece, parallel to one of the upper faces of a tetrahedron. To the right we see the patterns formed on $\partial M \times \{1\}$ (we sometimes refer to this as the {\bf boundary picture}). The orientation on the boundary curves is Yoshida's orientation, which within a triangle on the boundary torus, points from the $\infty$ corner to the 0 corner. This orientation may or may not agree with the orientation induced from the orientation on the triangle or twisted square. The labels "w" and "a" are to be read as "with the induced orientation" or "against the induced orientation". We show two parallel copies of each surface to demonstrate how the ordering of various copies translates to the boundary picture.\\

The curves on the boundary pictures will connect to curves on the boundary of neighbouring tetrahedra, passing through neighbourhoods of the vertex (corresponding to $\mathcal{N}_e$). Within these neighbourhoods of the vertices we see the boundary edges of the long thin strips. In order to respect the product structure on each $\mathcal{N}_e$, we require that the way in which the curves connect to each other through a vertex is consistent with the way in which curves connect at the other end of the edge passing through the manifold. The vertex at the other end of an edge through the manifold can be found by moving two steps along the boundary picture to the right or left, and consistency requires that the picture near one vertex be the mirror image of the picture near the vertex at the other end of its edge. The axis of the reflection is roughly vertical in the boundary pictures. \\

\begin{rmk}
Given a surface in the form described by Definition \ref{piece_types}, we can tell if a surface is orientable by looking at the boundary picture: The curve components on the boundary picture must be consistently oriented according to the induced orientation from a choice of orientation on each twisted square and triangle. A Yoshida form surface contains no triangle pieces, and the curve components are already each oriented with Yoshida's orientation. Showing that such a surface is orientable amounts to showing that half of the curve segment orientations can be reversed, in the ways allowed looking at the diagram, and still having the curve components be oriented.
\end{rmk}
\section{Translating Between the Two Forms of Surface}
\sectionmark{Translating Between the Two Forms}
\label{translation}
In this section we describe the Yoshida form surfaces that will correspond (after isotopies and removal of sphere components in section \ref{surfaces_isotope}) with a given Floyd-Hatcher form surface.\\

\begin{figure}[htb]
\begin{center}
\includegraphics[width=0.7\textwidth]{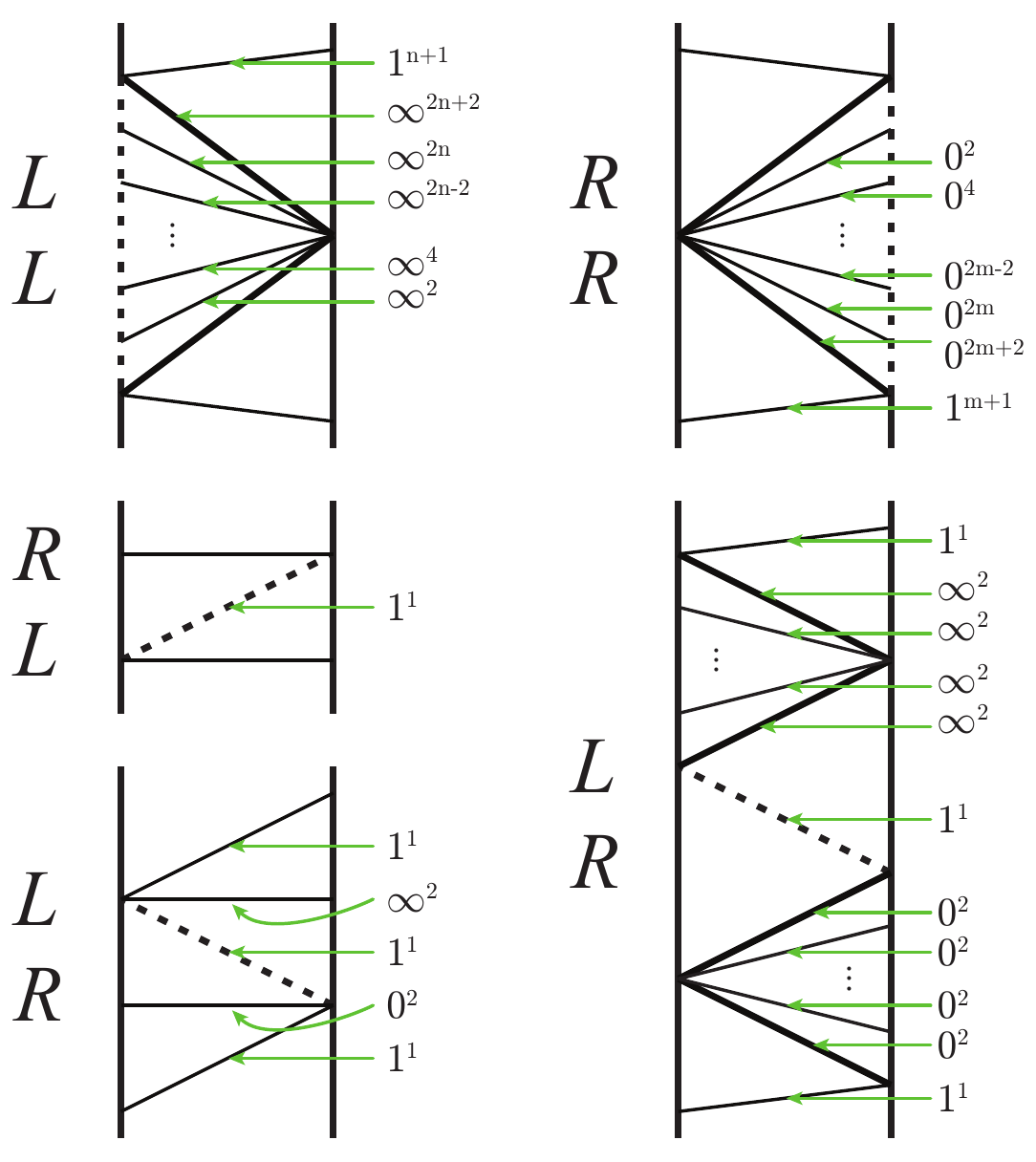}
\caption{The four types of path section, positions on the Farey strip.}
\label {section_types}
\end{center}
\end{figure}

We break a Floyd-Hatcher edge path $\gamma$ into four different types of section, labelled $\,^L_L, \,^R_R, \,^R_L$ and $\,^L_R$. We deal with the Yoshida form for the tetrahedra around each section separately.  The sections are divided at vertices of the Farey strip at which there is a possible choice in which route the path takes, i.e. not at a vertex in the middle of a fan with the path edges travelling along one side of the strip (in an $\,^L_L$ or $\,^R_R$ section). See Figure \ref{section_types}.\\

\begin{landscape}
\begin{figure}[hbtp]

\includegraphics[width=1.8\textwidth]{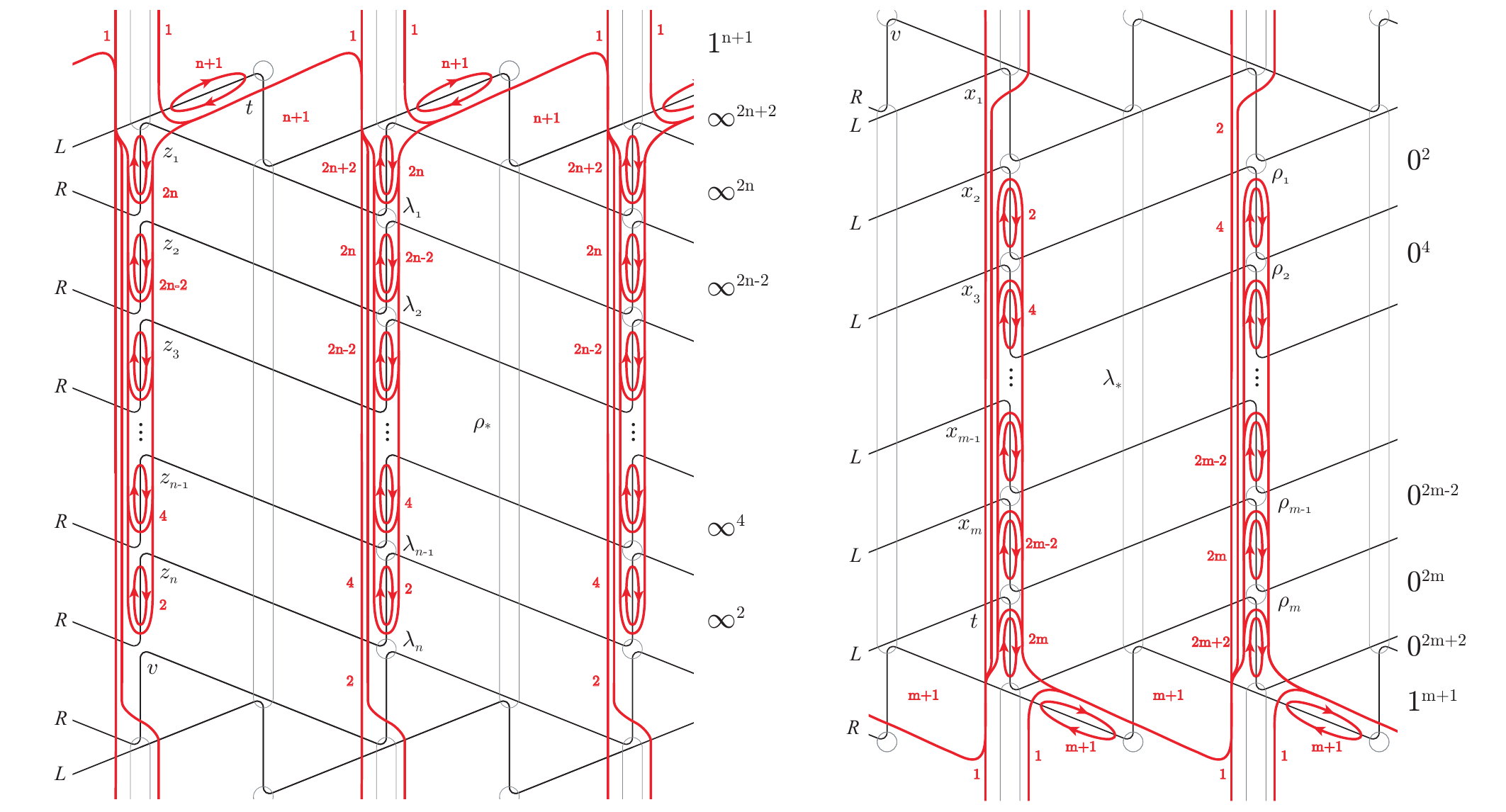}

  \vspace{0pt}
  \hbox{\hspace{1.60in} $\,^L_L$ \hspace{3.85in} $\,^R_R$}
  \vspace{0pt}

  \caption{Torus Boundary pictures for Yoshida form surfaces in sections $\,^L_L$ and $\,^R_R$.}
  \label{fig_L-L_R-R}
\end{figure}
\end{landscape}

\begin{landscape}
\begin{figure}[hbtp]
\includegraphics[width=1.8\textwidth]{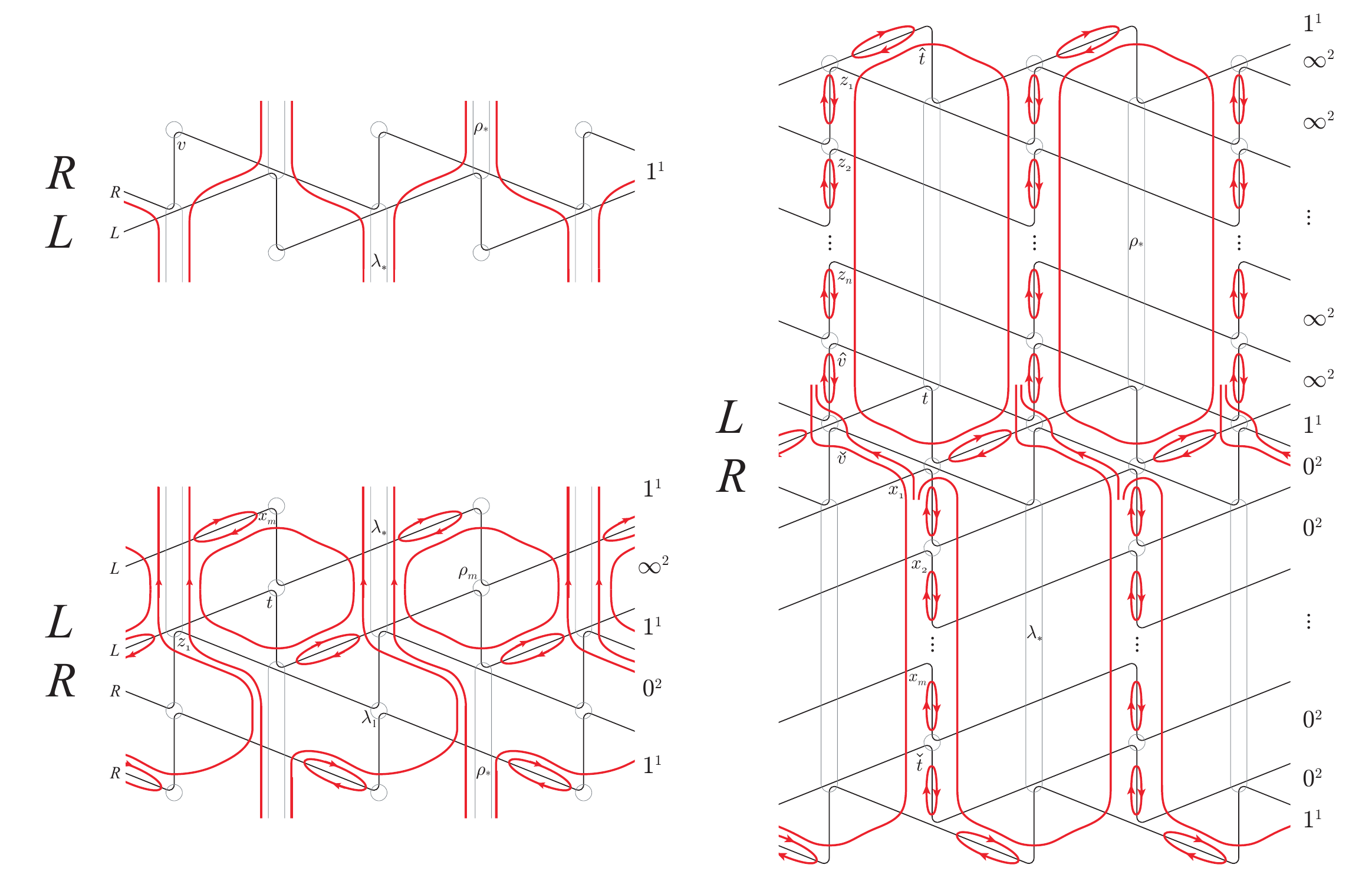}
  \caption{Torus Boundary pictures for Yoshida form surfaces in sections $\,^R_L$ and $\,^L_R$.}
  \label{fig_crossings}
\end{figure}
\end{landscape}

In these diagrams the edges in the Farey strip are solid and the Floyd-Hatcher edge path is dotted. Edges which are at the boundary of a fan are thicker than edges in the middle of a fan, although the edges at the top and bottom of each diagram may or may not be on a fan boundary. The number and type of twisted square ($0$, $1$ or $\infty$) required for the translation to Yoshida form are shown, the arrows point to the edge between the two triangles of the Farey strip that correspond to a tetrahedron with twisted squares. We will see the Yoshida forms of these surfaces in full detail in Figures \ref{fig_L-L_R-R} through \ref{fig_crossings}. The numbers $n$ and $m$ are the numbers of triangles in each respective fan. We give two examples of the $\,^L_R$ case, which depends more than the others on the surroundings within the Farey strip.  In the case for which the fans above and below an $\,^L_R$ have at least two triangles (so there is at least one non-hinge tetrahedron), the pattern is as in the lower left diagram for $\,^L_R$. In the lower right is the picture for the if fans both above and below the $\,^L_R$ have only one triangle. The patterns above and below the $\,^L_R$ edge are independent, so situations with a single triangle fan below and larger fan above look like the top of the lower left diagram joined to the bottom of the lower right diagram and so on.\\

It is worth noting at this point how convoluted the Yoshida forms of these surfaces are in comparison with the Floyd-Hatcher forms. We described the Floyd-Hatcher surfaces corresponding to the curve $\gamma$ on the Farey graph at the end of section \ref{fhform}: For edges of $\gamma$ that cross the strip we have a single type $1$ twisted square inside the hinge tetrahedron corresponding to the two triangles on the Farey strip that meet at the crossing edge. For edges that travel up either side of the strip we have triangle surface pieces on the boundary between each pair of neighbouring tetrahedra in the fan. In contrast, the Yoshida forms of these blocks are considerably more complicated. For one of the two crossing cases (the $\,^R_L$ case) nothing changes, we still have a single twisted square. For the other crossing case some extra surface parts need to be introduced in tetrahedra next to the hinge tetrahedra, it turns out, due to an issue of orientation. For the blocks in which $\gamma$ travels up on side of the strip it seems that a great amount of extra "scrunching up" has to happen. The surface parts that in Floyd-Hatcher form are evenly spread throughout the fan are increasingly bunched up to one side of the fan. We see in fact the number of twisted squares in each tetrahedra "ramping up" by two each time. This corresponds to the tetrahedra degenerating faster and faster as we look along the fan. It isn't intuitively clear why this needs to happen, but is somehow forced by the local gluing equations. \\

The path sections are strung together, and in the cases for which tetrahedra appear in both path section pictures, so are required to collapse by two neighbouring sections, the type of twisted square always agree and the number of each are additive. So, for example, if our path reads from the top $\,^R_R$,  $\,^R_L$,  $\,^L_L$ where the fans between sections have only a single triangle, then the number of twisted squares in the central hinge tetrahedron, at the  $\,^R_L$ is $(m+1) + 1 + (n+1)$. Another example is the lower right $\,^L_R$, which must have a $\,^L_L$ above and $\,^R_R$ below, since the minimality condition precludes the other possible option either side of an $\,^L_R$, namely an $\,^R_L$.\\

\begin{figure}[htb]
\begin{center}
\includegraphics[width=0.7\textwidth]{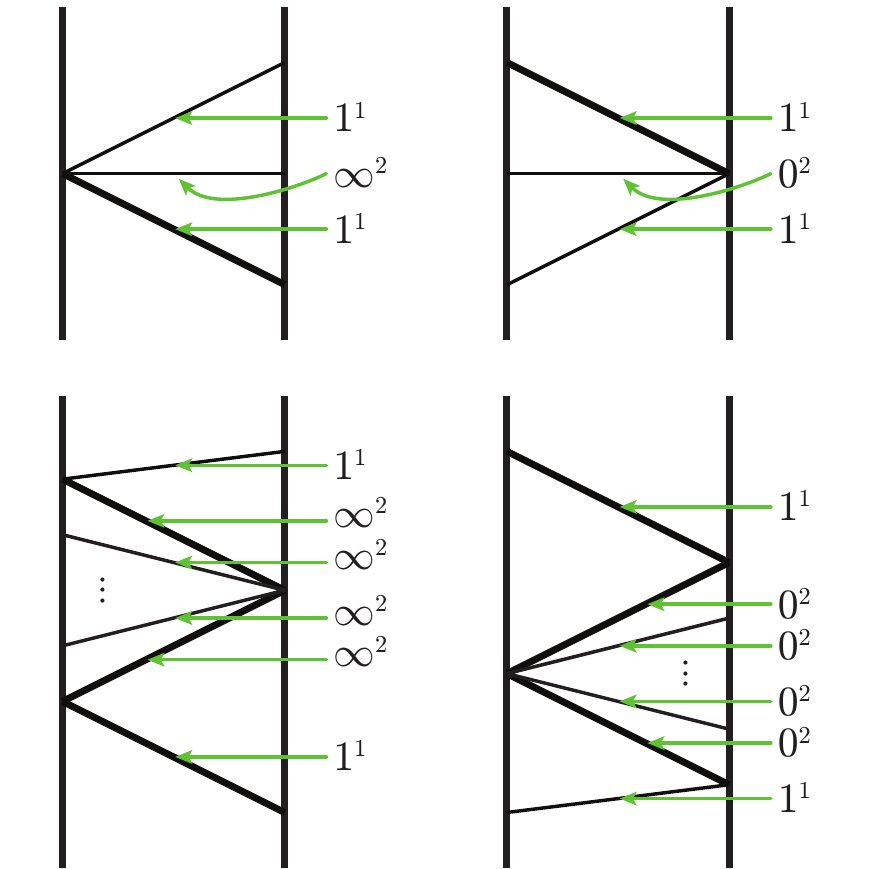}
\caption{Positions of spheres in the Farey strip.}
\label {sphere_types}
\end{center}
\end{figure}

We also introduce in Figure \ref{sphere_types} the numbers and types of twisted squares corresponding to spheres which we will have to add to the Floyd-Hatcher surface in order to solve equations that will come up later in finding ideal points corresponding to these surfaces (the Yoshida forms become yet more complicated). Some number of these spheres will be added either side of $\,^L_R$ sections, and again the numbers of twisted squares are additive. We will show in section \ref{surfaces_isotope} that these do indeed give surfaces isotopic to spheres.\\

We now show (in figures \ref{fig_L-L_R-R} and \ref{fig_crossings}) the boundary torus pictures of the Yoshida form surfaces we defined above, making our choices of which twisted squares are joined to which through an $\mathcal{N}_e$ to simplify the later conversion back to Floyd-Hatcher form. To the right of each diagram we again label the type and number of each twisted square in each tetrahedron. One can check that the picture near one vertex of these diagrams is the mirror image of the picture near the vertex at the other end of its edge (as required for the strip pieces to respect the product structure of $\mathcal{N}_e$). The orientations on the curves are the Yoshida orientation, going from $\infty$ corners to 0 corners of each triangle, and so are necessarily identical for parallel curves going through a triangle. A useful heuristic is that the orientation of an edge within each triangle of the torus boundary picture is always anti-clockwise relative to the center of the triangle.\\

There are often many parallel curve segments (coming from boundaries of the twisted squares) going through the same region on the boundary torus, and so we draw this as a single curve labelled with a number. There are the same number of curve segments going through each of the four triangular truncated ends of each tetrahedron, and the number of such segments entering a junction is equal to the number exiting, so one can quickly work out the number of parallel segments when an edge is not labelled. \\

We name the complex angles in the tetrahedra to match with the labelling in Figure \ref{tetrahedralisation_diag}. We do not label the first and last tetrahedra because those might or might not be hinge tetrahedra, for which we are following a different naming convention.\\

It should be mentioned that we are making a choice here, in that the Yoshida orientation of the surfaces we are constructing always enters our sections from below and exits above. The Floyd-Hatcher surfaces have no inherent orientation, and this choice of direction accounts for the apparent asymmetry between the pictures in, say, the $\,^L_R$ and $\,^R_L$ sections. Had we chosen the arrows to point downwards instead of up, we would effectively rotate all of our pictures 180 degrees and swap the roles of $L$ and $R$.\\

As noted at the end of section \ref{fhform}, we may need to double up the number of surface pieces in each tetrahedron, depending on whether the complete surface, after all blocks are joined together, is orientable. \\

\begin{landscape}
\begin{figure}[hbtp]
\includegraphics[width=1.6\textwidth]{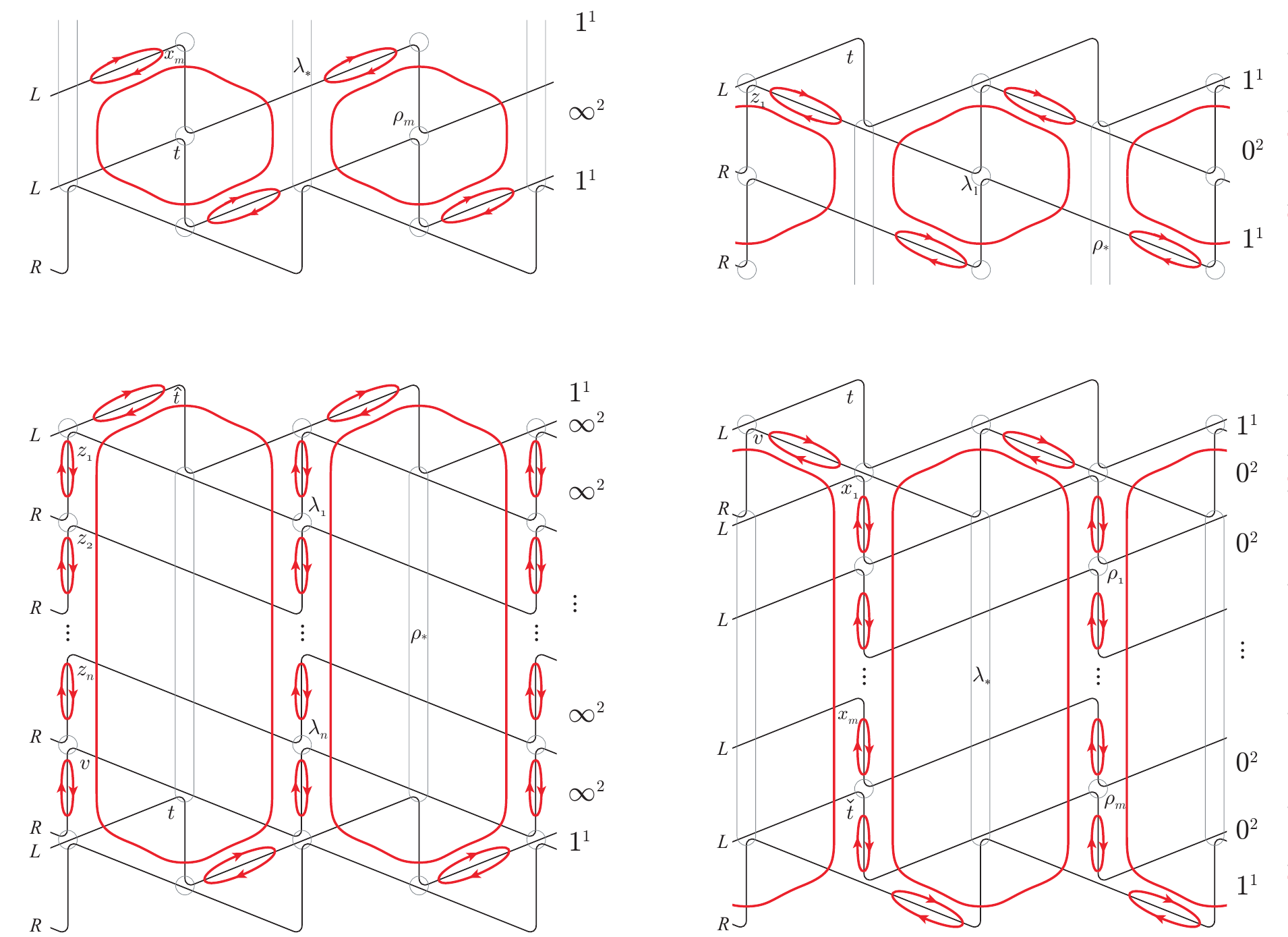}
  \caption{Torus Boundary pictures for Yoshida form spheres.}
  \label{fig_spheres}
\end{figure}
\end{landscape}

\section{Equivalence of the Two Forms of Surface}
\begin{figure}[ht]
\begin{center}
\includegraphics[width=1.0\textwidth]{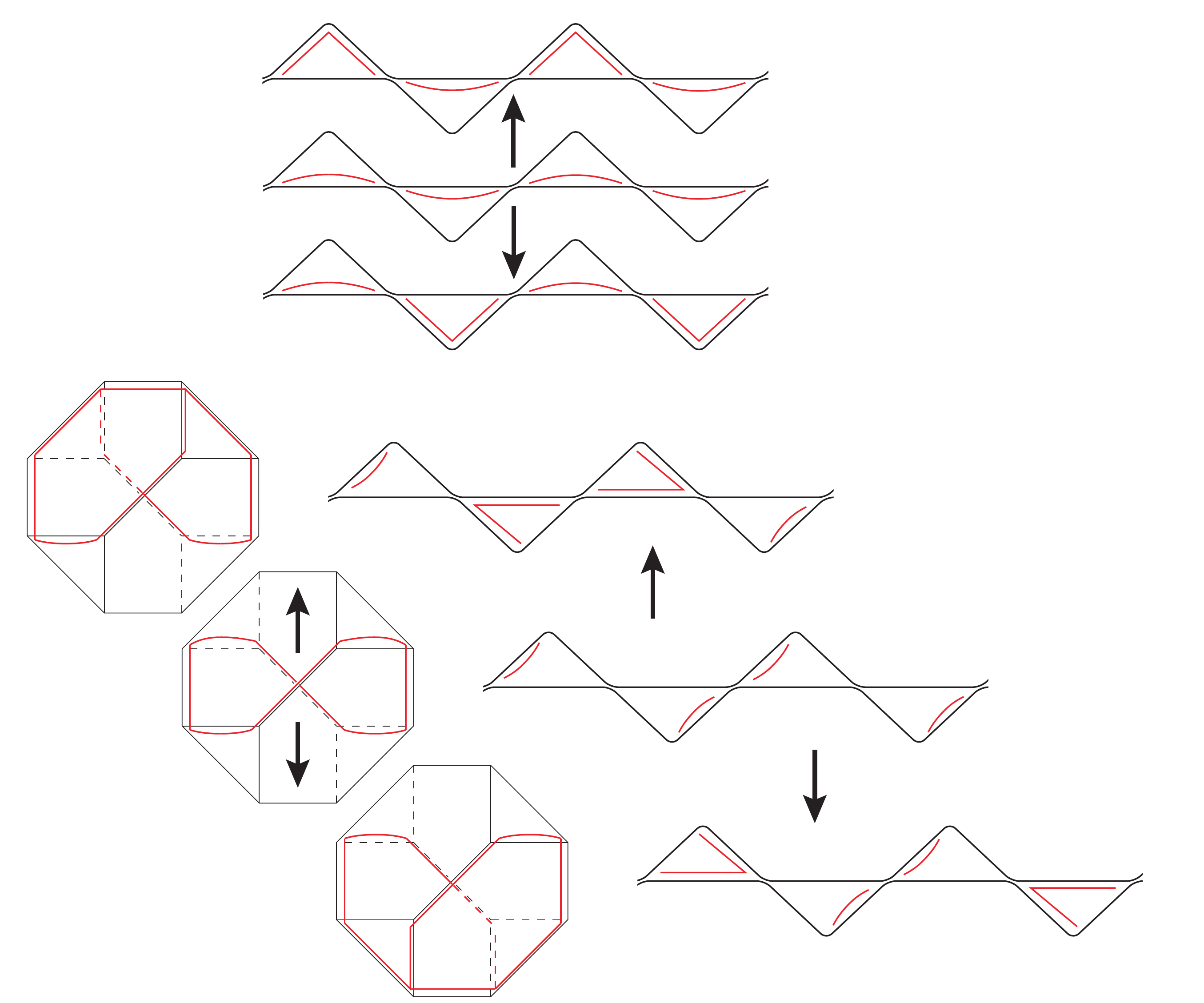}
\caption{Examples of type 1 moves. Above is the picture for twisted squares of type 1 (horizontal with respect to the torus bundle), below is what happens for the other types of twisted squares.}
\label{moves1}
\end{center}
\end{figure}

\begin{figure}[hp]
\begin{center}
\includegraphics[width=0.7\textwidth]{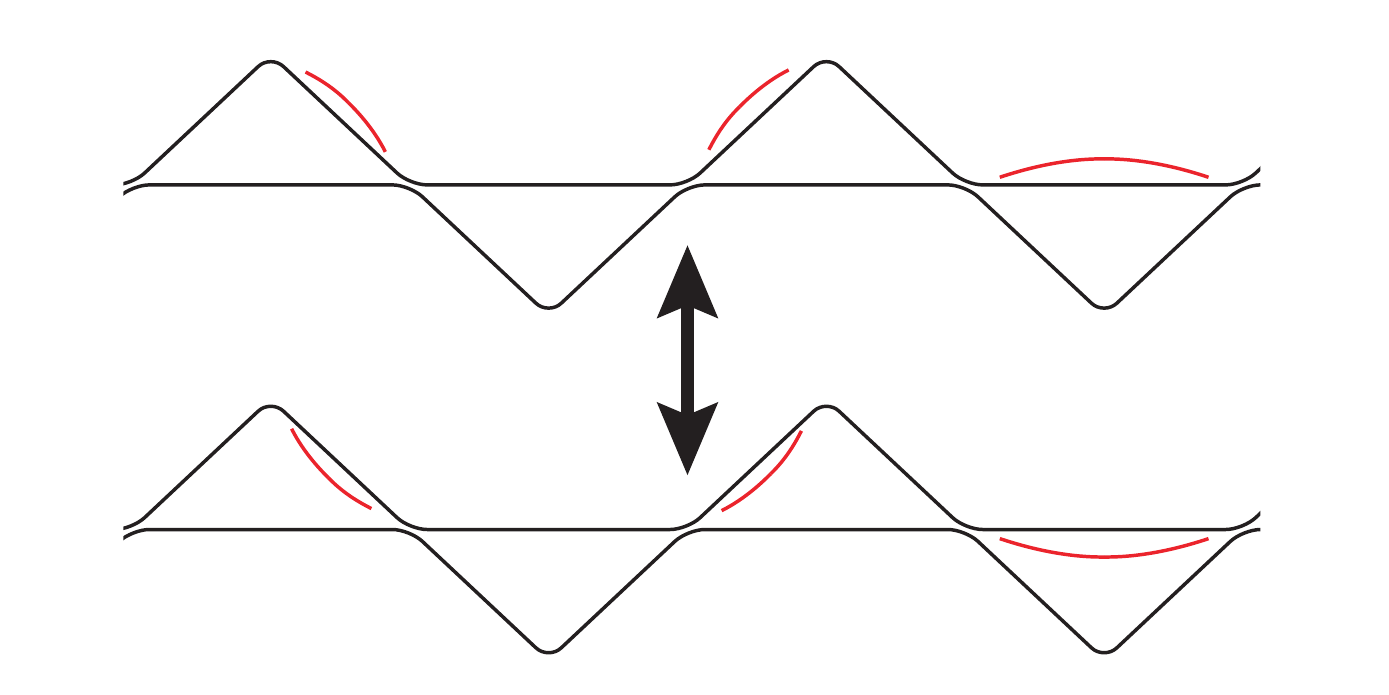}
\caption{An example of a type 2 move.}
\label{moves2}
\end{center}
\end{figure}

\begin{figure}[hp]
\begin{center}
\includegraphics[width=0.7\textwidth]{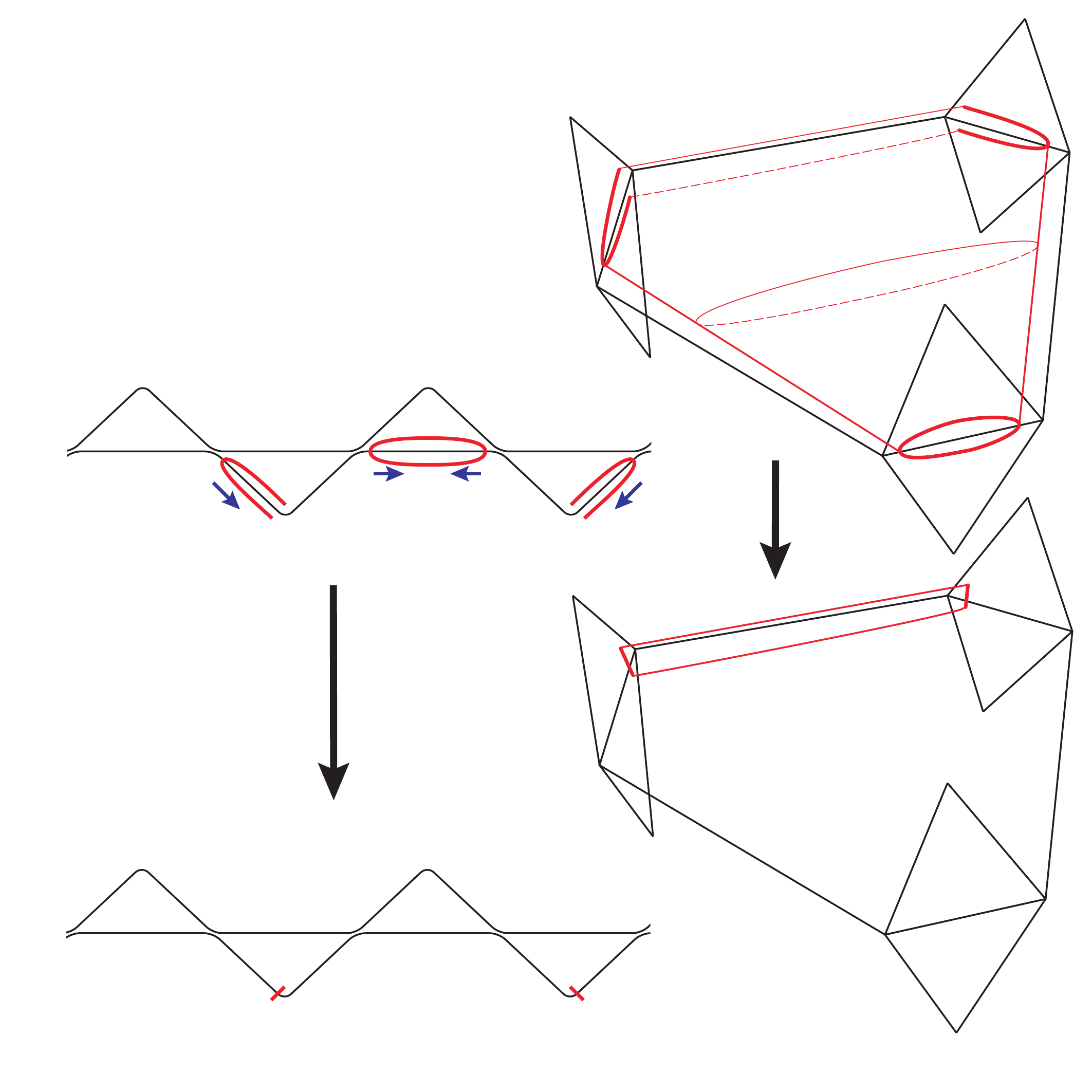}
\caption{An example of a type 3 move.}
\label{moves3}
\end{center}
\end{figure}

\sectionmark{Equivalence of the Two Forms}
\label{surfaces_isotope}
We will employ a number of moves to alter our surfaces. All but one of these kinds of move are isotopies of the surface, the last (Move 4) removes trivial sphere components from the surface.  Moves 1, 2 and 3 are illustrated in Figures \ref{moves1} through \ref{moves3}.

\begin{defn}
The allowed moves are of the following types, all of which may be reversed:
\begin{itemize}
\item{Move 1:} A twisted square in a $\mathcal{T}_z$ may be pushed in one of two directions to become two triangles parallel to the faces of the tetrahedron, with a long thin strip joining the two triangles. The boundary picture changes by an isotopy, and parts in $\partial M \times I$ follow along. 

\item{Move 2:} A triangle parallel to a face of a tetrahedron may be pushed through the tetrahedron face into the neighbouring tetrahedron, to become a triangle within and parallel to the face of the second tetrahedron. This move may introduce or remove strips through the $\mathcal{N}_e$ and does so in the obvious way. Again the parts in $\partial M \times I$ follow along.

\item{Move 3:} Boundary bigon removal: if we have a bigon curve component on the boundary picture (so it is capped off by a disk part in $\partial M \times I$), where the two surface pieces the bigon is a boundary of are both triangles, necessarily each parallel to the shared face between neighbouring tetrahedra, then we may isotope away the cap and two triangles to leave only a strip, as detailed in Figure \ref{moves3}. As well as removing a disk part, there are some isotopies of curves on $\partial M \times \{1\}$ and again parts in $\partial M \times I$ follow along. 

\item{Move 4:} Removal of sphere components: The sphere components we will see will be of the form of a cylinder of strips within a $\mathcal{N}_e$, surrounding the edge $e$, with two disk caps, one on either end. 
We delete the whole sphere.\\

\end{itemize}
\end{defn}

We consider the moves needed on the different types of path section individually.\\

\begin{figure}[ht]
\begin{center}
\includegraphics[width=0.95\textwidth]{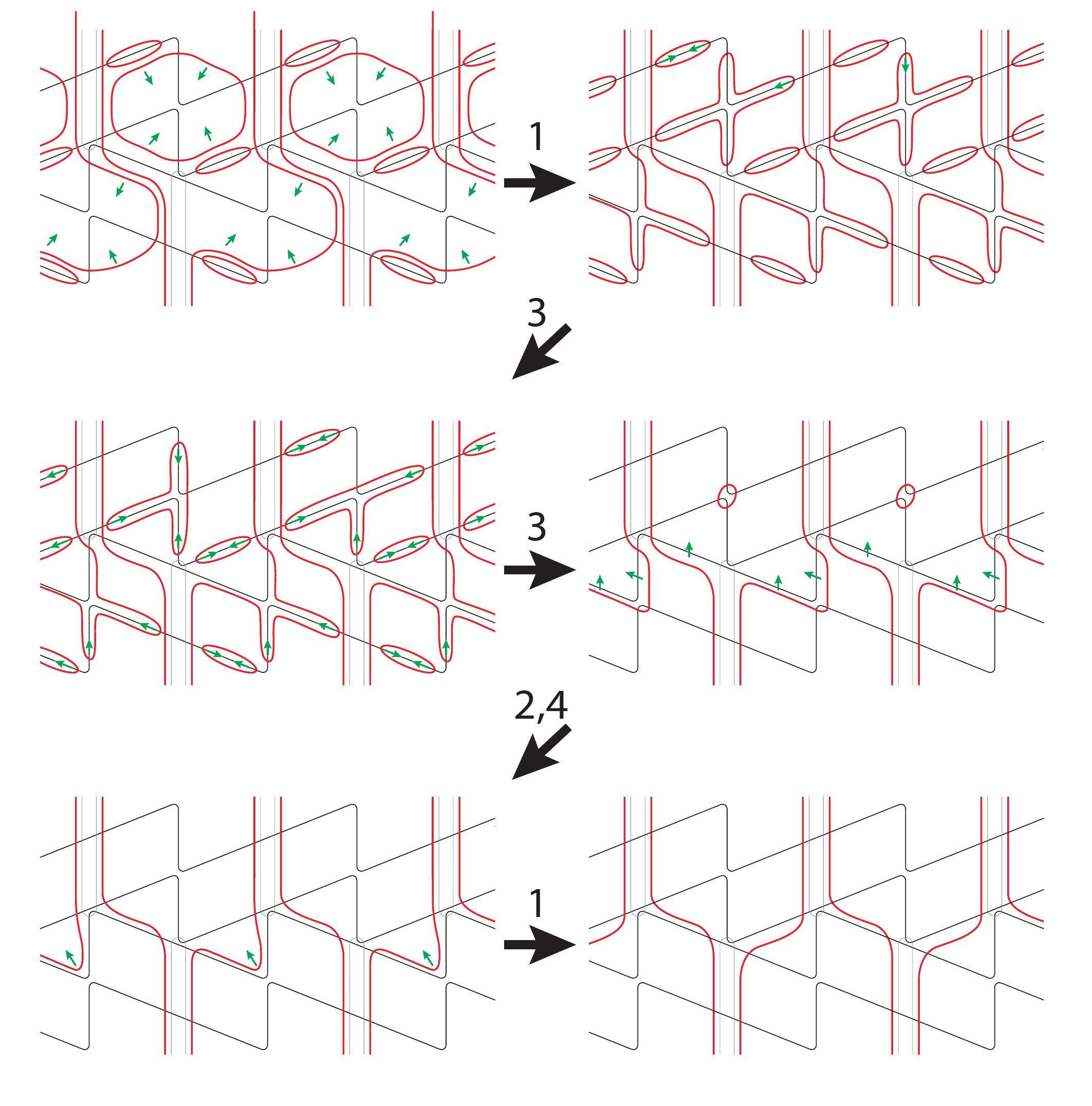}
\caption{Moves to convert a "small" $\,^L_R$ section from Yoshida form to Floyd-Hatcher form.}
\label{isotope_crossing}
\end{center}
\end{figure}

In the case $\,^R_L$, nothing need be done. The single twisted square in Yoshida form is already in the hinge tetrahedron, the correct place for the Floyd-Hatcher form. The $\,^L_R$ case needs some work. We describe the sequence of moves for the case in which the fans above and below the $\,^L_R$ have at least two triangles in Figure \ref{isotope_crossing}.\\

We draw arrows to show how each move is being used. From the first picture to the second, we apply move 1 seven times. From the second to the third we show one use of move 3, then from the third to the fourth we do the rest of the move 3s. From the fourth picture to the fifth we remove a sphere component by move 4 and use move 2 twice to push two triangles into the hinge tetrahedron. Finally, from the fifth picture to the sixth we use move 1 once more, to push two triangles inwards, to form a twisted square. \\

Note that the sphere component we removed in this process is exactly one of the types of sphere from Figure \ref{fig_spheres}. We include a sphere in the $\,^L_R$ section "for free" to simplify calculations later on. It turns out that at least one such sphere must always be present in that spot (if more than one, we mean concentric spheres). It is also worth noting that there is a symmetry between the top and bottom of the Yoshida form for an $\,^L_R$ section: we could have the sphere be the lower $1^1 0^2 1^1$ and the part of the incompressible surface be the upper $1^1 \infty^2$.\\

This case demonstrates the procedure for removing "small" spheres, and the required moves are analogous for the larger spheres. Likewise for the larger versions of the $\,^L_R$ section.\\

\begin{landscape}
\begin{figure}[hp]
\begin{center}
\includegraphics[width=1.7\textwidth]{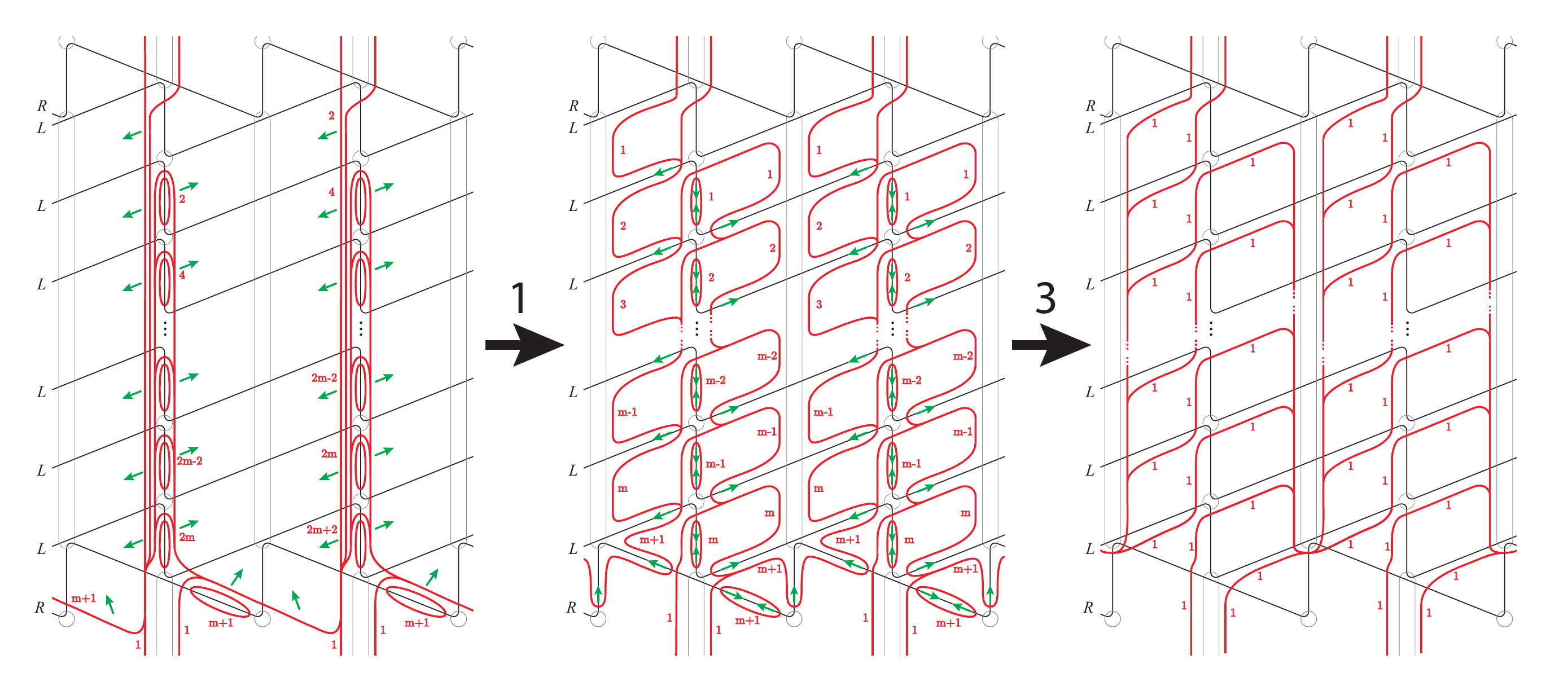}
\caption{Converting a $\,^R_R$ join from Yoshida form to Floyd-Hatcher form.}
\label{isotope_R-R}
\end{center}
\end{figure}
\end{landscape}

We now look at the $\,^L_L$ and $\,^R_R$ sections. We need only do one of them, since (ignoring the Yoshida orientation arrows) reflecting the diagram for  $\,^L_L$ across a horizontal line and translating horizontally gives us the diagram for $\,^R_R$. The moves we employ do not care which way up the diagram is. We demonstrate the sequence of moves in the $\,^R_R$ case in Figure \ref{isotope_R-R}.\\

As discussed at the end of section \ref{fhform} we expect to get the triangles between each pair of neighbouring tetrahedra in the fan, plus the triangles between each tetrahedron on the end of the fan and the hinge tetrahedron next to it, which is the state in the third diagram of the figure. All edges connect through the 4-valent vertices in the diagrams vertically (with respect to the torus bundle vertical).\\

We have shown that each section of the surface can be individually converted from Yoshida form to Floyd-Hatcher form, so the surface as a whole is equivalent to the corresponding surface in Floyd-Hatcher form. None of the preceding arguments are changed if we needed to double up each surface piece (if the complete surface would otherwise be non-orientable). By their result, Theorem 1.1 of \cite{floydhatcher82}, we can construct in this way all incompressible surfaces in the torus bundle, other than the fiber $T^2 \setminus \{0\}$ and the peripheral torus.\\

We now would like to show that the surfaces we have constructed correspond to ideal points of the deformation variety, but first we will need to identify places at which we may need to add extra spheres. These spheres will of course not alter the eventual surface we obtain when we convert back from Yoshida to Floyd-Hatcher form, since all spheres are removed in that process. They will however alter the numbers of twisted squares in some tetrahedra, and hence the rates at which those tetrahedra are supposed to degenerate as we approach our ideal point.\\

Before that however, we are now in a position to be able to classify which surfaces are in fact semi-fibers.

\section{Semi-fibers}
\sectionmark{Identifying Semi-fibers, Ideal Points}
\label{semi-fibers}

\begin{defn}
A \textbf{semi-fibration} is a 3-manifold formed by taking two copies of a twisted I-bundle over a non-orientable surface and gluing them to each other along their (orientable) boundaries. The orientable boundary becomes the \textbf{semi-fiber} in the semi-fibration.
\end{defn}

\begin{defn}\label{tight}
A connected subset of a path $\gamma$ in the Farey graph of at least two edges is \textbf{tight} if at each vertex of the sub-path the two edges leaving that vertex belong to neighbouring triangles.
\end{defn}
In other words, when the path reaches a vertex it takes either the "second left" or the "second right" turn. Taking the "first right" or "first left" is prohibited by the minimality condition on paths, so in this sense a tight sub-path turns as tightly as it possibly can.\\

\begin{prop}
An incompressible surface in a punctured torus bundle (other than the fiber) is a semi-fiber if the whole of the corresponding path in the Farey graph is tight.
\end{prop}

\begin{proof}
If the path $\gamma$ in the Farey graph is tight, then (after conjugating) the path near a given vertex looks like either the path from $\frac{-1}{1}$ to $\frac{1}{0}$ to $\frac{1}{1}$ or the path from $\frac{-1}{1}$ to $\frac{1}{0}$ to $\frac{1}{1}$ (see Figure \ref{farey_graph}). Following the Floyd-Hatcher construction as in section \ref{fhform}, we obtain two saddles corresponding to the two edges of the Farey graph here, which for the case of $\frac{-1}{1}$ to $\frac{1}{0}$ to $\frac{1}{1}$ are shown in the left diagram of Figure \ref{saddles_semi-fiber}. The picture for $\frac{-1}{1}$ to $\frac{1}{0}$ to $\frac{1}{1}$ is obtained by reflecting this picture in the line of slope $\frac{-1}{1}$ and the arguments go through similarly.
\begin{figure}[hp]
\includegraphics[width=1.0\textwidth]{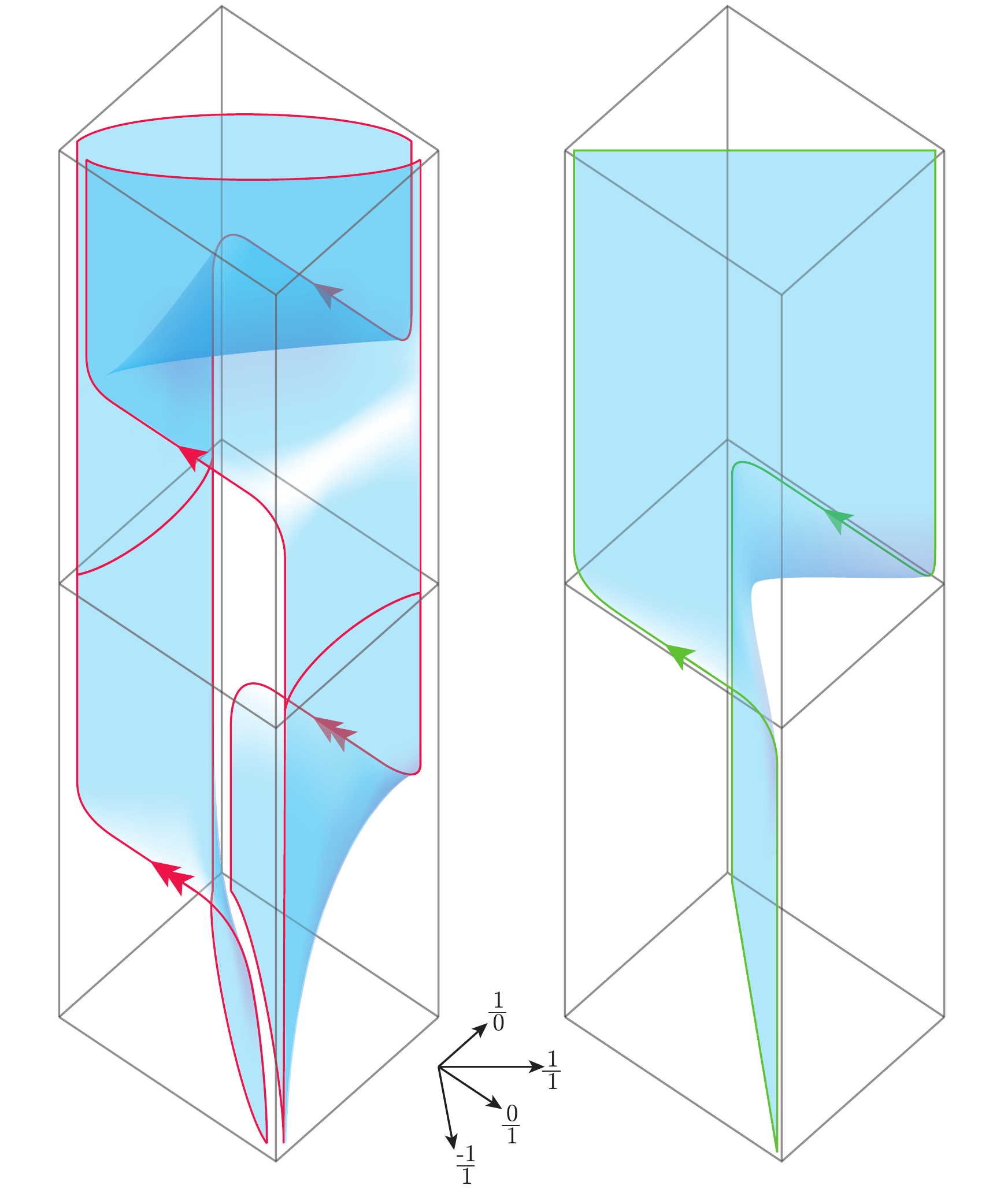}
\caption{Saddles coming from a tight path in the Farey graph, and the non-orientable surface obtained by flowing from these saddles "inwards".}
\label{saddles_semi-fiber}
\end{figure}
The two saddles have both been sheared from the saddle in Figure \ref{saddle} according to the slopes they are supposed to have above and below, so in this diagram each saddle is made up of two pieces, joined to each other under the identification of the front and back faces of the tall cuboid when we glue them together to form the punctured torus bundle.\\

The surface made up of the two saddles cuts the tall cuboid (with identifications on the four vertical faces) into two pieces. It isn't too hard to see how we can flow from the two saddles in the left diagram "inwards", until the surface meets itself resulting in the non-orientable surface in the right diagram. Flowing in the other direction similarly leads to the surface meeting itself. This is best seen by looking at the vertices of $\gamma$ one to either side of the vertex we are currently examining, and following the same observations as for flowing "inwards" here.\\

Consider first the case for which the path $\gamma$ has an even number of segments, and so the surface made from the saddles is orientable. Then the entire torus bundle is split into two pieces by the saddles surface, and each piece is a twisted I-bundle over a non-orientable surface formed by gluing together copies of the surface in the right diagram of Figure \ref{saddles_semi-fiber} vertically.\\

In the case that $\gamma$ has an odd number of segments, the Floyd-Hatcher construction has us take the double cover of the resulting non-orientable saddles surface, i.e. taking the boundary of a small neighbourhood of the saddles surface. In this case the non-orientable saddles surface does not separate the punctured torus bundle, but the double of it does. One of the twisted I-bundles is then the neighbourhood of the saddles surface and is a twisted I-bundle over the saddles surface. The other I-bundle is the rest of the punctured torus bundle, and is again a twisted I-bundle over the non-orientable surface formed from copies of the surface in the right diagram of 
Figure \ref{saddles_semi-fiber}.
\end{proof}

\begin{rmk}
Fibers and semi-fibers cannot be produced from ideal points of the deformation variety or via the Culler-Shalen construction from ideal points on components of the character variety consisting of irreducible representations. There is however a well known construction of fibers and semi-fibers from relatively simple actions on trees, which can be obtained from ideal points on components of the character variety consisting of reducible representations. In this sense the results of this paper show that all incompressible surfaces, including the fibers and semi-fibers, are detected by the character variety.  
\end{rmk}

\section{Adding spheres to the Yoshida form surfaces}
\label{adding_spheres}

In order to show that an ideal point corresponding to a certain Yoshida form surface exists, we need to show the existence of a sequence of finite tetrahedra shapes for the manifold that degenerate with the rates and directions of collapse corresponding to the Yoshida form surface. \\

In section \ref{tilde} we will examine limiting versions of the gluing equations (obtained by a blow up) as the tetrahedra corresponding to our Yoshida form surface degenerates, and find explicit solutions to those equations at the ideal point, which will later give us the existence of nearby finite tetrahedra shapes. At present however, we are only interested in the conditions these equations impose on the form of the twisted squares surface, not yet in actual numerical values. So far we have only paid attention to the following condition, which we give a name to here:

\begin{defn}
At an edge $e$ for which the surface passes through $\mathcal{N}_e$, the \textbf{"$0 - \infty$" matching condition} states that the number of twisted squares with an $\infty$-side in $\mathcal{N}_e$ must match the number of twisted squares with a $0$-side. 
\end{defn}

It is easy enough to check that this condition holds in Figures \ref{fig_L-L_R-R}, \ref{fig_crossings} and \ref{fig_spheres} above. If we are to be able to approach an ideal point via a sequence of finite tetrahedra shapes, this condition must hold since the product of angles around an edge must always be 1 (and hence be bounded, so cannot converge to 0 or $\infty$, which it would have to if the numbers of $0$ and $\infty$ angles were not balanced).\\

There are however some more subtle requirements resulting from edges $e$ for which no surface passes through $\mathcal{N}_e$, but for which many or all of the neighbouring tetrahedra are degenerating. The situation to consider is of an edge $e$ surrounded by tetrahedra, \em all but one \em of which are degenerating such that their complex dihedral angles at $e$ are converging to $1$ (the other two dihedral angles of each tetrahedron converging to $0$ and $\infty$, and hence the twisted square(s) in each tetrahedron connecting the corresponding edges). The one other tetrahedra we imagine is not degenerating. The gluing equation around $e$ then requires that the product of the complex dihedral angles around it be equal to $1$. Since all the degenerating tetrahedra contribute complex angles of $1$, the complex angle at the supposedly non-degenerating tetrahedron is also forced to be $1$, and it is forced to degenerate.\\

The more general condition that this example fails to satisfy is given in the following definition:

\begin{defn}
\label{non-unique_min}
At an edge $e$ for which no complex angle is degenerating to $0$ or $\infty$, the angles are either not degenerating or degenerating to $1$ at some integer rate (corresponding to the number of twisted squares in the tetrahedron). Viewing the non degenerating angles as "degenerating at rate $0$", the \textbf{non-unique minimum rate condition} states that the minimum degeneration rate over all tetrahedra around $e$ is not achieved at a unique tetrahedron.
\end{defn}

This condition is necessary to be able to solve the equations that will follow (see section \ref{sec_non-unique_min} for more detail on where this condition comes from), and in the punctured torus bundle case is enough (it turns out) to ensure a solution. However, this condition (together with the condition matching numbers of $0$-edges with $\infty$-edges) is not obviously sufficient in general. \\

In our case, the edges for which the non-unique minimum rate condition applies are seen as the vertices in the boundary diagrams either side of an $\,^L_R$ crossing, in fact the vertices in the center of the truncated ends of the spheres we will be adding. We call these vertices \textbf{sphere vertices}. The simplest example is that of the "small" $\,^L_R$ crossing (see Figure \ref{fig_crossings}), assuming no parts of surface above or below that section affect the rate of degeneration of the outer $1$s. In this situation the non-unique minimum rate condition already holds. The relevant vertices are labelled $\rho_m$ and $\lambda_1$, and the minimal rate is $1$ in both cases, achieved in two tetrahedra in both cases. \\

The next situation to consider is a "large" $\,^L_R$ crossing, for which we must add a  $\,^L_L$ above and  $\,^R_R$ below, resulting in the types and rates:

\[ 
\begin{array}{llllll}
\,^L_L & \,^L_R & \,^R_R & \text{sum} & \text{upper sphere} & \text{lower sphere} \\
\hline
1^{n+1}       & 1^1      &           & 1^{n+2}       & 1^1      &        \\
\infty^{2n+2} & \infty^2 &           & \infty^{2n+4} & \infty^2 &        \\
\infty^{2n}   & \infty^2 &           & \infty^{2n+2} & \infty^2 &        \\
\vdots        & \vdots   &           & \vdots        & \vdots   &        \\
\infty^2      & \infty^2 &           & \infty^4      & \infty^2 &        \\
              & \infty^2 &           & \infty^2      & \infty^2 &        \\
              & 1^1      &           & 1^1           & 1^1      & 1^1    \\
              & 0^2      &           & 0^2           &          & 0^2    \\
              & 0^2      &  0^2      & 0^4           &          & 0^2    \\
              & \vdots   &  \vdots   & \vdots        &          & \vdots \\
              & 0^2      &  0^{2m}   & 0^{2m+2}      &          & 0^2    \\
              & 0^2      &  0^{2m+2} & 0^{2m+4}      &          & 0^2    \\
              & 1^1      &  1^{m+1}  & 1^{m+2}       &          & 1^1
\end{array}
\]
Note first that the powers and types of degeneration in the "small" $\,^L_R$ crossing (see figure \ref{fig_crossings}) fits into the pattern in the "sum" column here if we allow $n$ or $m$ to be $-1$. Here the non-unique minimum rate condition does not yet hold either above or below (assuming $n \neq -1 \neq m$). However, we are allowed to add spheres to try to satisfy the condition. First note that with the allowed additions, the minimum rates at the vertex within, say, the upper sphere can only involve the upper $1$ tetrahedron, the middle $1$ tetrahedron and the $\infty$ tetrahedron immediately above it. All of the other $\infty$ tetrahedra necessarily have a faster rate since they start off with a faster rate, and we can only add to all $\infty$ rates equally. Similarly for the lower vertex, and we delete the irrelevant rows to reduce the problem to this form: 
\[ 
\begin{array}{lll}
\text{sum} & \text{upper sphere} & \text{lower sphere} \\
\hline
1^{n+2}       & 1^1      &        \\
\infty^2      & \infty^2 &        \\
1^1           & 1^1      & 1^1    \\
0^2           &          & 0^2    \\
1^{m+2}       &          & 1^1
\end{array}
\]

We rewrite this one further time, removing reference to the type of collapse, and setting $\alpha = n + 1$, $\beta = m + 1$, so that we need to be able to solve the problem (now essentially a problem in integer linear programming) of adding some numbers, $a$ and $b$ of upper and lower spheres to satisfy the non-unique minimum rate condition for $\alpha, \beta \geq 0$.

\[ 
\begin{array}{cccc}
\text{Sum} & \text{Upper Sphere} & \text{Lower Sphere} & \text{Sum} + a \text{(Upper Sphere)}  + b \text{(Lower Sphere)}\\
\hline
\alpha+1      & 1 &   & \alpha + a + 1\\
       2      & 2 &   & 2a + 2        \\
       1      & 1 & 1 & a + b + 1     \\
       2      &   & 2 & 2b + 2        \\
\beta +1      &   & 1 & \beta + b + 1 
\end{array}
\]

The solution:
\begin{enumerate}
\item If $\alpha < \beta$ set $a = \alpha + 1, b = \alpha$
\item If $\alpha = \beta$ set $a = \beta, b = \alpha$
\item If $\alpha > \beta$ set $a = \beta, b = \beta + 1$
\end{enumerate}
The rates for the three cases then look like:
\[ 
\begin{array}{rrr}
\alpha < \beta & \alpha = \beta & \alpha > \beta \\
\hline
2\alpha+2         & \alpha + \beta +1 & \alpha + \beta +1 \\
2\alpha+4         & 2\beta + 2        & 2\beta+2          \\
2\alpha+2         & \alpha + \beta +1 & 2\beta+2          \\
2\alpha+2         & 2\alpha + 2       & 2\beta+4          \\
\alpha + \beta +1 & \alpha + \beta +1 & 2\beta+2
\end{array}
\]
It isn't hard to check that the non-unique minimum rate condition is now satisfied for the two "sphere vertices". I.e. for each of the three types of solution, the minimum value in the first three entries in the column is achieved in at least two places, and the same for the last three. If $\alpha +1 = \beta$, the first two columns actually give identical parameters, and the minimum value for the last three rows is achieved in all 3 places. A similar observation can be made for if $\alpha = \beta + 1$.\\

More generally, immediately above or below an $\,^L_R$, or "extended" $\,^L_R$ (a "large" $\,^L_R$, including the $\,^L_L$ above and $\,^R_R$ below) could be a $\,^R_L$, with no non-degenerating tetrahedra between. The effect this has is to add one to $\alpha$ or $\beta$ (respectively), and can be solved using the above scheme. However, immediately next to the  $\,^R_L$ could be another "extended" $\,^L_R$, and now the spheres that we add to solve one $\,^L_R$ start to interfere with the solution of the other. To be explicit, the new problem we would have to solve would be of the form:

\[ 
\begin{array}{rrrrrr}
\text{Surface} & S_1 & T_1 & S_2 & T_2 & \text{Surface} + a_1 S_1  + b_1 T_1 + a_2 S_2  + b_2 T_2\\
\hline
\alpha+1      & 1 &   &   &   & \alpha + a_1 + 1\\
       2      & 2 &   &   &   & 2a_1 + 2        \\
       1      & 1 & 1 &   &   & a_1 + b_1 + 1     \\
       2      &   & 2 &   &   & 2b_1 + 2        \\
\beta +1      &   & 1 & 1 &   & \beta + b_1 + a_2 + 1 \\
       2      &   &   & 2 &   & 2a_2 + 2        \\
       1      &   &   & 1 & 1 & a_2 + b_2 + 1     \\
       2      &   &   &   & 2 & 2b_2 + 2        \\
\gamma+1      &   &   &   & 1 & \gamma + b_2 + 1 \\
\end{array}
\]
Here $\alpha, \gamma \geq 0$ and $\beta \geq 2$, $S_i, T_i$ are the upper and lower spheres for the $i$th $\,^L_R$, and we have to find values for $a_i, b_i$ so that in each block of three the minimum value is achieved in more than one place. In the most general case, we have a whole stack of $\,^L_R$s, and we are solving this problem:
\[ 
\begin{array}{rrrrrrrrr}
 \text{Surface} & S_1 & T_1 & S_2 & T_2 & \cdots & S_N & T_N &  \text{Surface} + \Sigma(a_i S_i  + b_i T_i)\\
 \hline
 \alpha_1+1    & 1 &   &   &   &        &   &   & \alpha_1 + a_1 + 1\\
          2    & 2 &   &   &   &        &   &   & 2a_1 + 2        \\
          1    & 1 & 1 &   &   &        &   &   & a_1 + b_1 + 1     \\
          2    &   & 2 &   &   &        &   &   & 2b_1 + 2        \\
 \alpha_2+1    &   & 1 & 1 &   &        &   &   & \alpha_2 + b_1 + a_2 + 1 \\
          2    &   &   & 2 &   &        &   &   & 2a_2 + 2        \\
          1    &   &   & 1 & 1 &        &   &   & a_2 + b_2 + 1     \\
          2    &   &   &   & 2 &        &   &   & 2b_2 + 2        \\
 \alpha_3+1    &   &   &   & 1 &        &   &   & \alpha_3 + b_2 + a_3 + 1 \\
    \vdots     &   &   &   &   & \ddots &   &   & \vdots \\
 \alpha_N+1    &   &   &   &   &        & 1 &   & \alpha_N + b_{N-1}  + a_N + 1 \\
          2    &   &   &   &   &        & 2 &   & 2a_N + 2        \\
          1    &   &   &   &   &        & 1 & 1 & a_N + b_N + 1     \\
          2    &   &   &   &   &        &   & 2 & 2b_N + 2        \\
 \alpha_{N+1}+1&   &   &   &   &        &   & 1 & \alpha_{N+1} + b_N + 1 \\
\end{array}
\]
Here $\alpha_1, \alpha_{N+1} \geq 0$, all other $\alpha_i \geq 2$. We additionally increase $\alpha_1$ and $\alpha_{N+1}$ if there is a $\,^R_R$ directly above the top of the "chain of spheres" or $\,^L_L$ directly below respectively.\\

\begin{prop}
\label{prop_add_spheres}
The non-unique minimum rate condition can be satisfied by adding spheres to the Yoshida forms of surfaces in punctured torus bundles as given in section \ref{forms_of_surface}.
\end{prop} 

A key observation is that another characterisation of this "stack of $\,^L_R$s" is as a tight sub-path (see definition \ref{tight}) containing some number of $\,^L_R$s. This is unlikely to be immediately obvious, but consideration of the combinations of the relevant diagrams in Figure \ref{fig_crossings} should convince the reader. Another key observation is that we may assume that the stack above has ends, i.e. that it does not wrap around the whole punctured torus bundle and join onto itself, for then we would be in the case already dealt with in section \ref{semi-fibers}. \\

\begin{rmk}
If our incompressible surface is a semi-fiber, it turns out to be impossible to add a finite number of spheres and have the non-unique minimum rate condition hold everywhere.
\end{rmk}

We will use the following lemma, which amounts to finding a "balance point" for the weights $\alpha_i$.

\begin{lemma}
\label{balance}
For any finite sequence of positive integers $\alpha_1, \alpha_2, \ldots, \alpha_{N+1}$, we can find one of the following:
\begin{enumerate}
\item $k$ such that $\Sigma_{i=1}^k \alpha_i = \Sigma_{i=k+1}^{N+1} \alpha_i$
\item $k$ such that $\alpha_k > \left |\left (\Sigma_{i=1}^{k-1} \alpha_i \right ) - \left (\Sigma_{i=k+1}^{N+1} \alpha_i \right ) \right |$
\end{enumerate}
\end{lemma}
\begin{proof}
Suppose there is no $k$ such that the first case occurs. Then there is some $k$ for which $\Sigma_{i=1}^{k-1} \alpha_i < \alpha_k  + \Sigma_{i=k+1}^{N+1} \alpha_i$ but $\Sigma_{i=1}^{k-1} \alpha_i + \alpha_k> \Sigma_{i=k+1}^{N+1} \alpha_i$. Rearranging these two equations gives the second case.\\

Here if $k$ is at either end of the sequence we allow "empty" sums with no terms. $k$ is uniquely determined in our case, since $\alpha_i \geq 2$ for $i \neq 1, N+1$, and so moving an $\alpha_i$ from one side of the weighing scales to the other must have an effect. The only time this might not happen is if we are looking at $\alpha_1$ or $\alpha_{N+1}$ and that value is 0. This situation cannot occur at the "balance point" unless there are only two weights in the list, and they are both 0. However this falls under case 1 of the lemma.
\end{proof}

Using this lemma we can solve the general problem. We will use the solutions we found for sequences of only two weights above. Note that in the "$\alpha < \beta$" case the actual value of $\beta$ is irrelevant to the non-unique minimum rate condition holding, as long as it is large enough. So we may later add spheres which raise that value, and as long as the value of "$\alpha + \beta +1$" is greater than (or equal to) "$2\alpha+2$" when the dust settles, the non-unique minimum rate condition will hold.\\

\begin{proof}[Proof of Proposition \ref{prop_add_spheres}]
First assume we are in case 2 of the lemma. The plan of action is to use the "$\alpha < \beta$" case of the two-weight solution to work from the top of the stack and the "$\alpha > \beta$" case to work from the bottom, until we meet in the middle, at $\alpha_k$ of the lemma. Looking only at the top two $\,^L_R$ blocks, the process looks like this:
\[ 
\begin{array}{rrr}
\text{Start} & \text{Step 1} & \text{Step 2} \\
\hline
 \alpha_1+1  & 2\alpha_1+2            & 2\alpha_1+2 \\
          2  & 2\alpha_1+4            & 2\alpha_1+4 \\
          1  & 2\alpha_1+2            & 2\alpha_1+2 \\
          2  & 2\alpha_1+2            & 2\alpha_1+2 \\
 \alpha_2+1  & \alpha_1 + \alpha_2 +1 & 2(\alpha_1 + \alpha_2) + 2 \\
          2  & 2                      & 2(\alpha_1 + \alpha_2) + 4 \\
          1  & 1                      & 2(\alpha_1 + \alpha_2) + 2 \\
          2  & 2                      & 2(\alpha_1 + \alpha_2) + 2 \\
 \alpha_3+1  & \alpha_3+1             & \alpha_1 + \alpha_2 + \alpha_3 + 1 \\
     \vdots  & \vdots                 & \vdots 
\end{array}
\]
Notice that the non-unique minimum rate condition now holds for the 3rd, 4th and 5th rows here no matter the values of $\alpha_1$ and $\alpha_2$. An analogous situation occurs at the bottom of the stack, and we continue this process until we meet at $\alpha_k$. After adding all these spheres, the non-unique minimum rate condition holds everywhere apart from possibly around $\alpha_k$, at which point the situation is:
\[ 
\begin{array}{r}
\vdots \\
2(\Sigma_1^{k-1}\alpha_i)+2 \\
2(\Sigma_1^{k-1}\alpha_i)+4 \\
2(\Sigma_1^{k-1}\alpha_i)+2 \\
2(\Sigma_1^{k-1}\alpha_i)+2 \\
\Sigma_{1}^{N+1}\alpha_i +1 \\
2(\Sigma_{k+1}^{N+1}\alpha_i)+2 \\
2(\Sigma_{k+1}^{N+1}\alpha_i)+2 \\
2(\Sigma_{k+1}^{N+1}\alpha_i)+4 \\
2(\Sigma_{k+1}^{N+1}\alpha_i)+2 \\
\vdots
\end{array}
\]
The non-unique minimum rate condition holds as long as: $$\Sigma_{1}^{N+1}\alpha_i +1 \geq 2(\Sigma_1^{k-1}\alpha_i)+2$$ and $$\Sigma_{1}^{N+1}\alpha_i +1 \geq 2(\Sigma_{k+1}^{N+1}\alpha_i)+2$$
The first equation rearranges to:
$$\alpha_k \geq \left (\Sigma_{i=1}^{k-1} \alpha_i \right ) - \left (\Sigma_{i=k+1}^{N+1} \alpha_i \right ) +1$$
and the second to:
$$\alpha_k \geq  \left (\Sigma_{i=k+1}^{N+1} \alpha_i \right ) - \left ( \Sigma_{i=1}^{k-1} \alpha_i \right ) +1$$
Combined, we obtain the condition from case 2 of the lemma:
$$\alpha_k > \left |\left (\Sigma_{i=1}^{k-1} \alpha_i \right ) - \left (\Sigma_{i=k+1}^{N+1} \alpha_i \right ) \right |$$
If we are in case 1 of the lemma, we work inwards from both ends as for case 1, and then use the "$\alpha = \beta$" case of the two weight solution for the section between $\alpha_k$ and $\alpha_{k+1}$. The situation then looks like:
\[ 
\begin{array}{rrr}
\text{Start} & \text{Penultimate Step} & \text{Final Step} \\
\hline
\vdots & \vdots & \vdots \\
\alpha_{k-1}+1 & 2(\Sigma_1^{k-1}\alpha_i)+2       &  2(\Sigma_1^{k-1}\alpha_i)+2 \\
2              & 2(\Sigma_1^{k-1}\alpha_i)+4       &  2(\Sigma_1^{k-1}\alpha_i)+4 \\
1              & 2(\Sigma_1^{k-1}\alpha_i)+2       &  2(\Sigma_1^{k-1}\alpha_i)+2 \\
2              & 2(\Sigma_1^{k-1}\alpha_i)+2       &  2(\Sigma_1^{k-1}\alpha_i)+2 \\
\alpha_k + 1   & \Sigma_1^k \alpha_i +1            &  2(\Sigma_1^{k}\alpha_i) +1 \\
2              & 2                                 &  \Sigma_1^{N+1}\alpha_i + 2 \\
1              & 1                                 &  \Sigma_1^{N+1}\alpha_i + 1 \\
2              & 2                                 &  \Sigma_1^{N+1}\alpha_i + 2 \\
\alpha_{k+1}+1 & \Sigma_{k+1}^{N+1}\alpha_i +1     &  2(\Sigma_{k+1}^{N+1}\alpha_i)+1 \\
2              & 2(\Sigma_{k+2}^{N+1}\alpha_i)+2   &  2(\Sigma_{k+2}^{N+1}\alpha_i)+2 \\
1              & 2(\Sigma_{k+2}^{N+1}\alpha_i)+2   &  2(\Sigma_{k+2}^{N+1}\alpha_i)+2 \\
2              & 2(\Sigma_{k+2}^{N+1}\alpha_i)+4   &  2(\Sigma_{k+2}^{N+1}\alpha_i)+4 \\
\alpha_{k+2}+1 & 2(\Sigma_{k+2}^{N+1}\alpha_i)+2   &  2(\Sigma_{k+2}^{N+1}\alpha_i)+2 \\
\vdots & \vdots & \vdots \\
\end{array}
\]
In the middle five rows of course, $2(\Sigma_1^{k}\alpha_i) = \Sigma_1^{N+1}\alpha_i = 2(\Sigma_{k+1}^{N+1}\alpha_i)$. Looking at the non-unique minimum in the 3rd, 4th and 5th rows shown here, note that $\alpha_k \geq 2$ (if $k \neq 1$) implies that $2(\Sigma_1^{k}\alpha_i) +1 > 2(\Sigma_1^{k-1}\alpha_i)+2$ (and if $k = 1$ then there is no condition to satisfy, as the sequence starts at the $\alpha_1 + 1$ line). Similarly for the condition below.
\end{proof}

We now have Yoshida form surfaces (with added spheres) that satisfy both the "$0 - \infty$" matching condition and the non-unique minimum rate condition. We are now ready to show that the corresponding rates and directions of tetrahedra collapses actually correspond to ideal points of the tetrahedron variety. That is, we need to be able to approach a proposed set of rates and directions of degenerations with finite shapes of tetrahedra.

\begin{rmk}
A result of Kabaya~\cite{kabaya07} is able to construct these ideal points in some cases, although not if there are 'sphere' gluing equations or tetrahedra that do not degenerate at the ideal point.
\end{rmk}

\section{Tilde Equations and Solutions}
\label{tilde}
\subsection{Changing variables}
\label{changing_variables}
These are the gluing equations for part of a punctured torus bundle: a fan of Ls followed by a fan of Rs. See Figure \ref{tetrahedralisation_diag}.
\small
\[
\begin{array}{rrrr}
\lambda_*: & \hat{z}_{\hat{n}}(\frac{1}{1-\hat{v}})^2(\prod_{j=1}^{m}(\frac{1}{1-x_j})^2)(\frac{1}{1-t})^2\check{z}_1 = 1 & \rho_1: & \hat{v}(\frac{x_1-1}{x_1})^2 x_2 = 1\\
\lambda_1: & t(\frac{1}{1-z_1})^2 z_2 = 1 & \rho_2: & x_1(\frac{x_2-1}{x_2})^2 x_3 = 1\\
\lambda_2: & z_1(\frac{1}{1-z_2})^2 z_3 = 1 & \rho_3: & x_2(\frac{x_3-1}{x_3})^2 x_4 = 1\\
\lambda_3: & z_2(\frac{1}{1-z_3})^2 z_4 = 1 & \vdots & \vdots\\
\vdots & \vdots & \rho_{m-1}: & x_{m-2}(\frac{x_{m-1}-1}{x_{m-1}})^2 x_m = 1\\
\lambda_{n-1}: & z_{n-2}(\frac{1}{1-z_{n-1}})^2 z_n = 1 & \rho_m: & x_{m-1}(\frac{x_m-1}{x_m})^2 t = 1\\
\lambda_n: & z_{n-1}(\frac{1}{1-z_n})^2 \check{v} = 1 & \rho_*: & x_m(\frac{t-1}{t})^2(\prod_{j=1}^{n}(\frac{z_j-1}{z_j})^2)(\frac{\check{v}-1}{\check{v}})^2\check{x}_1 = 1
\end{array}
\]
\normalsize
If $n=0$, then that block has no $\lambda_i$ equations, and the $\lambda_*$ equation above has the $z_1$ term replaced by $\check{v}$. The $\check{\lambda}_*$ below has the $z_n$ term replaced by $t$. Similarly for if $m=0$.\\

The whole torus bundle may contain many such blocks, with different numbers of tetrahedra in each fan (so different values of $m$ and $n$). Since we are breaking down the problem into sections of the torus bundle, it is notationally convenient to \emph{not} specify which $L^{m+1}R^{n+1}$ block a particular variable or gluing equation is from. We use notation such as $\hat{z}_{\hat{n}}$ to denote a variable from the next block above the one currently in focus (and in this case the "$n$" above is denoted $\hat{n}$), or $\check{v}$ for variables below, and generally use such symbol accents whenever they are needed for clarity.\\

We know from section \ref{surfaces_isotope} that the surfaces and spheres given in figures \ref{fig_L-L_R-R} through \ref{fig_spheres} are indeed equivalent to our original Floyd-Hatcher surfaces. Fix such a surface, then running Yoshida's construction backwards we know that the orientation of twisted squares in each tetrahedron (if there are any) tells us how it is degenerating (which angle is supposedly converging to 0, which to $\infty$ and which to 1). The number of twisted squares is supposed to tell us the relative rates of degeneration.\\ 

For example, if $(z_1, z_2, z_3) \rightarrow (0, \infty, 1)$ at some rate, then as a point of Bergman's logarithmic limit set (see definition \ref{tetra_ideal_point}) we expect to see our point on $S^{3N-1}$ with $(-x, x, 0)$ in the slots corresponding to $z_1$ through $z_3$, where $x > 0$ depends on all of the other degenerations of tetrahedra.\\

Our fixed surface gives us a subset of the complex dihedral angles which are supposed to be converging to 0 (at most one angle from each tetrahedron of course, but possibly none from a given tetrahedron if that tetrahedron does not degenerate). The question is, can we find a sequence in $\mathfrak{D}(M)$ so that the appropriate tetrahedra degenerate in the appropriate ways for our surface? If so, are the relative rates of degeneration correct to give an ideal point in the logarithmic limit set which will correspond (running Yoshida's construction forward) to the surface?\\

To simplify the coming calculations, we will choose to (re)label our tetrahedralisation such that the complex dihedral angle chosen in each tetrahedron is the one that converges to $0$, assuming that tetrahedron degenerates at all. We will obtain a different $\mathfrak{T}(M;\mathcal{T})$, because our labelling of $\mathcal{T}$ has changed. \\

Although we defined $\mathfrak{T}(M)$ as an affine variety in $(\mathbb{C}\setminus\{0,1\})^N$, there is no reason to not consider trying to extend the variety over some point $p \in \mathbb{C}^N$ at which all the variables corresponding to degenerate tetrahedra are 0 (by definition, no variable is going to $\infty$ there). Assuming there is a path in our variety which actually approaches such a $p$, it should correspond to the point we hope to be able approach on the logarithmic limit set, which we will call $\bar{p}$. If there is such a $p$, $\mathfrak{T}(M)$ will be singular at this point, and of course will not correspond to a point of $\mathfrak{D}(M)$, although one could consider the corresponding added point for $\mathfrak{D}(M)$ in $(\mathbb{C}P^1)^{3N}$. \\

We need some more machinery before we can show the existence of such a $p$. We consider the equations corresponding to a variety obtained from $\mathfrak{T}(M)$ by performing a weighted blow up, the weights given to us by the supposed relative rates of degeneration of the different tetrahedra. We will denote the blow up of $\mathfrak{T}(M)$ by $\widetilde{\mathfrak{T}}(M)$\footnote{To be precise, $\widetilde{\mathfrak{T}}(M)$ is derived from the blow up of $\mathfrak{T}(M)$ in a way which we will describe later, see definition \ref{def_tilde_variety}.}. If we can extend $\mathfrak{T}(M)$ to $p$ then we should be able to extend $\widetilde{\mathfrak{T}}(M)$ to some $\tilde{p}$. It will turn out to be possible to solve the equations of $\widetilde{\mathfrak{T}}(M)$ explicitly at the blown up point, giving us $\tilde{p}$.\\

An algebraic geometry result will give us the existence of solutions to these equations near to $\tilde{p}$. A continuity argument will tell us that those nearby solutions are finite, in the sense that they correspond to points of $\mathfrak{T}(M)$ and hence to $\mathfrak{D}(M)$. Therefore $\tilde{p}$ gives us a point $\bar{p}$ in the logarithmic limit set, as opposed to being some algebraic artefact in $\widetilde{\mathfrak{T}}(M)$ (if one could only approach $\tilde{p}$ via degenerate points), and we have an ideal point. The structure of the weighted blow up in $\widetilde{\mathfrak{T}}(M)$ will give us the extra data needed to show that we have the appropriate relative rates as we approach $\tilde{p}$, and so gives us the $\bar{p}$ corresponding to our surface.\\

In the sets of equations that follow we will perform in detail the manipulations discussed above: With similar equations for other parts of the punctured torus bundle, the equations of (\ref{equations_2}) will be the polynomial equations defining $\mathfrak{T}(M)$. As we will discuss in section \ref{ideal_point}, if the punctured torus bundle has $N$ tetrahedra this variety will have $N-1$ equations in $N$ complex variables, and so will be generically 1-dimensional. The equations of (\ref{equations_3}) will be the blow up of $\mathfrak{T}(M)$, and with the extra variable $\zeta$, the blow up is generically 2-dimensional and projects down to $\mathfrak{T}(M)$. In the "tilde equations" (\ref{equations_4}), we will remove the exceptional divisor and still have a generically 2-dimensional variety. We will later add one extra normalising equation (which we will describe in section \ref{hol_semi-merid}), taking a slice of the blow up to give $\widetilde{\mathfrak{T}}(M)$, a generically 1-dimensional affine variety.\\

If we then slice again with the equation $\zeta = 0$ (the equations (\ref{equations_5})) then we expect to see a generically 0-dimensional variety. A point $\tilde{p}$ of this will correspond to a set of shapes of tetrahedra, some of which are degenerate. We will show firstly that such a point $\tilde{p}$ exists, and that this point is isolated in the slice of $\widetilde{\mathfrak{T}}(M)$ with $\zeta = 0$. Since $\widetilde{\mathfrak{T}}(M)$ is 1-dimensional however, there must be points near to $\tilde{p}$ for which $\zeta \neq 0$. Such points, at least in a neighbourhood of $\tilde{p}$, will project down to points of $\mathfrak{T}(M)$ (and hence to $\mathfrak{D}(M)$), and so give us a sequence of points approaching $\tilde{p}$, which thus corresponds to an ideal point $\bar{p}$ of $\mathfrak{D}(M)$.\\

\begin{rmk}We suspect that the stronger result that $\tilde{p}$ is a regular point of $\widetilde{\mathfrak{T}}(M)$ is generally true, although we have not proved this in all cases.
\end{rmk}

\begin{figure}[htbp]
\begin{center}
\includegraphics[width=0.9\textwidth]{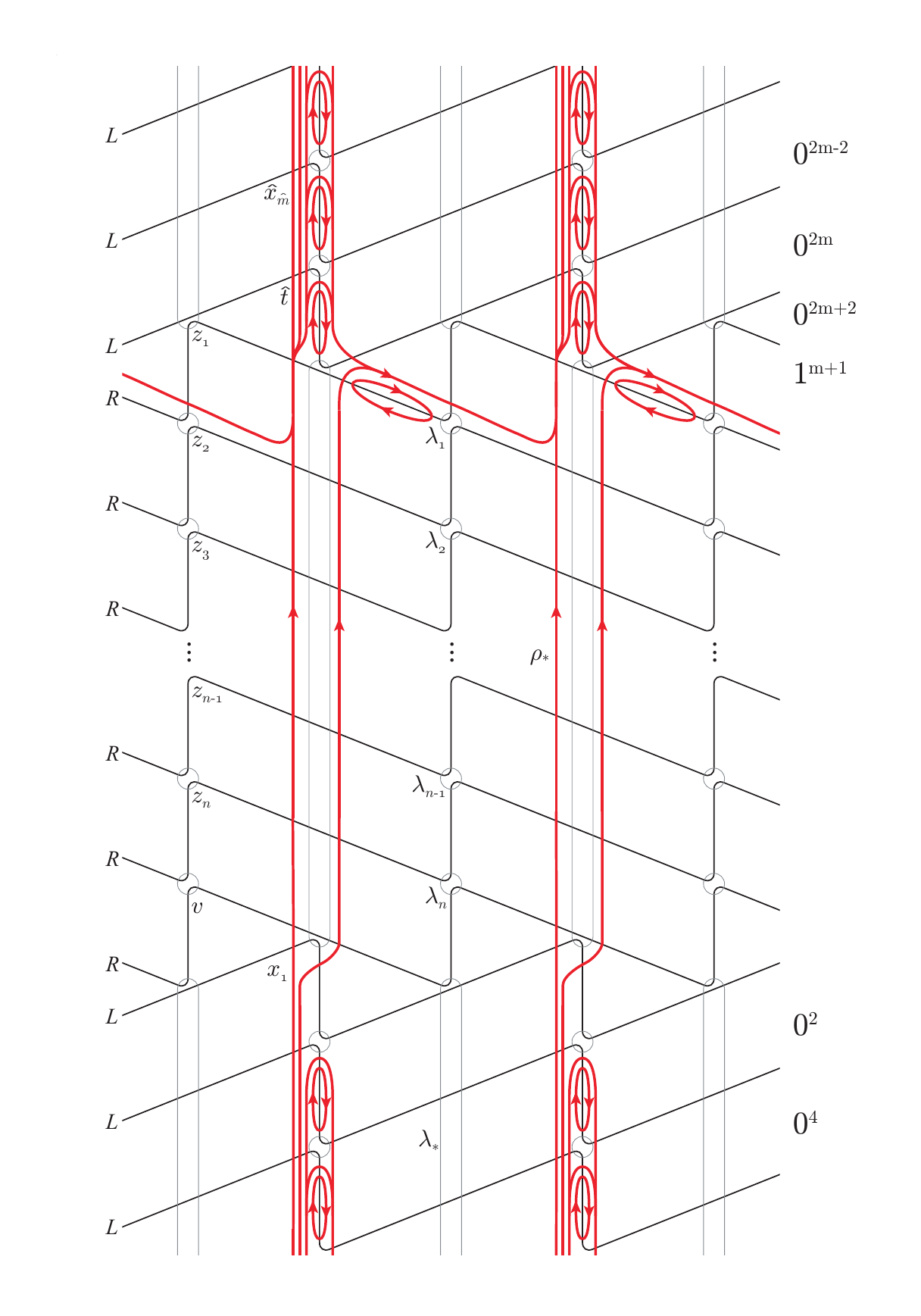}
\caption{Tetrahedra between two $\,^R_R$ sections.}
\label {R_S^S}
\end{center}
\end{figure}

So, we first make a change of variables. We will follow the manipulations with an example (see Figure \ref{R_S^S}) following the tetrahedra in between two $\,^R_R$ sections, then consider further cases.\\

First we change variables (replacing lower case with upper case) so that in each tetrahedron that is degenerating, we use the angle that is converging to 0, rather than to $\infty$ or 1. So in our example, we replace $z_1$ (which is converging to 1 at rate $m+1$) with $Z_1 = \frac{z_1 - 1}{z_1}$ (and so $z_1 = \frac{1}{1 - Z_1}$). We do not change variables corresponding to tetrahedra that are not degenerating. We will be interested here in the gluing equations $\rho_*, \lambda_1, \lambda_2, \ldots, \lambda_n$ and $\lambda_*$. After making these changes, those gluing equations look like this:

\begin{equation}
\label{equations_1}
\begin{array}{rr}

\rho_*': & \hat{X}_{\hat{m}} (\frac{ \hat{T} - 1}{ \hat{T} })^2 Z_1^2 (\prod_{j=2}^{n}(\frac{z_j-1}{z_j})^2)(\frac{v-1}{v})^2 X_1 = 1\\

\lambda_1': & \hat{T}(\frac{Z_1 - 1}{Z_1})^2 z_2 = 1\\
\lambda_2': & \frac{1}{1- Z_1}(\frac{1}{1-z_2})^2 z_3 = 1\\
\lambda_3': & z_2(\frac{1}{1-z_3})^2 z_4 = 1\\
\vdots & \vdots \\
\lambda_{n-1}': & z_{n-2}(\frac{1}{1-z_{n-1}})^2 z_n = 1\\
\lambda_n': & z_{n-1}(\frac{1}{1-z_n})^2 v = 1\\

\check{\lambda}_*': & z_n(\frac{1}{1-v})^2(\prod_{j=1}^{m}(\frac{1}{1-X_j})^2)(\frac{1}{1-\check{T}})^2\frac{1}{1-\check{Z}_1} = 1
\end{array}
\end{equation}

A shortcut for seeing what the gluing equations will look like in this form is to notice the following: Because the orientation of the boundary curve within one triangle is such that the curve goes from the $\infty$ vertex to the 0 vertex, we need only look at the directions of the arrows on parts of the curve touching that vertex. If a curve part enters a vertex of the torus bundle boundary torus from a given triangle (with variable $z$, say), then that angle of the triangle must converge to 0, and so the relevant term in the gluing equation for that vertex must be $Z$. If the curve leaves that vertex, then the angle must converge to $\infty$, and the relevant term must be $\frac{Z - 1}{Z}$. Lastly, if the curve does not enter or exit the vertex, but goes between the other two corners of the triangle, then the angle must converge to 1, and the relevant term is $\frac{1}{1-Z}$.\\

Next, we convert the equations into polynomial equations, by multiplying up by denominators and moving all terms to the left side:

\begin{equation}
\label{equations_2}
\begin{array}{rr}
\rho_*'': & \hat{X}_{\hat{m}} (\hat{T} - 1 )^2 Z_1^2 (\prod_{j=2}^{n}(z_j-1)^2)(v-1)^2 X_1 -  \hat{T}^2 (\prod_{j=2}^{n}z_j^2)v^2 = 0 \\

\lambda_1'': & \hat{T}(Z_1 - 1)^2 z_2 - Z_1^2 = 0\\
\lambda_2'': & z_3 - (1- Z_1)(1-z_2)^2  = 0\\
\lambda_3'': & z_2 z_4 -(1-z_3)^2 = 0\\
\vdots & \vdots \\
\lambda_{n-1}'': & z_{n-2} z_n -(1-z_{n-1})^2=0\\
\lambda_n'': & z_{n-1} v -(1-z_n)^2=0\\

\check{\lambda}_*'': & z_n - (1-v)^2(\prod_{j=1}^{m}(1-X_j)^2)(1-\check{T})^2(1-\check{Z}_1) = 0
\end{array}
\end{equation}

Now we introduce a new variable $\zeta$, which will be the parameter which converges to 0, and to which all other rates of convergence are relative to. If the variable $Z$ is (according to the number of twisted squares in the corresponding tetrahedron) supposed to converge to 0 at rate $k$ (i.e. there are $k$ parallel copies of the twisted square in that tetrahedron), then we set $Z = \zeta^k y$. In general we follow the same procedure for all variables, replacing the upper case letter with the lower case of the alphabetically previous letter. This is the blow up. The idea is to remove the singularity at our point by replacing variables with the "directions" ($y$ above) near the point.  \\

If we had to double up the surfaces due to non-orientability, then the only effect this has on the equations is to replace $\zeta$ with $\zeta^2$. This change will have no effect on the results of the following calculations.\\

\begin{defn} \textbf{Angle variables} are those which have not been changed in the preceding steps, that are not supposed to be going to 0, $\infty$ or 1 as we approach our supposedly ideal point. \textbf{Direction variables} are the replacements for variables that are converging to 0 (e.g. $y$ in the above example).
\end{defn}

For reference, $t, x_i, v$ and $z_j$ are angle variables and $s, w_i, u$ and $y_j$ are direction variables.\\

After these changes, the gluing equations become:
\begin{equation}
\label{equations_3}
\begin{array}{rr}

\rho_*''': & \zeta^{2\hat{m}}\hat{w}_{\hat{m}} (\zeta^{2\hat{m} + 2}\hat{s} - 1)^2 (\zeta^{\hat{m} + 1}y_1)^2 (\prod_{j=2}^{n}(z_j-1)^2)(v-1)^2 \zeta^2 w_1\\

           &  - ( \zeta^{2\hat{m} + 2}\hat{s} )^2 (\prod_{j=2}^{n}z_j^2) v^2 = 0\\
\lambda_1''': & \zeta^{2\hat{m} + 2}\hat{s}(\zeta^{\hat{m} + 1}y_1 - 1)^2 z_2 - (\zeta^{\hat{m} + 1}y_1)^2 = 0  \\
\lambda_2''': &  z_3 - (1-\zeta^{\hat{m} + 1}y_1)(1-z_2)^2 = 0 \\
\lambda_3''': & z_2 z_4 -(1-z_3)^2 = 0 \\
\vdots & \vdots \\
\lambda_{n-1}''': & z_{n-2} z_n -(1-z_{n-1})^2 = 0\\
\lambda_n''': & z_{n-1} v -(1-z_n)^2 = 0 \\

\check{\lambda}_*''': & z_n - (1-v)^2(\prod_{j=1}^{m}(1-\zeta^{2j}w_j)^2)(1-\zeta^{2m+2}\check{s})^2(1-\zeta^{m+1}\check{y}_1) = 0
\end{array}
\end{equation}
Note that in some of the equations, a power of $\zeta$ factors out. All of the above equations are of the form $A - B = 0$, and in fact, for any vertex of the torus boundary which the curve $\gamma$ passes through, the power of $\zeta$ that factors from $A$ is the number of edges of $\gamma$ entering the vertex, whereas that from $B$ is the number of edges exiting the vertex. These are of course equal. We factor out this power of $\zeta$ and delete it from our equations to obtain:

\begin{equation}
\label{equations_4}
\begin{array}{rr}
\widetilde{\rho_*}: & \hat{w}_{\hat{m}} (\zeta^{2\hat{m} + 2}\hat{s} - 1)^2 y_1^2 (\prod_{j=2}^{n}(z_j-1)^2)(v-1)^2 w_1 - \hat{s}^2 (\prod_{j=2}^{n}z_j^2) v^2 = 0\\

\widetilde{\lambda_1}: & \hat{s}(\zeta^{\hat{m} + 1}y_1 - 1)^2 z_2 - y_1^2 = 0  \\
\widetilde{\lambda_2}: &  z_3 - (1-\zeta^{\hat{m} + 1}y_1)(1-z_2)^2 = 0 \\
\widetilde{\lambda_3}: & z_2 z_4 -(1-z_3)^2 = 0 \\
\vdots & \vdots \\
\widetilde{\lambda_{n-1}}: & z_{n-2} z_n -(1-z_{n-1})^2 = 0\\
\widetilde{\lambda_n}: & z_{n-1} v -(1-z_n)^2 = 0 \\

\widetilde{\check{\lambda}_*}: & z_n - (1-v)^2(\prod_{j=1}^{m}(1-\zeta^{2j}w_j)^2)(1-\zeta^{2m+2}\check{s})^2(1-\zeta^{m+1}\check{y}_1) = 0
\end{array}
\end{equation}

Next, we slice the variety, setting $\zeta = 0$, and the tilde equations become "bar equations":

\begin{equation}
\label{equations_5}
\begin{array}{rr}

\overline{\rho_*}: & \hat{w}_{\hat{m}} y_1^2 (\prod_{j=2}^{n}(z_j-1)^2)(v-1)^2 w_1 - \hat{s}^2 (\prod_{j=2}^{n}z_j^2) v^2 = 0\\

\overline{\lambda_1}: & \hat{s} z_2 - y_1^2 = 0  \\
\overline{\lambda_2}: &  z_3 - (1-z_2)^2 = 0 \\
\overline{\lambda_3}: & z_2 z_4 -(1-z_3)^2 = 0 \\
\vdots & \vdots \\
\overline{\lambda_{n-1}}: & z_{n-2} z_n -(1-z_{n-1})^2 = 0\\
\overline{\lambda_n}: & z_{n-1} v -(1-z_n)^2 = 0 \\

\overline{\check{\lambda}_*}: & z_n - (1-v)^2 = 0
\end{array}
\end{equation}

We are looking for a solution $\tilde{p}$ to these equations for which none of the direction variables are 0, and none of the angle variables are 0, 1 or $\infty$. This will later allow us to show that nearby points are non degenerate.\\

It should now be clear why we included $\check{\lambda}_*$ in our set of equations for this part of the punctured torus bundle. The equations $\overline{\lambda_2}$ through $\overline{\lambda_n}$ then $\overline{\check{\lambda}_*}$ for the angle variables, $z_2$ through $z_n$ then $v$, form a clear pattern. If we imagine two extra variables, one at either end of the list of angle variables, which are set to have value 1, then all of these equations are of the form $a_{k-1}a_{k+1} - (1 -  a_k)^2 = 0$. \\

\subsection{Solving for angle variables}
\label{angle_variables}

\begin{lemma}
$a_k = \frac{1 - \cos{k \beta}}{1 - \cos{\beta}}$  is a solution of  $a_{k-1}a_{k+1} - (1 - a_k)^2 = 0$
\label{recursion_a}
\end{lemma}

\begin{proof}
Let $\alpha = \frac{1}{1 - \cos{\beta}}$ \big(so $1 - \frac{1}{\alpha} = \cos\beta$\big). Then:

\begin{equation*}
\begin{split}
&a_{k-1}a_{k+1} - (1 - a_k)^2 = \alpha\left(1 - \cos{\big((k-1)\beta\big)}\right)\alpha\left(1 - \cos{\big((k+1) \beta\big)}\right) - (1 - \alpha(1 - \cos{k \beta}))^2\\
&= \alpha^2\left( (1 - \cos{(k \beta - \beta)})(1 - \cos{(k \beta + \beta)}) - \left(\frac{1}{\alpha} - 1 + \cos{k \beta}\right)^2 \right)\\
&= \alpha^2\left( \big(1 - (\cos{k \beta}\cos{\beta} + \sin{k \beta}\sin{\beta})\big)\big(1 - (\cos{k \beta}\cos{\beta} - \sin{k \beta}\sin{\beta})\big) - (\cos{k \beta} - \cos{\beta})^2 \right)\\
&= \alpha^2\left( 1 - 2\cos{k \beta}\cos{\beta} + (\cos{k \beta}\cos{\beta})^2 - (\sin{k \beta}\sin{\beta})^2 - (\cos{k \beta} - \cos{\beta})^2 \right)\\
&= \alpha^2\left( 1 - 2\cos{k \beta}\cos{\beta} + (\cos{k \beta}\cos{\beta})^2 - (1 - \cos^2{k \beta})(1 - \cos^2{\beta}) - (\cos{k \beta} - \cos{\beta})^2 \right)\\
&= \alpha^2\left( 1 - 2\cos{k \beta}\cos{\beta} + (\cos{k \beta}\cos{\beta})^2 - (1 - \cos^2{k \beta} - \cos^2{\beta} + (\cos{k \beta}\cos{\beta})^2)\right.\\
& \hspace{0.15 in} \left.- (\cos{k \beta} - \cos{\beta})^2 \right)\\
&= \alpha^2\left(- 2\cos{k \beta}\cos{\beta} + (\cos^2{k \beta} + \cos^2{\beta}) - (\cos{k \beta} - \cos{\beta})^2 \right)\\
&= 0
\end{split}
\end{equation*}  
\end{proof}

In our case, we want solutions with $a_1 = 1$ and $a_{N+1} = 1$ (the two "extra" variables). The first equation is automatically true for this form of solution, and the second may be satisfied by choosing $\beta = \frac{2 \pi}{N+2}$. There are other possible choices for $\beta$ that give a solution, but this solution is easy to understand. We note the following feature of this solution for future reference:

\begin{lemma}
For the solution $a_k$ as above with $\beta = \frac{2 \pi}{N+2}$, $a_k \in \mathbb{R}$ and $a_k > 1$ for $2 \leq k \leq N$. In particular, no $a_k$ is equal to $0$.
\label{recursion_a2}
\end{lemma}

\begin{rmk}
It is worth noting as an aside that this explicit solution, generalised slightly to $a_k = \frac{1 - \cos{(k \beta + \theta})}{1 - \cos{\beta}}$ is a solution to these blocks of gluing equations independent of our looking at a degenerate point. These blocks of equations are solved as a unit by specifying $\beta, \theta \in \mathbb{C}$ (which determine what happens at either end of the fan of tetrahedra) for all points of the tetrahedron variety. It seems likely that this observation could be useful in studying torus bundles as collections of fans in contexts other than this. For instance, it is plausible that this formula should give complex versions of the "concave" sequences of angles Gu\'{e}ritaud~\cite{gueritaud} finds in fans of a torus bundle.
\end{rmk}

We record also some similar results for another sequence of equations that result from a set of tetrahedra between two $\,^L_L$ sections:
\begin{lemma}
$b_k = \frac{1 - \cos{\beta}}{1 - \cos{k \beta}}$  is a solution of  $b_{k-1}(b_k - 1)^2 b_{k+1} - b_k^2 = 0$.
\label{recursion_b}
\end{lemma}
\begin{proof}
If we set $b_k = \frac{1}{a_k}$ then the equation becomes $$ a_k^2 (\frac{1}{a_k} - 1)^2 - a_{k-1} a_{k+1}= 0$$ which is just
$$(1 - a_k)^2 - a_{k-1} a_{k+1}= 0$$
We do not worry about division by zero as all of the solutions we are interested in are positive.
\end{proof}
\begin{lemma}
For $b_1 = 1$ and $b_{N+1} = 1$ we may choose $\beta = \frac{2 \pi}{N+2}$, then  $b_k \in \mathbb{R}$ and $0 < b_k < 1$ for $2 \leq k \leq N$. 
\label{recursion_b2}
\end{lemma}

We also note the following for future use:
\begin{lemma}
The sequence of equations $a_{k-1}a_{k+1} - (1 - a_k)^2 = 0$ with $a_1 = 1$, $a_{N+1} = 1$ and  $a_k \in \mathbb{C}$ have only finitely many solutions. The same is true for the $b_k$ equations.
\label{finite_solns}
\end{lemma}
\begin{proof}
First note that having chosen a value for $a_2$, and fixing $a_1 = 1$ but leaving $a_{N+1}$ free, all values for $a_k$ are fixed, even if we extend the sequence in the obvious manner to $k > N+1$. In fact they are rational functions of $a_2$:
\begin{equation*}
\begin{split}
a_3 = &(1 - a_2)^2\\
a_4 = &\frac{(1 - a_3)^2}{a_2}\\
a_5 = &\frac{(1 - a_4)^2}{a_3}\\
&\vdots \\
a_k = &\frac{(1 - a_{k-1})^2}{a_{k-2}}\\
&\vdots
\end{split}
\end{equation*}
(If we ever had to divide by zero in this sequence then we are not at a solution to the original equations.) Thus solving the equations in the statement of this lemma is equivalent to finding solutions to $a_{N+1}(a_2) = 1$ (where we view $a_{N+1}$ as a rational function of $a_2$). Multiplying up by the denominator of this rational function we see that we have the number of possible solutions equal to the number of roots of a polynomial. The only way this can be infinite is if the polynomial is identically zero, or equivalently if $a_{N+1}(a_2)$ is identically 1. This is clearly untrue, as from Lemma \ref{recursion_a}, we have the existence of solutions to such sequences of equations with $a_{N+1} \neq 1$.\\

As we saw in Lemma \ref{recursion_b}, the equations for $b_k$ are essentially the same as those for $a_k$, and a similar argument goes through.
\end{proof}

These lemmas apply more generally in our context, in fact for all non-degenerating tetrahedra throughout the punctured torus bundle: The other possibilities for sections in Figures \ref{fig_L-L_R-R} through \ref{fig_spheres} that could surround some stack of non-degenerating  $z_j$ tetrahedra (i.e. angle variables) are $\,^L_R$ above, and/or $\,^R_L$ below instead of $\,^R_R$. Analysis of those cases (analogous to our analysis in section \ref{changing_variables}) shows that the "boundary equations" for the $z_j$ are the same as in our example, and lemma \ref{recursion_a} will apply again. Similarly, the patterns of degeneration surrounding a stack of $x_i$ variables result in the same boundary equations, and lemma \ref{recursion_b} applies for all of those cases.\\

With these observations we have found solutions for the angle variables at $\tilde{p}$.

\subsection{The holonomy of the semi-meridian}
\label{hol_semi-merid}
We now need to solve for the direction variables. As mentioned before, we need to add a normalising equation for the direction variables. To see why, consider a change of variables, setting $\zeta = a\zeta'$, and for each direction variable $y$ such that $z = \zeta^k y$, set $y' = a^ky$. Then $z = \zeta^k y = \zeta'^k a^k y =  \zeta'^k y'$. Then $\zeta'$ and the $y'$ give a different solution to the tilde equations, but one which corresponds to the same point of $\mathfrak{T}(M)$. To remove this slack, we could set one direction variable to be, say,  $1$ and solve for all the rest. However a more symmetrical and cleaner way to do things is to introduce a new variable, related to the holonomy of the meridian, or rather half of it.\\

\begin{figure}[htb]
\begin{center}
\includegraphics[width=0.9\textwidth]{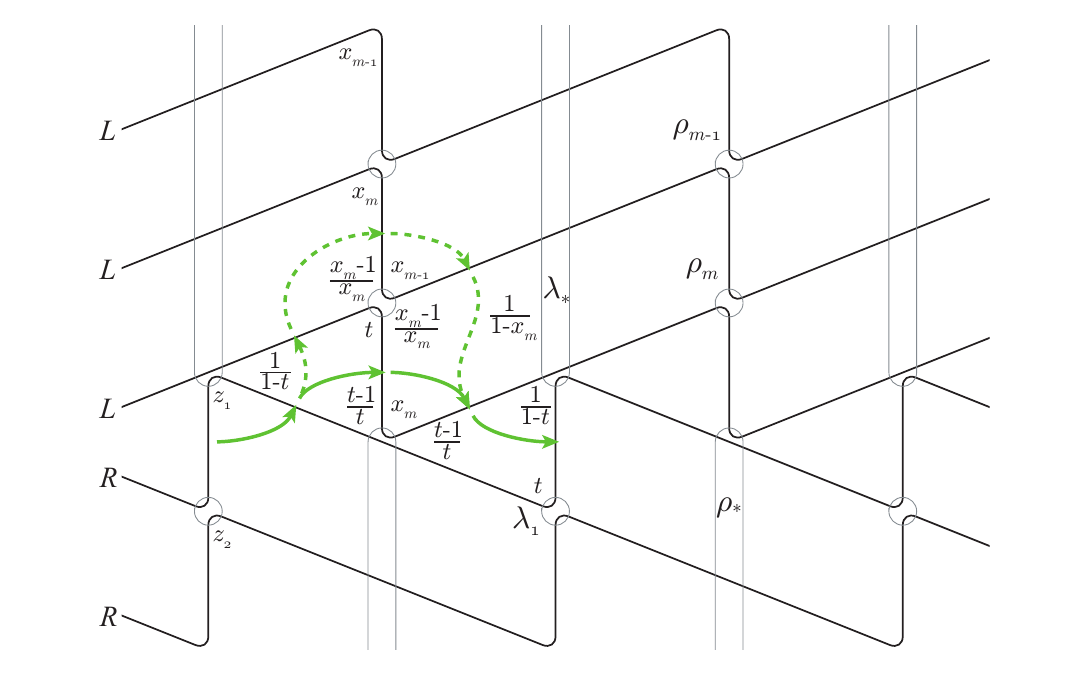}
\caption{Two ways to measure the holonomy of the semi-meridian.}
\label {fig_semimerid}
\end{center}
\end{figure}

Figure \ref{fig_semimerid} shows two ways to measure the holonomy of the \textbf{semi-meridian}. The \em meridian \em of the punctured torus bundle is the curve that wraps horizontally across all of our punctured torus bundle boundary pictures. There is a 2-fold translation symmetry in the boundary picture, and so it makes sense to talk about the semi-meridian as the curve on the quotient space of the boundary torus by that translation. We however, will only be using its holonomy as a way to simplify the algebra.\\

The holonomy of a curve on the boundary torus may be read off from the picture by taking the product of the complex angles we turn around anti-clockwise, and the inverses of the complex angles we turn around clockwise (assume for now that nothing is yet degenerate). So in Figure \ref{fig_semimerid}, the "holonomy of the semi-meridian" measured on the solid curve is: $$z_1 \frac{t}{t-1} \frac{1}{x_m} \frac{1}{1-t} = -\frac{z_1 t}{ (1-t)^2 x_m}$$

The holonomy of the meridian is of course just the square of this, but we then lose some sign information we can retain if we look at the semi-meridian. The holonomy as measured by the dotted curve on the other hand is:

$$z_1 \frac{1}{1-t} \frac{x_m}{x_m - 1} \frac{1}{x_{m-1}} \frac{1}{1-x_m} \frac{1}{1-t}\\
 = -\frac{z_1 x_m}{(1-t)^2 (x_m-1)^2x_{m-1}}$$

That these two expressions are equal is precisely the content of the gluing equation around which the two paths differ, namely $\rho_m$. The same is true in general: all measurements of the holonomy of the semi-meridian give the same answer, because they differ by the products of complex angles we see in the gluing equations. However, we want to be able to use the same relations after passing to tilde equations (and so using our new variables).\\

If we calculate the holonomy of the semi-meridian using the versions of the variables at the step just before we divide out by powers of $\zeta$ in the equations, we will (it turns out) always obtain an expression of the form $\zeta^2 f$, where $f$ is a function of angle and direction variables, and $\zeta$, but such that no power of $\zeta$ factors out further (and were we to set $\zeta = 0$, the expression for $f$ would not be zero)\footnote{If the whole surface happens to be non-orientable then when we take the double cover $\zeta^2 f$ becomes $\zeta^4 f$, and then everything goes through identically.}. The way to see this is to measure it along some version of the semi-meridian (and see that the power is indeed $2$), then show that the leading power doesn't change between two measurements of the holonomy of the semi-meridian that differ by travelling around opposite sides of a vertex. If they differed in the leading power of $\zeta$ then the terms in the gluing equation around the vertex, which are just one measurement divided by the other, would have a remaining leading $\zeta$ term. This would imply that the boundary curve of the surface entered a neighbourhood of the vertex a different number of times than it left, since each entrance contributes a $\zeta$, and each exit a $\frac{1}{\zeta}$. This is clearly impossible.\\

\begin{defn}
$$\mu := \left.\frac{\text{(holonomy of the semi-meridian)}}{\zeta^2}\right|_{\zeta = 0}$$  
This is our new variable, which we will use to normalise the direction variables in section \ref{norm_dir_vars}.
\end{defn}
There is a distinction to notice here between the equations for sphere vertices and those for non-sphere vertices. Measuring $\mu$ along paths either side of a sphere vertex result in the exact same expression of angle and direction variables. This is to do with the fact that all the complex angles around a sphere vertex go to $1$ when we set $\zeta = 0$. At all other vertices, for which the measurement of $\mu$ does change, we may effectively reconstruct the gluing equation from two "measurement equations" of $\mu$, for paths that differ by going opposite sides of the vertex. Finding a solution for these gluing equations can therefore be achieved by finding a solution to the measurement equations for $\mu$.\\

\subsection{The non-unique minimum rate condition}
\label{sec_non-unique_min}

Analysis of the gluing equation around a sphere vertex will show us the origin of the non-unique minimum rate condition that required us to add various sphere components to our incompressible surface, from Definition \ref{non-unique_min}. We illustrate with the example of the gluing equation around the vertex $\lambda_1$ in the "small" $\,^L_R$ diagram in Figure \ref{fig_crossings}. Assume that the tetrahedron below is not a hinge tetrahedron, but continues the fan, so it is assigned the variable name $z_2$. Then (with the same notation as in equations (1) through (4) from section \ref{changing_variables}):
\[
\begin{array}{rr}
\lambda_1': & \frac{1}{1 - \hat{T}}(\frac{1}{1 - Z_1})^2 \frac{1}{1 - Z_2} = 1\\
\lambda_1'': & 1 - (1 - \hat{T})(1 - Z_1)^2 (1 - Z_2) = 0\\
\lambda_1''': & 1 - (1 - \zeta\hat{s})(1 - \zeta^2 y_1)^2 (1 - \zeta y_2) = 0\\
\widetilde{\lambda_1}: & -s -y_2 + \zeta(-2y_1 + s y_2) + \zeta^2(\cdots) = 0\\
\end{array}
\]
The 1's cancel, then we remove a factor of $\zeta$. We could, but will not need to calculate the higher order terms in this equation. When we set $\zeta = 0$, we get:
\[
\overline{\lambda_1}:  -s -y_2 = 0
\]
This is a perfectly valid equation, because the powers of $\zeta$ on the different variables had a non-unique minimum. Were there a unique minimum however, then we would have reached the conclusion that some direction variable were zero, and we would apparently have the wrong degeneration rates.\\

\begin{rmk}
Satisfying the non-unique minimum rate condition is enough, in the case of punctured torus bundles, to ensure that a solution for $\tilde{p}$ exists, as we will see later. This seems unlikely to be enough more generally, as it is still possible for the equations to require that some direction variable be zero by some global combination of these local relationships. The non-unique minimum rate condition does however save us from immediate local failure.\\
\end{rmk}

The equation we eventually obtain from gluing equations about a sphere vertex is determined by which of the variables have the minimum rate of degeneration. We will be more specific about this in section \ref{tight_sections_calculations}.\\

\subsection{Normalising the direction variables}
\label{norm_dir_vars}

We now have almost enough to begin finding a solution for the direction variables (we already have the angle variables from section \ref{angle_variables}). The last ingredient is to normalise the direction variables, and we do this by setting $\mu = -1$.\\

This is an ad hoc choice for punctured torus bundles, but in this case is a good choice for a number of reasons. We could have normalised by setting one of the direction variables to be $1$ say, and solved for the other direction variables in terms of it. However the equations, like the vertices they come from, are very localised to a small number of the variables. $\mu$ on the other hand is closely related to all of the direction variables, and those relations are easily read off by "taking measurements" of $\mu$ along different paths. Additionally, this choice simplifies the behaviour of the variables in $\,^L_L$ and $\,^R_R$ sections greatly, as we will see in section \ref{solve_dir_vars}.\\

We now give a concise definition of the tilde equations and a precise definition of $\widetilde{\mathfrak{T}}(M)$.
\begin{defn}\label{def_tilde_equs}
Given a tetrahedron variety $\mathfrak{T}(M;\mathcal{T})$ corresponding to a labelled ideal tetrahedralisation of a 3-manifold $M$ and a proposed set of degeneration types and rates for the tetrahedra, we construct the {\bf tilde equations} from the defining equations of $\mathfrak{T}(M;\mathcal{T})$ as follows: 
\begin{enumerate}
\item Replace in each gluing equation every instance of a complex angle variable $z$ that is proposed to degenerate with the appropriate choice (in the new variables $\zeta$ and $\til{z}$) between $\zeta^k\til{z}, \frac{\zeta^k\til{z}-1}{\zeta^k\til{z}}$ and $\frac{1}{1-\zeta^k\til{z}}$ where $k$ is the proposed degeneration rate and the choice is given by the type of degeneration relative to the labelling of the tetrahedron.
\item Multiply up by denominators and rearrange so that the equations are of the form of a polynomial being equal to 0.
\item Remove from each polynomial equation any factor of a power of $\zeta$.
\end{enumerate}

\end{defn}

\begin{defn}
\label{def_tilde_variety}
$\widetilde{\mathfrak{T}}(M) = \widetilde{\mathfrak{T}}(M;\mathcal{T})$ is the affine variety in $\mathbb{C}^{N+2}$ ($N$ angle and direction variables, $\zeta$ and $\mu$) defined by the tilde equations, $\mu = -1$, and equations we get from measurements of $\mu$.
\end{defn}

\subsection{Solving for direction variables}
\label{solve_dir_vars}

We need to find values for all of the direction variables in the four different sections, $\,^L_L, \,^R_R, \,^R_L$ and $\,^L_R$. We will also have to deal with tight sequences separately. We consider $\,^L_L$ first, and assume that it is not part of an "extended" $\,^L_R$.\\

\subsubsection{$\,^L_L$ path section}
\label{L-L}

\begin{figure}[htb]
\begin{center}
\includegraphics[width=1.0\textwidth]{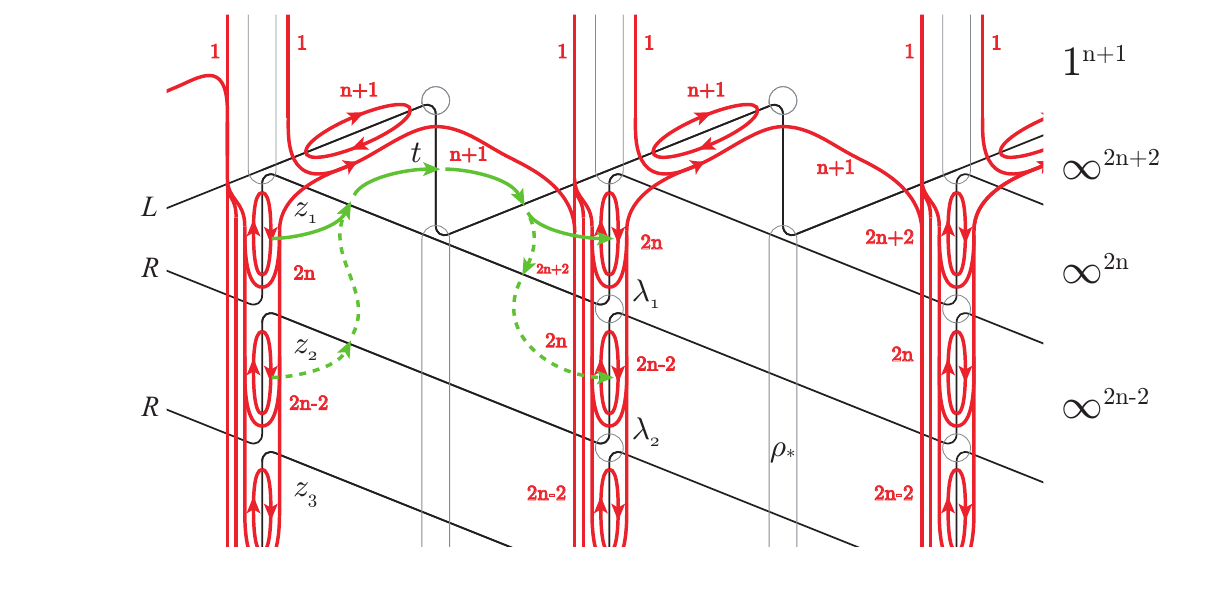}
\caption{Semi-meridians in the $\,^L_L$ section.}
\label {fig_L-L_semimerid}
\end{center}
\end{figure}

Notice first a particularly nice measurement of $\mu$ at the top of the diagram for $\,^L_L$ (see Figure \ref{fig_L-L_R-R}, and the solid path shown in Figure \ref{fig_L-L_semimerid}), passing through the $z_1$ and $t$ tetrahedra and the tetrahedron above $t$, which is not labelled as it could be either $x_m$ (part of a fan) or $\hat{v}$ (a hinge). It turns out to not matter which it is, but let us assume it is $x_m$ for now. The holonomy of the semi-meridian is:

\begin{equation*}
\begin{split}
&\left(\frac{Z_1 -1}{Z_1}\right) \left( \frac{1}{1 - T}\right) \left( \frac{1}{1 - X_m}\right) T\\
         &= \left( \frac{\zeta^{2n}y_1 -1}{\zeta^{2n}y_1}\right) \left( \frac{1}{1 - \zeta^{2n+2}s}\right) \left( \frac{1}{1 - \zeta^{n+1}w_m}\right) \zeta^{2n+2}s\\
         &= \left(\frac{\zeta^{2n}y_1 -1}{y_1}\right) \left( \frac{1}{1 - \zeta^{2n+2}s}\right) \left( \frac{1}{1 - \zeta^{n+1}w_m} \right)\zeta^{2}s\\
\end{split}
\end{equation*}

$$
-1 = \mu = \left(\frac{0 -1}{y_1}\right) \left( \frac{1}{1 - 0}\right) \left( \frac{1}{1 - 0} \right)s\\
$$
So $s = y_1$. Neither $w_m$ nor $\hat{u}$ (the direction variable for $\hat{v}$) would appear, which is why it doesn't matter which it is. A fast way to calculate $\mu$ is to follow the path along, multiplying by the direction variables for the tetrahedra we cross through, inverted if the angle is crossed in a clockwise direction and inverted and taking a minus sign whenever it crosses the boundary curve behind an orientation arrow (that is, the path crosses a corner of a triangle that the boundary curve is leaving rather than entering).\\ 

The calculation for the dotted line path in Figure \ref{fig_L-L_semimerid} is virtually the same, and gives us $y_1 = y_2$. The same is true throughout the fan, by the same argument, and so $s = y_1 = y_2 = \dots = y_n$.\\

If the $v$ tetrahedron (at the bottom of $\,^L_L$ in Figure \ref{fig_L-L_R-R}) is non-degenerate (recall that we already have non-zero solutions for non-degenerate (i.e. angle) variables) then again measuring $\mu$ by the path that loops over the $\rho_*$ vertex gives $v y_n = -1$. Thus  $s = y_1 = y_2 = \dots = y_n =-\frac{1}{v}$. If $v$ is degenerate then we have either another $\,^L_L$ or the top of a large $\,^L_R$ below. In both of these cases, measuring $\mu$ gives us simply $s = y_1 = y_2 = \dots = y_n =-1$. (*)\\

As for the variable above $t$ in the diagram (either $\hat{v}$ or $x_m$):
\begin{enumerate}
\item If it is $x_m$ then (looking at the possibilities of path section above, $\,^R_L$ or another $\,^L_L$) the variable above that must be non-degenerate (if it is $x_{m-1}$) or be degenerating to 1 (if it is $\hat{v}$). In these cases the gluing equation around the vertex there (which in this case is $\rho_m$) gives 

$$x_{m-1} w_m^2 \left( \frac{-1}{s} \right) = 1$$

(here we allow "$x_{m-1} = 1$" if it is $\hat{v}$) and so $w_m^2 = \frac{s}{x_{m-1}}$, and since we know $s$ and $x_{m-1}$, we know $w_m$ up to sign. There are no equations that involve $w_m$ other than as a square, so either root will do for our solution, and of course both will be non-zero.

\item If it is $\hat{v}$ then we look at the measurement of $\mu$ through the three tetrahedra with $\hat{v}$ in the middle. The top tetrahedron of these three is either $\hat{z}_n$ or $\hat{t}$. We will deal with the case of $\,^R_L$ above later, so assume another $\,^L_L$ is above for now. In this case $\hat{z}_n$ and $\hat{t}$ act in exactly the same way and we have:
$$-1 = \mu = (-\hat{u})\left(\frac{-1}{s}\right)u(-\hat{y}_n)$$
(or the same with $\hat{y}_n$ replaced by $\hat{s}$).
So $\hat{u}^2 = \frac{s}{\hat{y}_n}$ (or $\frac{s}{\hat{s}}$). In this case, we have already solved for the $\hat{y}_n$ or $\hat{s}$ above (see (*)), and they must be -1, so we get $\hat{u}^2 = -s$. Again there are no equations that involve $\hat{u}$ other than as a square, and so again either root will do.
\end{enumerate}

\subsubsection{$\,^R_R$ path section}
\label{R-R}
The situation for the $\,^R_R$ picture is very similar. We obtain $s = w_m = w_{m-1} = \dots = w_1 = v$, or 1 if $v$ is degenerate.\\

For the variable below $t$ (either $\check{v}$ or $z_1$):
\begin{enumerate}
\item If it is $z_1$ then the gluing equation around $\lambda_1$ gives:
$$z_2 \left(-\frac{1}{y_1}\right)^2 s = 1$$
(as above we allow "$z_2 = 1$" in the case of $\check{v}$ directly below), and so $y_1^2 = z_2 s$. Again either root will do.
\item If it is $\check{v}$, the measurement of $\mu$ (assuming $\,^R_R$ below rather than $\,^R_L$) gives:
$$-1 = \mu = \check{u}\left(\frac{1}{s}\right)(-\check{u})\check{w}_1$$ 
(similarly to above, $\check{w}_1$ could be $\check{s}$) and so $\check{u}^2 = \frac{s}{\check{w}_1}$, and again $\check{w}_1$ (or $\check{s}$) must be 1 in this case, so $\check{u}^2 = s$ and either root will do.
\end{enumerate}

\subsubsection{$\,^R_L$ path section}

First assume that the tetrahedra either side of $v$ are non-degenerate. They must then be $z_n$ above and $x_1$ below (there is no way to put a hinge tetrahedron in one of those spots and not have it degenerate). Measuring $\mu$ here gives us $-1 = \mu = (-u)x_1 u \left(\frac{1}{z_n}\right)$, so $u^2 = \frac{z_n}{x_1}$, and since the angle variables are all known and non-zero, we obtain the value of $u$, up to sign. As before, the sign doesn't matter. \\

Now if one or both of the tetrahedra either side of $v$ are degenerate because of a $\,^R_R$ above or $\,^L_L$ below we have the same equation, with $x_1$ replaced by $-\frac{1}{\check{s}}$ and/or $\frac{1}{z_n}$ replaced by $\frac{1}{\hat{s}}$. Assuming the $\,^R_R$ above or $\,^L_L$ below are not part of extended $\,^L_R$s, we already have non-zero solutions for $\hat{s}, \check{s}$ and so as before, we are done and the sign doesn't matter.

\subsubsection{$\,^L_R$ path section, spheres, tight sections}
\label{tight_sections_calculations}
\begin{figure}[htb]
\begin{center}
\includegraphics[width=0.9\textwidth]{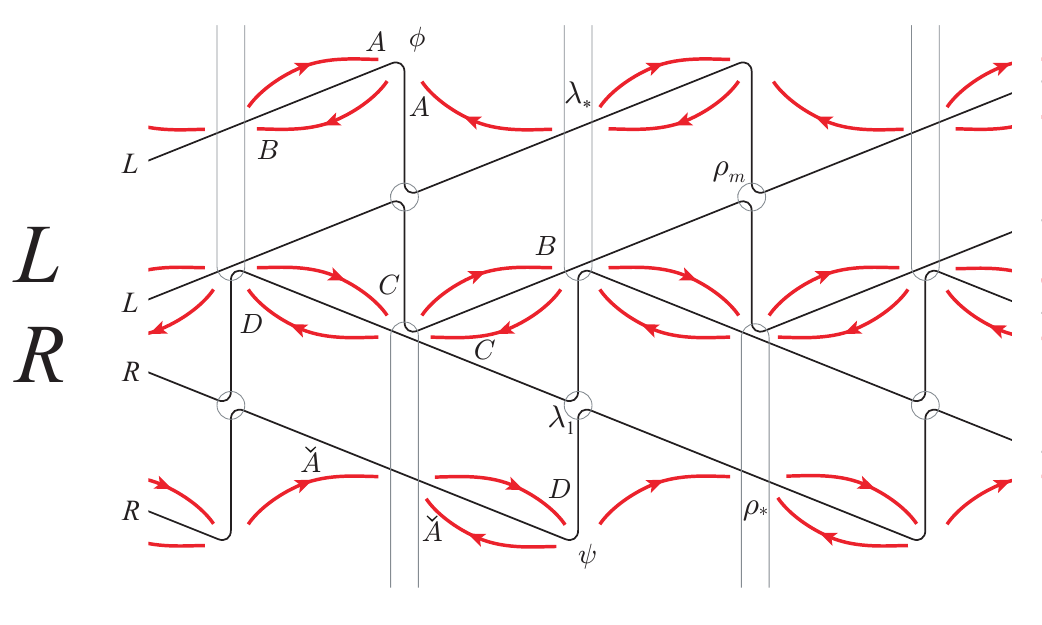}
\caption{Simplified boundary picture of the small $\,^L_R$ section (ends of thin strips not shown).}
\label {fig_sphere_chain_small}
\end{center}
\end{figure}

\begin{figure}[htbp]
\begin{center}
\includegraphics[width=0.9\textwidth]{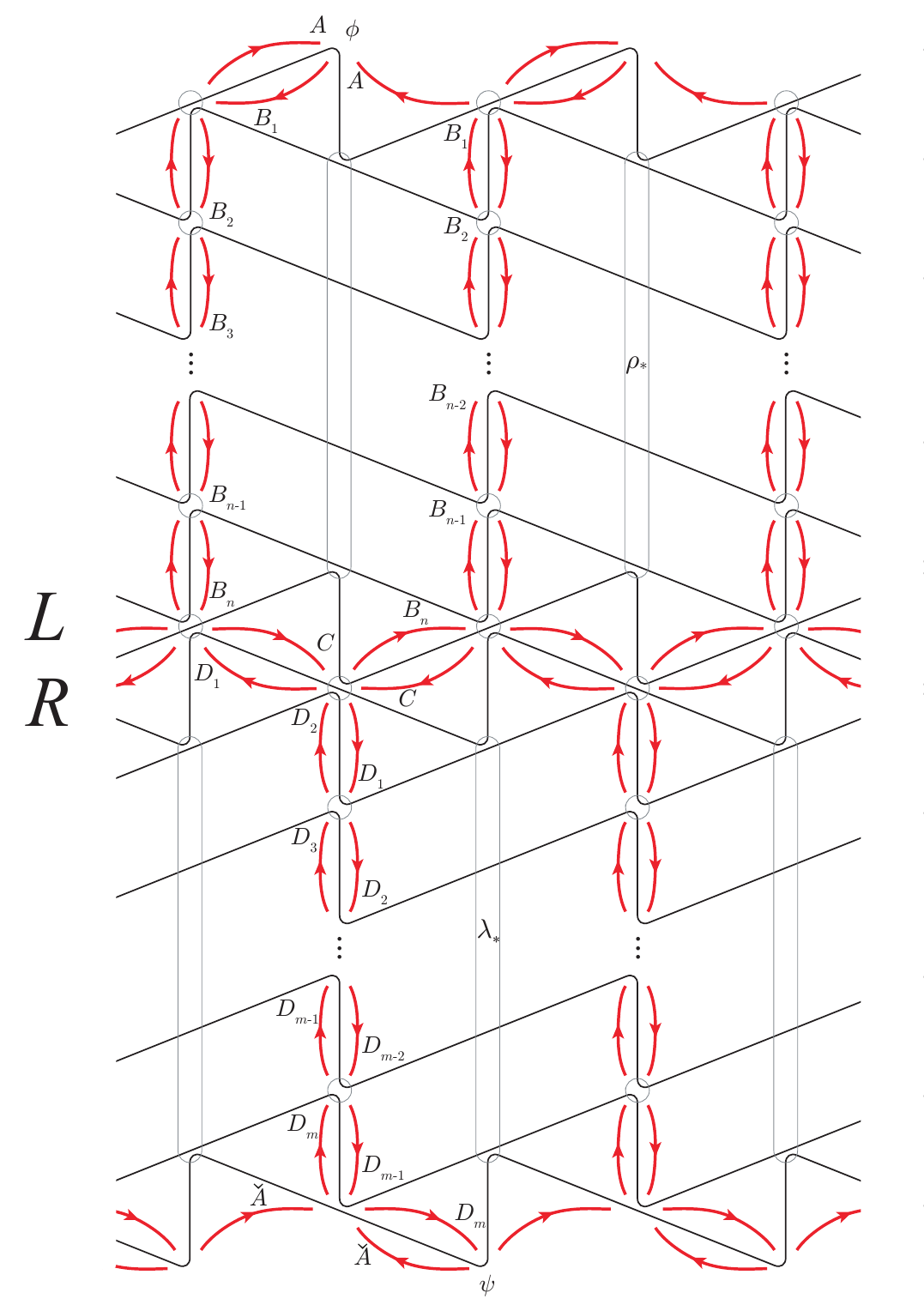}
\caption{Simplified boundary picture of the extended $\,^L_R$ section (ends of thin strips not shown).}
\label {fig_sphere_chain_large}
\end{center}
\end{figure}

In Figures \ref{fig_sphere_chain_small} and \ref{fig_sphere_chain_large} we see in simplified form the small and large versions of an $\,^L_R$ path section, we assume with the appropriate numbers of spheres added, and with $\,^R_R$ and $\,^L_L$ sections added to make what we have been calling an "extended" $\,^L_R$. We have changed the labelling of the tetrahedra to correspond more closely with the behaviour of the surface rather than the punctured torus bundle. In particular $A$ and $\check{A}$ may or may not be hinge tetrahedra. The labelled complex angles in the tetrahedra are all those degenerating to $0$. The $\phi$ and $\psi$ we will use later to refer to the product of terms contributing to the equation around the vertex which are not otherwise labelled on the diagram.\\

Consider the effect that adding spheres has to the equations when we change variables and set $\zeta = 0$. No tetrahedron that was non degenerate before is degenerate after, all that changes are the rates at which tetrahedra degenerate. The corresponding powers of $\zeta$ are larger as a result, but then a larger power is factored out (for gluing equations around non-sphere vertices). When we set $\zeta = 0$, all trace of the added sphere is gone, apart from the effect it can have on the sphere vertices. We have already seen why the non-unique minimum rate condition is required, but the exact way in which it is satisfied determines the equation we obtain from the sphere vertices.\\

Recall the solution to the problem of adding spheres to a single $\,^L_R$ from section \ref{adding_spheres}, to which we have added two more columns and eliminated some variables when we can for clarity,  and labelled with our new variable names:
\[
\begin{array}{crrrrr}
      & \alpha + 1< \beta  & \alpha+1 = \beta    & \alpha = \beta    & \alpha = \beta+1  & \alpha > \beta +1 \\
\hline
A     & 2\alpha+2          & 2\alpha+2           & 2\alpha +1        & 2\beta +2           &  \alpha + \beta +1 \\
B_n   & 2\alpha+4          & 2\alpha+4           & 2\alpha + 2       & 2\beta+2            & 2\beta+2          \\
C     & 2\alpha+2          & 2\alpha+2           & 2\alpha  +1       & 2\beta+2            & 2\beta+2          \\
D_1   & 2\alpha+2          & 2\alpha+2           & 2\alpha + 2       & 2\beta+4            & 2\beta+4          \\
\check{A}& \alpha + \beta +1  & 2\alpha+2           & 2\alpha +1        & 2\beta+2            & 2\beta+2
\end{array}
\]
We will continue to refer to these five types as "$\alpha + 1 < \beta$" etc., although the connection to the numbers "$\alpha$" and "$\beta$" is rather tenuous at this point. When we pass to tilde equations, and then the bar equations (setting $\zeta = 0$), we get the following equations from the gluing equations around sphere vertices. Here we break our "alphabetically previous" convention on variable names, and simply set $A = \zeta^k a$ and so on:

\[
\begin{array}{rrrrr}
\alpha + 1< \beta  & \alpha+1 = \beta    & \alpha = \beta    & \alpha = \beta+1  & \alpha > \beta +1 \\
\hline
-a -c = 0          &  -a -c = 0          &  -a -c = 0        & -a -2b_n - c  = 0 & -2b_n - c = 0     \\
-c -2d_1 = 0       &-c-2d_1-\check{a} = 0&  -c -\check{a} = 0&  -c -\check{a} = 0&  -c -\check{a} = 0
\end{array}
\]
One of these five possibilities occurs around each $\,^L_R$, and which one occurs we determined in the proof of Proposition \ref{prop_add_spheres}. Note from that proof that we are always either in the "$\alpha + 1< \beta$" case (working from above) or the "$ \alpha > \beta +1$" case (working from below) apart possibly from where we meet in the middle, at which any of the five possibilities can happen.\\

The differences between the five cases are expressed only in the equations coming from sphere vertices. The following equations hold in all cases:

\renewcommand{\labelenumi}{\roman{enumi})}

\begin{enumerate}
\item $b_1 = b_2 = \ldots = b_n$ (and hence we suppress the subscripts from now on).
\item $d_1 = d_2 = \ldots = d_m$ (ditto).
\item $c^2 = - b d$.
\item $b = -a^2 \phi$.
\item $d = \frac{\check{a}^2}{\psi}$.

\end{enumerate}
i) and ii) come from the same calculations as were made for $\,^L_L$ and $\,^R_R$. iii) is the measurement of $\mu$ through $b_n, c, d_1$. iv) and v) come from the gluing equations in the diagrams labelled with $\phi$ and $\psi$.\\

Note that seen from this perspective, the "small" and "large" versions of $\,^L_R$ are part of the same inclusive scheme. We now solve for the variables in the different cases, starting with the "$ \alpha + 1 < \beta$" case, working in from above. Assume for now that we know the value of $\phi$.\\

\paragraph{The "$ \alpha + 1 < \beta$" case.}
$b = -a^2 \phi$ and $d = -\frac{c}{2} = \frac{a}{2}$, so $c^2 = a^2 = -bd = a^2 \phi \frac{a}{2}$. We are interested in solutions for $a \neq 0$, so $1 = \frac{\phi a}{2}$ and $a = \frac{2}{\phi}$. This also gives us values for $b, c$ and $d$ in terms of $\phi$: $b = \frac{-4}{\phi}$, $c = -\frac{2}{\phi}$ and $d = \frac{1}{\phi}$. Assuming that the chain of spheres continues below with another (possibly extended) $\,^L_R$ section (and so the $\check{A}$ in this section is the $A$ for the section below), then we can also calculate $\check{\phi}$ (the $\phi$ for the section below). If the $\,^L_R$ is large on its lower half, we have $\check{\phi} = \frac{d_{m-1}}{(-d_m)^2} = \frac{1}{d} = \phi$. If it is small on its lower half, we get the same result: $\check{\phi} = \frac{c^2}{(-d)^2(-b)} = \frac{-bd}{-bd^2} =  \frac{1}{d} =\phi$.\\

\paragraph{The "$ \alpha > \beta + 1$" case.}
Working from below upwards with the "$\alpha > \beta +1$" case is very similar. We assume we know the value of $\psi$ and obtain: $\check{a} = -2\psi$, $d = 4\psi$, $c =  2\psi$, $b = -\psi$ and $\hat{\psi} = \psi$.\\

We work in from both sides, and eventually arrive in the middle. If we do not have one of the three cases $\alpha+1 = \beta, \alpha = \beta$ or $ \alpha = \beta+1$ in the middle, we apply the appropriate same procedure as above and have one final value to determine: the $a$ sitting at the "balance point". We know $\phi$ and $\psi$ throughout, and have a gluing equation from which we get the equation $\phi a^2 \psi = 1$. Thus $a = \pm \frac{1}{\sqrt{\phi \psi}}$. It doesn't matter which sign we choose.\\

\paragraph{The "$ \alpha = \beta$" case.}
In this case: $a = -c = \check{a}$, $b = -a^2 \phi$, $d = \frac{\check{a}^2}{\psi}$, $c^2 = -bd = \frac{a^2 \phi \check{a}^2}{\psi} = c^4 \frac{\phi}{\psi}$. Again we are looking for non-zero solutions, so we can divide out to get $1 = c^2 \frac{\phi}{\psi}$ so $c = \pm\sqrt{\frac{\psi}{\phi}} = -a = -\check{a}$ and $b = -\frac{\psi}{\phi}\phi = -\psi$, $d = \frac{\psi}{\phi}\frac{1}{\psi} = \frac{1}{\phi}$.\\

\paragraph{The "$ \alpha + 1 =  \beta$" case.}
The final two cases are a little more complicated, as for the first time we have an equation involving three terms added. For the $ \alpha+1 = \beta $ case:   $b = -a^2 \phi$, $d = \frac{\check{a}^2}{\psi}$, $c^2 = -bd = \frac{a^2 \phi \check{a}^2}{\psi}$. $-a = c$ so we have $1 = \frac{\phi}{\psi}\check{a}^2$ and so $\check{a} = \pm\sqrt{\frac{\psi}{\phi}}$ and $d = \frac{1}{\phi}$. However $- a = c = -2d - \check{a} = -\check{a}(2\frac{\check{a}}{\psi} + 1) = \mp \sqrt{\frac{\psi}{\phi}}\left(\pm 2\frac{\sqrt{\frac{\psi}{\phi}}}{\psi} + 1\right) = \mp \sqrt{\frac{\psi}{\phi}}\left(1 \pm \frac{2}{\sqrt{\phi \psi}} \right)$, and $b = -a^2\phi = -\frac{\psi}{\phi}\left( 1 \pm \frac{4}{\sqrt{\phi \psi}} + \frac{4}{\phi \psi}   \right)$.\\

We should be concerned now, that it is possible to get a direction variable being $0$ if we choose the wrong sign and have $\pm \frac{2}{\sqrt{\phi \psi}} = -1$. Of course, we can just choose the other sign if one causes trouble.\\

\paragraph{The "$ \alpha = \beta + 1$" case.}
This is similar, we obtain $a = \pm \sqrt{\frac{\psi}{\phi}}, b = -\psi, -\check{a} = c = \mp\sqrt{\frac{\psi}{\phi}}\left( 1 \mp 2\sqrt{\phi \psi} \right), d = \frac{\psi}{\phi} \left( 1 \mp 4\sqrt{\phi\psi} + 4\phi\psi \right)$.\\

It remains to calculate the values of $\phi$ and $\psi$. 
\subsubsection{$\phi$ and $\psi$.}

We look first at $\psi$, at the bottom end of a chain of spheres:\\

If the $\check{A}$ tetrahedron is not a hinge tetrahedron, then the vertex at which $\psi$ sits is 4-valent, and $\psi$ is either an angle variable, or possibly 1 (if the $\psi$ tetrahedron is a hinge tetrahedron). In either case the value is determined and non-zero. If the $\check{A}$ tetrahedron is a hinge tetrahedron we consider the three possible cases for what path sections are below this (possibly extended) $\,^L_R$: 
\begin{enumerate}
\item $\,^R_R$
\item $\,^R_L$ then $\,^L_L$
\item $\,^R_L$ then $\,^L_R$
\end{enumerate}

\begin{figure}[htbp]
\begin{center}
\includegraphics[width=0.9\textwidth]{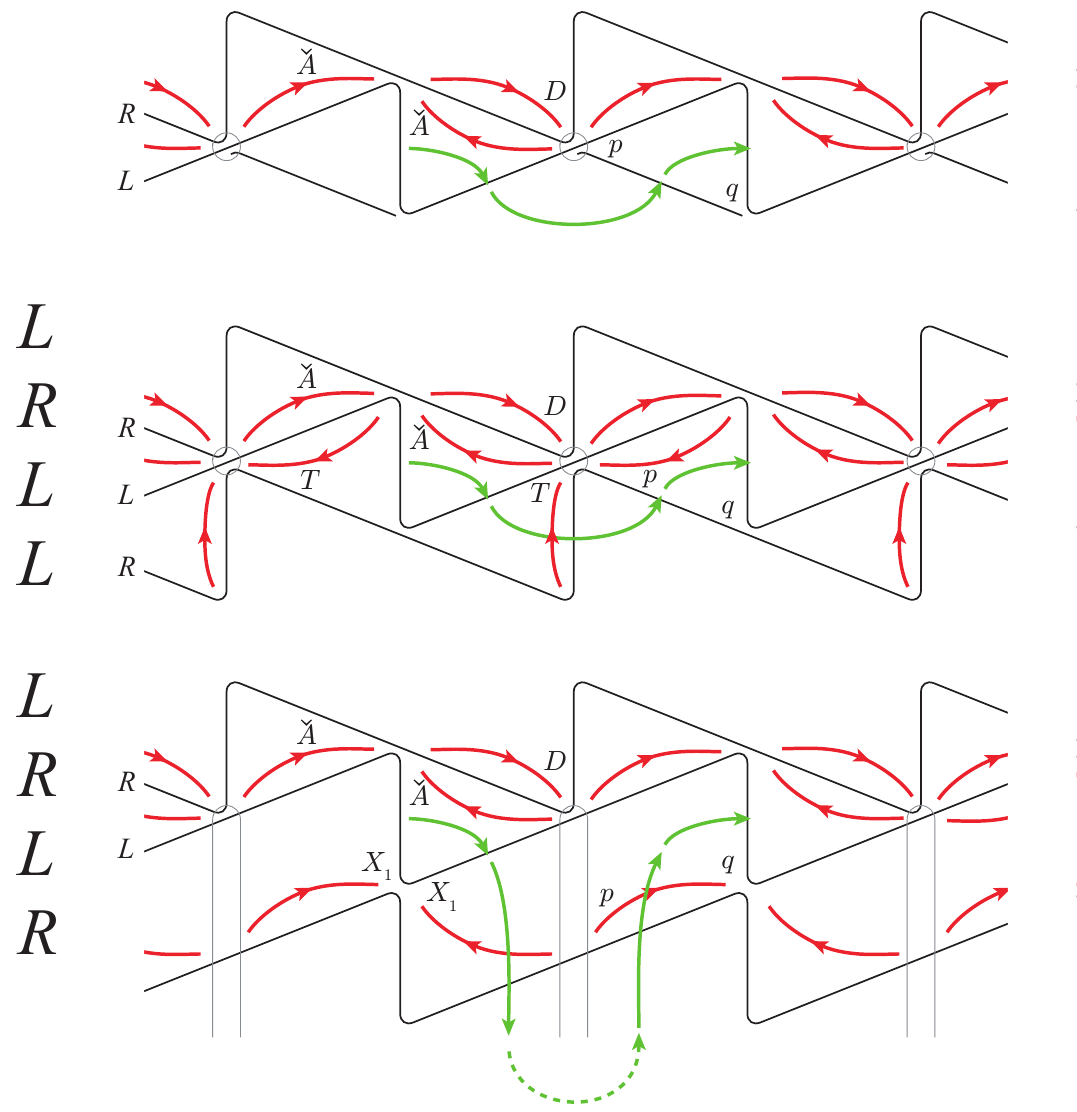}
\caption{The general situation and two of the possibilities at the bottom of a chain of spheres.}
\label {fig_below_chain}
\end{center}
\end{figure}

If we have an $\,^R_R$ below, all the complex angles that multiply to form $\psi$ are 1 (see Figures \ref{fig_L-L_R-R} and \ref{fig_crossings}), and so $\psi = 1$. In the other two cases the $\,^R_L$ simply adds one to the degeneration rate at $\check{A}$. We show the situations in Figure \ref{fig_below_chain}.

Note that the semi-meridian shown is very nearly covering the same angles as what we want, $\psi$. In fact one can see that $\mu = \frac{\psi}{p q}$ ($p$ and $q$ are anti-clockwise "angles" as usual). Since $\mu = -1$ we have that $\psi = -pq$.\\

In the $\,^R_L$ then $\,^L_L$ case, we get $p = s$ (recall $s$ is the direction variable for $T$) and $q = 1$ so $\psi = -s$. We know the value of $s$ from section \ref{L-L}, noting that this $\,^L_L$ cannot be part of an extended $\,^L_R$ (as required in section \ref{L-L}) since if it were, we would not be at the bottom of the chain of spheres.\\

For the $\,^R_L$ then $\,^L_R$ case, Figure \ref{fig_below_chain} shows the extreme case of the top of the sphere above the $\,^L_R$ right next to the $\check{A}$ tetrahedron. It cannot be any higher (i.e. overlap with $\check{A}$) for that would again mean that we are not at the bottom of the chain of spheres. It can be lower however and we would have some angle variables between the sphere and $\check{A}$. In this case the tetrahedron with $p$ and $q$ marked is labelled $x_1$ at the uppermost vertex, $p = \frac{1}{1-x_1}, q = \frac{x_1-1}{x_1}$, and $\phi = \frac{1}{x_1}$. We of course already have non-zero solutions for the angle variables. For the extreme situation as in Figure \ref{fig_below_chain} we have $p = -\frac{1}{w_1}$ (recall $w_1$ is the direction variable for $x_1$), $q = w_1$ and so $\psi = 1$.\\

We have covered all cases for calculating $\psi$ at the bottom of the chain of spheres. The calculations for $\phi$ at the top of the chain of spheres are very similar:\\

If $A$ is not a hinge tetrahedron then $\phi$ is either an angle variable or $1$. If $A$ is a hinge tetrahedron then one of three possibilities can happen above:
\begin{enumerate}
\item $\,^L_L$
\item $\,^R_R$ above $\,^R_L$
\item $\,^L_R$ above $\,^R_L$
\end{enumerate}

If we have a $\,^L_L$ then as for the corresponding case at the bottom of the chain, $\phi = 1$.\\

\begin{figure}[htbp]
\begin{center}
\includegraphics[width=0.9\textwidth]{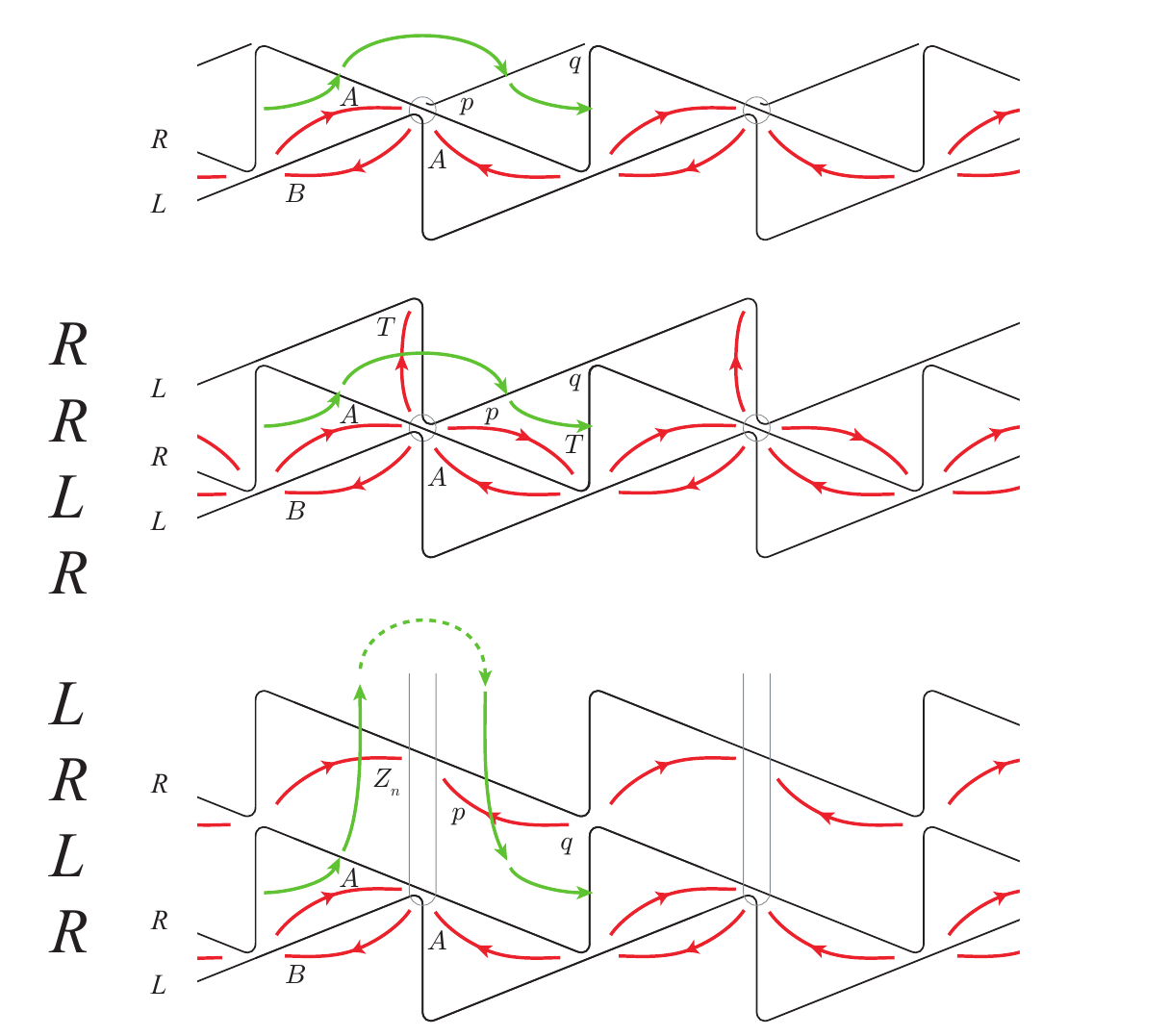}
\caption{The general situation and two of the possibilities at the top of a chain of spheres.}
\label {fig_above_chain}
\end{center}
\end{figure}

We show the situations in Figure \ref{fig_above_chain}. This time we get $-1 = \frac{1}{\mu} = \frac{\phi}{p q}$ so $\phi = - p q$.\\

In the $\,^R_R$ above $\,^R_L$ case, we get $p = -\frac{1}{s}$ and $q = 1$ so $\phi = \frac{1}{s}$. We know the value of $s$ from section \ref{R-R}, and again this $\,^R_R$ cannot be part of an extended $\,^L_R$ since if it were, we would not be at the top of the chain of spheres.\\

For the $\,^L_R$ above $\,^R_L$ case we again get in the extreme case of Figure \ref{fig_above_chain} $\phi = - y_n \left(-\frac{1}{y_n}\right) = 1$. If $Z_n$ is a non-degenerating tetrahedron we use the variable $z_n$, the complex angle at the top of the tetrahedron and get $\phi = - \frac{z_n-1}{z_n} \frac{1}{1-z_n} = \frac{1}{z_n}$.  \\

\subsection{$\tilde{p}$ corresponds to an Ideal Point}
\label{ideal_point}

We have now determined values, up to sign in some cases, for all angle and all direction variables. These values, plus $\zeta = 0$ solve the tilde equations by construction. Whenever we had a choice of sign, either option gives a solution to the equations for $\zeta = 0$.\\ 

We have the existence of a point of $\widetilde{\mathfrak{T}}(M)$, the variety defined by the tilde equations and the equation $\mu = -1$, with the extra condition that $\zeta = 0$. We now want to show that there are other points of the variety nearby, and moreover that we have nearby points that correspond to (finite) points of the variety $\mathfrak{T}(M)$ and hence $\mathfrak{D}(M)$. The following discussion and result prove the first part:\\

Suppose we have $N$ tetrahedra in our torus bundle. The torus bundle is made up of some number of $L^{m+1}R^{n+1}$ sections, and $N$ is the sum of all of those $m+1$s and $n+1$s. We begin with $N$ equations (the gluing equations), and the $N$ variables (the original complex angles). One of the gluing equations is dependent on the other $N-1$. This is a standard fact for tetrahedralisations of 3-manifolds with a single boundary component (starting with the original gluing equations, multiply $N-1$ of them together, and use the identities between the 3 angles in each tetrahedron to obtain the $N$th). Thus we can remove one gluing equation and now have $N-1$ equations in $N$ variables. Next we convert all these to tilde equations, and add a variable, $\zeta$. We then add the variable $\mu$, add the equation $\mu = -1$, and effectively add one more equation, a "measurement" of $\mu$ in terms of some direction (and possibly angle) variables. All other measurements of $\mu$ we use are derived from this measurement and the tilde equations. This brings us to $N+1$ equations in $N+2$ variables. \\

\begin{prop}
If $x \in \mathbb{C}^{K+1}$ satisfies polynomial equations: $$f_1, f_2, \ldots, f_{K} \in \mathbb{C}[x_1, x_2, \ldots x_{K+1}]$$ then there exist other solutions to these equations arbitrarily close to $x$.
\end{prop}
The heuristic reason for this is that starting from $\mathbb{C}^{K+1}$, every polynomial we add to our set of equations cuts down the dimension of the set of solutions by at most one (unless it results in an inconsistent set of equations). Since we only make $K$ cuts, and started with $K+1$ dimensions, we will have at least one left by the end. The existence of $x$ demonstrates the consistency of the equations.\\

Here is a more formal proof:

\begin{proof}
By corollary 3 of section I.\S 7 of Mumford \cite{mumford} (page 44), the codimension of any component of the variety defined by the $K$ equations is less than or equal to $K$. Recall the fact that irreducible varieties over an algebraically closed field are connected. Suppose for contradiction that $x$ is isolated. If this is so, then ${x}$ consists of an entire irreducible component which therefore has dimension 0 and so codimension $K+1$. This contradicts the fact that the codimension is less than or equal to $K$.
\end{proof}

We have shown the existence of a point $\tilde{p} \in \widetilde{\mathfrak{T}}(M)$ with $\zeta = 0$. We now also know that $\widetilde{\mathfrak{T}}(M)$ contains points arbitrarily near to $\tilde{p}$. Moreover, we know something about what such a nearby point looks like:
\begin{prop}
There exist points of $\widetilde{\mathfrak{T}}(M)$ arbitrarily near $\tilde{p}$ which are finite (when we convert them back to points of $\mathfrak{T}(M)$, no angle is 0, $\infty$ or 1).
\label{p_ideal1}
\end{prop}
\begin{proof}
We want to show first that for points near enough to $\tilde{p}$, $\zeta \neq 0$. In other words, that $\tilde{p}$ is isolated among elements of $\widetilde{\mathfrak{T}}(M)$ with $\zeta = 0$. In order to show this, we consider the steps we took to find $\tilde{p}$ (i.e. a solution with $\zeta = 0$). We first set $\zeta = 0$, then chose among finitely many solutions for each sequence of angle variables (see Lemma \ref{finite_solns}). If we look for points with $\zeta = 0$ and near enough to $\tilde{p}$, then the choices of angle variables must be the same as for $\tilde{p}$, since there are only finitely many such choices, and any two choices will have some distance between them. Then the only choices we had for direction variables were some signs. Choosing a different sign again puts us at some distance from $\tilde{p}$ and so when we look near enough to $\tilde{p}$, the only solution to the equations with $\zeta = 0$ is $\tilde{p}$ itself. Therefore we must have points nearby for which $\zeta \neq 0$. \\

We should also consider if at any point an assumption we made about not dividing by zero when finding a solution to $\tilde{p}$ could be false now that we are interested in any solution with $\zeta = 0$. If however such a solution does exist, it is not near to $\tilde{p}$, since our solution $\tilde{p}$ has no variable near 0. Thus we can ignore these possible solutions when trying to find solutions near $\tilde{p}$.\\

Suppose $\tilde{q}$ is a point of $\widetilde{\mathfrak{T}}(M)$ near $\tilde{p}$ for which $\zeta$ is near to but not equal to 0. Then by continuity, we can ensure that for $\tilde{q}$ all angle variables are bounded away from $0$, $\infty$ or $1$, since they are so for $\tilde{p}$. When we change variables back to the original angles of the original gluing equations (i.e. consider the $q \in \mathfrak{T}(M)$ corresponding to $\tilde{q} \in \widetilde{\mathfrak{T}}(M)$), all direction variables become $\zeta^ky$ for some $k>0$ and $y$ a direction variable. By continuity, $y$ is near whatever value it had at $\tilde{p}$, that is, bounded away from 0 and $\infty$. For $\tilde{q}$ sufficiently close to (but not equal to) $\tilde{p}$ then, $0<|\zeta^ky|<1$, and this angle is also finite.\\
\end{proof}
We can therefore construct a sequence of finite points of $\widetilde{\mathfrak{T}}(M)$ which converge to $\tilde{p}$ as $\zeta \rightarrow 0$. Now consider the points corresponding to this sequence in $B^{3N}$, as in definition \ref{tetra_ideal_point}. For each complex dihedral angle\footnote{Here the subscript indexes the tetrahedra rather than the 3 dihedral angles of a single tetrahedron, unlike in definition \ref{tetra_ideal_point}.} $z_j$ which converges to 0 as we approach $\tilde{p}$, 
$$z_j = \zeta^{k_j} y_j$$ 
for some integer $k_j$ and corresponding direction variable $y_j$. $y_j$ approaches some finite non-zero value (the value it attains at $\tilde{p}$), and so 
$$\log|z_j| = \log|\zeta^{k_j} y_j| = k_j \log|\zeta| + \log|y_j|$$ 
approaches negative infinity at speed $k_j$. The complex dihedral angles that converge to $\infty$ are all of the form $\frac{z_j - 1}{z_j}$, and:
$$\log\left|\frac{z_j - 1}{z_j}\right| = \log|z_j - 1| - \log|z_j|$$
The first term converges to $\log|1| = 0$, and the second to positive infinity, again at rate $k_j$. When we divide by the denominator as in definition \ref{tetra_ideal_point}, we scale our sequence to be within the unit ball $B^{3N}$, and the sequence converges to a point $\bar{p}$ on the boundary sphere $S^{3N-1}$, the exact location determined only by the relative rates $k_j$, and the directions of collapse within each tetrahedron.\\

Very similar arguments to those in this section give us the following theorem (the only difference is that our choice of $\mu = -1$ and a measurement of $\mu$ to normalise rather than a single equation give us one more variable and one more equation):

\begin{thm}\label{soln_of_tilde_gives_ideal_pt}
Given a 3-manifold $M$ with torus boundary and with ideal triangulation $\mathcal{T}$, and a proposed set of degeneration types and rates for the tetrahedra, if the tilde equations corresponding to this degeneration together with a normalising equation have a solution with $\zeta = 0$, all angle variables non-degenerate and all direction variables non-zero and that solution is an isolated point then the solution corresponds to an ideal point of the deformation variety $\mathfrak{D}(M;\mathcal{T})$.
\end{thm}

We have shown:
\begin{thm}\label{yoshida_from_defm}
For each incompressible surface of a punctured torus bundle that is not the fiber or a semi-fiber, we can find an ideal point $\tilde{p}$ of the deformation variety which corresponds to the incompressible surface under Yoshida's construction.
\end{thm}

Finally, in order to prove theorem \ref{main_result} we apply this result:

\begin{numless_thm}[theorem 5.2 of \cite{segerman_twsq_spnn}]
Let $M$ be an oriented 3-manifold with $\bdry M$ a union of tori with ideal triangulation $\mathcal{T}$ and $T$ a two-sided twisted squares surface obtained via Yoshida's construction from an ideal point of the deformation variety $\mathfrak{D}(M, \mathcal{T})$ which corresponds to an ideal point of the character variety. Then any essential surface obtained from $T$ by compressions is detected by the character variety.
\end{numless_thm}
This result uses a result of Tillmann (\cite{tillmann_degenerations}, Proposition 23) which gives the analogous statement for spun-normal surfaces to show that the same is true for twisted squares surfaces (which we refer to as Yoshida form surfaces here). Correspondence between an ideal point of the deformation and character varieties means that as we approach the ideal point in the deformation variety, the corresponding characters approach an ideal point of the character variety. We can see this by considering the trace of the meridian of the boundary torus:\\

As we saw in section \ref{hol_semi-merid}, the holonomy of the semi-meridian is of the form $\zeta^2 f$, where $f$ approaches some non-zero value as $\zeta \rightarrow 0$. The holonomy of the meridian is then $\zeta^4 f^2$. That is, going around the meridian acts on the sphere at infinity of $\mathbb{H}^3$ by $z \mapsto \zeta^4 f^2 z$. As an element of $PSL_2(\mathbb{C})$ we can write this as $\left(\begin{array}{cc} \zeta^2 f & 0 \\ 0 & (\zeta^{2} f)^{-1} \end{array} \right)$, and as $\zeta \rightarrow 0$ the trace of this goes to $\infty$.

\bibliographystyle{hamsplain}
\bibliography{/Users/h/Math/research/henrybib}

\providecommand{\bysame}{\leavevmode\hbox to3em{\hrulefill}\thinspace}
\providecommand{\href}[2]{#2}
\begin{thebibliography}{10}

\bibitem{bergman71}
G.M. Bergman, \emph{The logarithmic limit-set of an algebraic variety}, Trans.
  Amer. Math. Soc. \textbf{157} (1971), 459--469.

\bibitem{chesebrotillmann05}
Eric Chesebro and Stephan Tillmann, \emph{{Not all boundary slopes are strongly
  detected by the character variety}}, \mbox{arXiv:math.GT/0510418}.

\bibitem{cullershalen83}
Marc Culler and Peter Shalen, \emph{Varieties of group representations and
  splitting of 3-manifolds}, Annals of Mathematics \textbf{117} (1983),
  109--146.

\bibitem{floydhatcher82}
W.~Floyd and A.~Hatcher, \emph{Incompressible surfaces in punctured-torus
  bundles}, Topology and its Applications \textbf{13} (1982), 263--282.

\bibitem{gueritaud}
Fran\c{c}ois Gu\'eritaud and David Futer, \emph{{On canonical triangulations of
  once-punctured torus bundles and two-bridge link complements}},
  \mbox{arXiv:math.GT/0406242}.

\bibitem{kabaya07}
Yuichi Kabaya, \emph{{A method to find ideal points from ideal
  triangulations}}, \mbox{arXiv:math.GT/0706.0971}.

\bibitem{lackenby03}
Marc Lackenby, \emph{The canonical decomposition of once-punctured torus
  bundles}, Comment. Math. Helv. \textbf{78} (2003), no.~2, 363--384.

\bibitem{cullerjacorubinstein}
W.~Jaco M.~Culler and H.~Rubinstein, \emph{Incompressible surfaces in
  once-punctured torus bundles}, Proc. London Math. Soc. \textbf{45(3)} (1982),
  385--419.

\bibitem{mumford}
David Mumford, \emph{The red book of varieties and schemes}, second ed.,
  Springer, 1999.

\bibitem{ohtsuki94}
Tomotada Ohtsuki, \emph{Ideal points and incompressible surfaces in two-bridge
  knot complements}, J. Math. Soc. Japan \textbf{46} (1994), no.~1, 51--87.

\bibitem{schanuelzhang01}
S.~Schanuel and X.~Zhang, \emph{Detection of essential surfaces in 3-manifolds
  with $\textrm{SL}_2$ trees}, Annals of Mathematics \textbf{320(1)} (2001),
  149--165.

\bibitem{segerman_twsq_spnn}
Henry Segerman, \emph{On spun-normal and twisted squares surfaces},
  \mbox{arXiv:0810.1256v1}.

\bibitem{thurston}
W.~Thurston, \emph{Geometry and topology of 3-manifolds}.

\bibitem{tillmann_degenerations}
Stephan Tillmann, \emph{{Degenerations of ideal hyperbolic triangulations}},
  \mbox{arXiv:math.GT/0508295}.

\bibitem{yoshida91}
Tomoyoshi Yoshida, \emph{On ideal points of deformation curves of hyperbolic
  3-manifolds with one cusp}, Topology \textbf{30} (1991), no.~2, 155--170.

\end{thebibliography}

\end{document}